\documentclass{aims} %
\usepackage{amsmath}
\usepackage{paralist}
\usepackage[misc]{ifsym}
\usepackage{epsfig,tikz,pgfplots} 
\usepackage{epstopdf,enumitem} 
\usepackage[colorlinks=true]{hyperref}
\hypersetup{urlcolor=blue, citecolor=red}
\allowdisplaybreaks

\usepackage{booktabs,array,float}
\newcolumntype{M}[1]{>{\centering\arraybackslash}m{#1}}

\def\currentvolume{X}
 \def\currentissue{X}
  \def\currentyear{200X}
   \def\currentmonth{XX}
    \def\ppages{X-XX}
     \def\DOI{}
\pdfpagesattr
{ /CropBox [55 40 546 730] }

\newtheorem{theorem}{Theorem}[section]

\theoremstyle{definition}

\newtheorem{assumption}[theorem]{Assumption}

\newcommand{\eps}[1]{{#1}_{\varepsilon}}

\title[On the well-posedness of the Cahn--Hilliard--Biot model]
{On the well-posedness of the Cahn--Hilliard--Biot model and its applications to tumor growth} %

\author[Marvin Fritz]{}


\thanks{The author is supported by the state of Upper Austria.}

\thanks{$^*$Corresponding author: Marvin Fritz}

\usepackage[T1]{fontenc}

\usepackage{xcolor,bm, amssymb, mathtools,subfigure,comment}
\usepackage[capitalise,nameinlink]{cleveref}

\usepackage{tikz,bm}
\usepackage[sort,compress]{cite}

\crefname{subsection}{Subsection}{Subsections}
\crefformat{equation}{(#2#1#3)}
\crefformat{align}{(#2#1#3)}
\crefformat{enumi}{(#2#1#3)}
\crefname{figure}{Figure}{Figures}

\crefname{theorem}{Theorem}{Theorems}
\crefname{definition}{Definition}{Definitions}
\crefname{lemma}{Lemma}{Lemmas}
\crefname{assumption}{Assumption}{Assumptions}

	\definecolor{color00}{HTML}{4E79A7}
	\definecolor{color01}{HTML}{F28E2B}
	\definecolor{color02}{HTML}{E15759}

\makeatother

\newcommand{\pt}{\partial_t}
\renewcommand{\eps}{\varepsilon}

\newcommand{\T}{\mathcal{T}}

\newcommand{\dx}{\,\textup{d}x}

\renewcommand{\O}{\mathcal{O}}

\renewcommand{\div}{\textup{div}}

\newcommand{\dd}{\mathop{}\!\mathrm{d}}

\newcommand{\ddt}{\frac{\dd}{\dd\mathrm{t}}}
\newcommand{\dt}{\,\textup{d}t}
\newcommand{\p}{\partial}

\newcommand{\F}{\mathcal{F}} %
\newcommand{\E}{\mathcal{E}} %
\newcommand{\R}{\mathbb{R}} %

\newcommand{\red}{\textcolor{red}}

\newcommand{\N}{\mathbb{N}} %

\renewcommand{\rho}{\varrho}

\renewcommand{\L}{\mathcal{L}} %
\newcommand{\C}{\mathbb{C}} %

\newcommand{\OT}{{\Omega_T}}

\newcommand{\pf}{\varphi}

\renewcommand{\div}{\textup{div}}

\renewcommand{\phi}{\varphi}

\begin{document}
\maketitle

\centerline{\scshape
Marvin Fritz$^{{\href{mailto:marvin.fritz@ricam.oeaw.ac.at}{\textrm{\Letter}}}*1}$}

\medskip

{\footnotesize
 \centerline{$^1$Johann Radon Institute for Computational and Applied Mathematics,} \centerline{Austrian Academy of Sciences, Linz, Austria}
} %

\bigskip

\begin{abstract}
We study the Cahn--Hilliard--Biot model with respect to its mathematical well-posedness. The system models flow through deformable porous media in which the solid material has two phases with distinct material properties. The two phases of the porous material evolve according to a generalized Ginzburg--Landau energy functional, with additional influence from both viscoelastic and fluid effects. The flow-deformation coupling in the system is governed by Biot's theory. This results in a three-way coupled system that can be viewed as an extension of the Cahn--Larché equations by adding a fluid flowing through the medium. We distinguish the cases between a spatially dependent and a state-dependent Biot--Willis function. In the latter case, we consider a regularized system.  In both cases, we use a Galerkin approximation to discretize the system and derive suitable energy estimates. Moreover, we apply compactness methods to pass to the limit in the discretized system. In the case of Vegard’s law and homogeneous elasticity, we show that the weak solution depends continuously on the data and is unique. Lastly, we present some numerical simulations to highlight the features of the system as a tumor growth model.
\end{abstract}

\section{Introduction }

In this paper, we undertake a comprehensive well-posedness analysis of the Cahn--Hilliard--Biot model that is capable of capturing complex scenarios involving fluid flow through deformable porous media at the Darcy scale. This model has recently been derived in \cite{storvik2022cahn} and exhibits the ability to adapt to changing material characteristics, including variations in stiffness, permeability, compressibility, and poroelastic coupling strength. These variations arise from Cahn--Hilliard-type phase changes occurring within the solid matrix, which is a phenomenon with diverse applications.

One prominent application domain for this model is the study of solid tumor evolution. The stress effects resulting from tumor growth have been argued to exert a profound influence on the evolution of tumors themselves \cite{lima2016, lima2017}. Furthermore, stress has the potential to promote and inhibit tumor growth \cite{tumorstresscheng,tumorstresshelminger, tumorstressstylianopoulos}. Moreover, most solid malignant tumors exhibit elevated interstitial fluid pressure and alterations in the elastic properties of their surrounding matrix \cite{milosevictumor}. Conceptually, the porous medium can be viewed as a representation of the coexistence of healthy and malignant cells enveloped by an extracellular matrix, the fluid mimicking the interstitial fluid within this intricate biological environment. Several additional phenomena may be added to the model, such as chemotaxis \cite{fritz2019unsteady},  haptotaxis \cite{fritz2019localnonlocal} and angiogenesis \cite{fritz2021analysis}. We refer to \cite{fritz2023tumor} for an overview of tumor modeling using phase-field equations.

The proposed mathematical framework represents an extension that combines elements from the Cahn--Hilliard model and the Biot equations. In \cite{storvik2022cahn}, it was shown that the mathematical model possesses a generalized gradient flow structure, affirming the thermodynamic consistency of the model. The Cahn--Hilliard component governs solid phase changes within the system through a smooth phase--field variable, while the Biot equations oversee fluid flow and elasticity. The Cahn--Hilliard equation \cite{miranville2019cahn} employs interfacial free energy to model phase separation. The combination of the Cahn--Hilliard model with elasticity is often referred to as the Cahn--Larché model \cite{garcke2005mechanical}, and it has found applications in various areas, including lithium-ion batteries \cite{areias} and tumor evolution \cite{garcketumormechanics,fritz2020subdiffusive}. Other models have considered viscoelastic coupling; see \cite{grasselli2023phase,garcke2023approximation}. In the Cahn--Hilliard--Biot model, the fluid dynamics is incorporated into the system, assuming Biot-type coupling between flow and elasticity \cite{coussy}.

We prove the existence of weak solutions to the Cahn--Hilliard--Biot model based on a spatial discretization and deriving suitable energy estimates. We may conclude the existence of weakly/weakly-$\star$ convergent subsequences, which will turn out to be the weak solution to the original problem. In this step, we apply the Aubin--Lions compactness lemma to achieve strong convergences so that we can pass to the limit in the nonlinear terms of the model. Moreover, we prove, under several more restrictive assumptions, that the weak solution is well-posed, that is, it is unique and depends continuously on the data. In the analysis, we distinguish the cases between a spatially dependent and a state-dependent Biot--Willis function; we note that in most applications the Biot--Willis \textit{constant} is considered that is included in our analysis. In the latter case of a state-dependent function, we add a regularizing term to the system, which brings additional regularity that is required to pass to the limit in the proof. 

During the review period of our work, three preprints authored by independent research groups emerged, shedding further light on the analysis and solution strategies of the Cahn--Hilliard--Biot model (see \cite{abels2024existence,storvik2024sequential,riethmuller2023well,brunk2024}). In particular, the work by Abels et al. \cite{abels2024existence} explores additional challenges by considering different boundary conditions, which introduce complexities at specific points in their analysis. To tackle these complexities, they utilize a sequence of diverse techniques, incorporating concepts such as maximal regularity of solutions and a nonstandard application of the Aubin--Lions lemma. In addition, they regularize the volumetric fluid content and introduce a convolution in the deformation equation to enhance the robustness of their model. They investigate a limit process to return to the original model.
It should be noted that our previous investigation uncovered a consistent uniform energy estimate for the discretized solution, even without resorting to viscoelastic regularization. However, ensuring convergence in the elastic energy requires robust convergence of the deformation gradient. Consequently, the necessity for regularization techniques becomes apparent. This aspect had not been thoroughly addressed previously, prompting us to introduce regularization methods akin to those in subsequent works, albeit with a different approach to deriving strong convergence. 
Both works \cite{abels2024existence,riethmuller2023well} initially establish the strong convergence of the pressure before achieving the necessary convergence for the deformation gradient. In contrast to these studies, our approach tackles the problem in reverse order.

The structure of the article is as follows: In \cref{sec:pre}, we introduce the relevant function spaces and inequalities that appear in the well-posedness analysis. The Cahn--Hilliard--Biot model is presented in \cref{sec:model}. After introducing conservation laws for each of the three coupled processes (phase-field evolution, elasticity, fluid flow), we propose free energies and constitutive relations  to conclude the Cahn--Hilliard--Biot system.  Moreover, we discuss the case of a state-dependent Biot--Willis function and its consequences for the mathematical analysis of the system. In \cref{sec:well}, we state the definition and the existence of a weak solution for the Cahn--Hilliard--Biot system with a spatially dependent Biot--Willis function. This is the main result of this work, and we prove such a theorem in several steps. Moreover, we show the uniqueness and the continuous dependence of the weak solution under suitable assumptions. In \cref{sec:wellreg}, we prove the existence of a weak solution for a regularized Cahn--Hilliard--Biot system that allows the case of a state-dependent Biot--Willis function. Lastly, we show some numerical simulations in \cref{Sec:Numerics} and highlight the various effects of the elasticity and flow equations on the tumor's growth and shape.
\section{Notation and preliminaries} \label{sec:pre}

In the following, we assume that $\Omega \subset \mathbb{R}^d$, $d \in \{1,2,3\}$, is a bounded domain with a sufficiently smooth boundary $\partial\Omega$ and $T>0$ is a given fixed time horizon. The space-time cylinder is denoted by $\Omega_t=[0,t]\times \Omega$ for $t\in (0,T]$. Notationally, we write $(f,g)\mapsto(f,g)_\Omega:=\int_\Omega f(x)g(x) \dx$ for the integration of two functions $f\in L^p(\Omega)$, $g\in L^{q}(\Omega)$ with $\frac{1}{p}+\frac{1}{q}=1$ for $p,q \in [1,\infty]$. Moreover, we shall use standard notation for Lebesgue, Sobolev, and Bochner spaces. When denoting norms, we shall omit the spatial domain when no confusion is possible. 

We denote a generic constant simply by $C>0$ (which may change from line to line), and for brevity we may write $x \lesssim y$ instead of $x \leq Cy$.
We recall the H\"older, Poincar\'e--Wirtinger, Korn and Sobolev inequalities:
\begin{equation} \begin{aligned}
        \Vert  uv\Vert _{L^{\bar r}} & \leq \Vert  u\Vert _{L^{\bar p}} \Vert  v\Vert _{L^{\bar q}} && \forall u \in L^{\bar p}(\Omega), ~v \in L^{\bar q}(\Omega), \\ 
		\Vert  u-\langle u\rangle\Vert _{L^p} & \lesssim \Vert \nabla u\Vert _{L^p} && \forall u \in W^{1,p}(\Omega),   \\
  \Vert  u\Vert _{L^p} & \lesssim \Vert \nabla u\Vert _{L^p} && \forall u \in W_0^{1,p}(\Omega),   \\
		\Vert \nabla u\Vert _{L^p}^p              & \lesssim  \Vert \eps(u)\Vert _{L^p}^p   && \forall u \in W_0^{1,p}(\Omega), \\
		\Vert u\Vert _{W^{m,\tilde q}} & \lesssim \Vert u\Vert _{W^{k,\tilde p}} && \forall u\in W^{k,\tilde p}(\Omega),
	\end{aligned} \label{Eq:SobolevInequality} \end{equation}
 where the exponents satisfy the relationship $\frac{1}{\bar p}+\frac{1}{\bar q} = \frac{1}{\bar r}$, $\frac{1}{\hat p}+\frac{1}{\hat q} = 1 + \frac{1}{\hat r}$ and $k-\frac{d}{\tilde p}  \geq m-\frac{d}{\tilde q}$ for $k\geq m$, respectively.
Here, $\langle u\rangle =\frac{1}{\vert \Omega\vert } (u,1)_{L^2(\Omega)}$ denotes the mean of $u$ with respect to $\Omega$.

\section{Mathematical modeling}
\label{sec:model}

In this section, we briefly present the mathematical formulation of the Cahn--Hilliard--Biot model, which describes a saturated porous medium containing one fluid phase and two solid phases with distinct material properties. The solid phases are represented using a diffuse interface approach of Cahn--Hilliard type, where surface tension, solid material deformation, and pore pressure act as driving forces. We follow the original derivation in \cite{storvik2022cahn}.

\subsection{Cahn--Hilliard--Biot system} Let $T>0$ be a fixed final time and $\Omega\subset \mathbb{R}^d$, $d \in \{1,2,3\}$,  a bounded domain with boundary $\Gamma  := \partial \Omega$ that is of $C^{1,1}$-regularity or is convex. 
Furthermore, $\varphi: [0,T] \times \Omega \rightarrow [0,1]$ describes the phase-field representing the two phases $0$ and $1$. Specifically, we consider an application in tumor growth modeling by considering $\phi$ as the tumor volume fraction, that is, $\phi(x)=1$ means that all (100\%) of the cells at $x \in \Omega$ are cancerous. On the other hand, $\phi(x)=0$ corresponds to healthy cells. Furthermore, $u$ denotes the infinitesimal displacement with $|\nabla  u| \ll 1$, $\eps(u) =\frac12 (\nabla u + \nabla u^\top)$ the strain measure of $u$, $p$ the pore pressure and $q$ the fluid flux. The phase-field equation accounts for phase-change conservation through a phase-field flux $J_\phi$ and reactions $S_\phi$ possibly depending on the other variables:
\begin{equation} \label{Deriv:CH}
\partial_t\varphi + \div {J_\phi} = S_\phi.
\end{equation}

As the tumor expands, it exerts mechanical forces on the surrounding tissue, leading to deformation and organ remodeling. In addition, tumor cells actively remodel the extracellular matrix and exert mechanical forces on surrounding tissues, influencing their mechanical properties and deformation behavior.  We consider the elasticity equation
\begin{equation} 
    -\div  \sigma = S_u,
\end{equation}
where $\sigma$ is called stress tensor and $S_u$ models external body forces.  The usual assumptions of the Cahn--Larch\'e equations in \cite{garcke2005mechanical,garcketumormechanics} exclude external body forces. Moreover, these works and the original derivation \cite{storvik2022cahn} consider the quasi-static case under the assumption that mechanical equilibrium is rapidly achieved. Following the work \cite{bociu2023mathematical} that studied quasi-static Biot models with viscoelastic effects, we assume that the stress tensor $\sigma$ is split into elastic stress $\sigma_e$ and viscoelastic part $\sigma_v$. In line with the assumptions made in Abels et al.~\cite{abels2024existence}, we adopt the viscoelastic term in the form $\C_v \pt \eps(u)$ although we neglect the state dependency of the modulus of viscoelasticity. The addition of a viscoelastic term proves highly advantageous in ensuring the regularity of the deformation in our subsequent analysis. Specifically, when we test the variational form with $\pt u$, it leads to an estimate of $\pt \eps(u)$ in the space $L^2(\OT)$.

We consider a volume balance law for the fluid:
\begin{equation} \label{Deriv:Theta}
    \partial_t \theta + \div  J_\theta = S_\theta,
\end{equation}
where $\theta$ represents the volumetric fluid content, influenced by fluid flux $J_\theta$ and a source term $S_\theta$, which possibly couples $\theta$ to the tumor volume fraction $\phi$.

The model closure relies on its free energy, which encompasses three components: the surface energy, elastic energy, and fluid energy:
\begin{equation} \label{eq:totenergy}
    \mathcal{E}(\varphi, u, \theta) = \mathcal{E}_\phi(\varphi) + \mathcal{E}_u(\varphi,  u) +\mathcal{E}_\theta(\varphi, u, \theta).
\end{equation}
The surface energy is expressed as:
\begin{equation}
    \mathcal{E}_\phi(\varphi) := \int_\Omega \Psi(\varphi) + \frac{\gamma}{2}|\nabla\varphi|^2\dx,
\end{equation}
where $\Psi(\varphi)$ is a double-well potential that penalizes deviations from pure phases, and the second term accounts for interfacial energy. Here, $\gamma$ represents the interfacial tension.
The elastic energy, following the Cahn--Larché equations \cite{garcketumormechanics}, takes the form:
\begin{equation} \label{Def:EnergyU}
    \mathcal{E}_u(\varphi, u) = \int_\Omega W(\phi, \eps(u)) \dx= \int_\Omega \frac{1}{2} \big({ \varepsilon}( u) - \mathcal{T}(\varphi)\big)\!:\!\mathbb{C}(\varphi)\big({ \varepsilon}( u) - \mathcal{T}(\varphi)\big) \dx,
\end{equation}
where $\C(\phi)$ is the elasticity tensor and $\mathcal{T}(\varphi)$ is the eigenstrain (or: stress free strain) at $\varphi$, accounting for solid phase-specific properties.
Lastly, the fluid energy is expressed as:
\begin{equation} \label{Def:EnergyTheta}
    \mathcal{E}_\theta(\varphi, u, \theta) = \int_\Omega \frac{M(\varphi)}{2}\left(\theta - \alpha(\phi)\div  u\right)^2 \dx, %
\end{equation}
where $M$ is the compressibility function and $\alpha$ the  Biot--Willis function.
In classical Biot models, the Biot--Willis function is set to a constant value and is determined by the expression $\alpha = 1 - K / K_s$, where $K_s$ denotes the bulk modulus of the solid material and $K$ signifies the bulk modulus of the porous matrix \cite{showalter2000diffusion}. Therefore, it is usually referred to as the Biot--Willis constant in the literature. However, opting for a constant value significantly simplifies the analytical treatment. We will elaborate on the implications of assuming $\alpha$ to be spatially dependent (including the constant case) or dependent on the phase-field variable $\phi$. %

The source functions $S_\phi$ and $S_\theta$ appearing in \cref{Deriv:CH} and \cref{Deriv:Theta}, are short forms of $S_\beta=S_\beta(\phi,\eps(u),\theta)$ for $\beta \in  \{\phi,\theta\}$. When effect from the Biot model were not included, it was proposed in \cite{garcketumormechanics} that $$S_\phi=-\lambda_a k(\phi)+\lambda_p f(\phi) \eta/(1+|D_{\eps(u)} W|),$$ for $\lambda_a$ and $\lambda_p$ being the apoptosis and proliferation factors, respectively, $\eta$ the nutrients (which is not part of the model in this work), and $f$, $k$ bounded functions. In practice, it is common to employ the logistic growth function $f(\phi) = \phi(1 - \phi)$, often substituting $\phi$ with a cutoff operator, which ensures the boundedness of $f$. From a biological point of view, it is imperative that $\phi$ remains within the interval $[0,1]$, rendering the cutoff operator essentially identical to the identity operator.

The model relies on various constitutive relations by relating $J_\phi$, $J_u$, $J_\theta$ to the variational derivatives of $\E$ with respect to $\phi$, $u$ and $\theta$, respectively. For instance, Fick's law for non-ideal mixtures connects the flux $J_\phi$ to the negative gradient of the chemical potential, that is,
   $ J_\phi = -m(\varphi)\nabla \mu$
with $m(\varphi)$ representing the chemical mobility. The chemical potential $\mu$ is obtained as the variational derivative of the free energy with respect to $\varphi$, that is, $$\begin{aligned} \mu=\delta_\phi \E=\delta_\phi \E_\phi+\delta_\phi \E_u+\delta_\phi \E_\theta, \end{aligned}$$
with
\begin{equation}\label{Eq:EnergyDiff}\begin{aligned}\delta_\phi \E_\phi&=\Psi'(\phi)-\gamma\Delta \phi, \\
\delta_\phi \E_u&=\frac{1}{2}\left({ \varepsilon}( u) - \mathcal{T}(\varphi)\right)\!:\!\mathbb{C}'(\varphi)\left({ \varepsilon}( u)  - \mathcal{T}(\varphi)\right) - \mathcal{T}'(\varphi)\!:\!\mathbb{C}(\pf)\left({ \varepsilon}( u) - \mathcal{T}(\varphi)\right), \\
\delta_\phi \E_\theta &=\frac{M'(\varphi)}{2}(\theta-\alpha(\phi)\div  u)^2  - M(\pf)(\theta-\alpha(\phi)\div u)\alpha'(\phi)\div u.
    \end{aligned}\end{equation}
    In the case of a constant or spatially dependent Biot--Willis function, the last term in $\delta_\phi \E_\theta$ vanishes.
    
The elastic stress $\sigma_e$ and the pore pressure $p$ are defined based on the rate of change of energy with respect to strain and volumetric fluid content, respectively, that is, $$\begin{aligned}
p=\delta_{\theta} \E&=M(\phi)(\theta-\alpha(\phi) \div u), \\
\sigma_e=\delta_{\eps(u)} \E&=\mathbb{C}(\varphi)\left({ \varepsilon}( u)-\mathcal{T}(\varphi)\right)-  \alpha(\phi) M(\phi)(\theta-\alpha(\phi) \div u) I \\ &=\mathbb{C}(\varphi)\left({ \varepsilon}( u)-\mathcal{T}(\varphi)\right)-  \alpha(\phi) p I, 
    \end{aligned}$$
    where we inserted $p$ in the last term. We note that this approach has not been employed in other works, such as \cite{abels2024existence,riethmuller2023well}, on the analysis of the Cahn--Hilliard--Biot model. While we will need additional steps to establish an energy balance in the discrete setting later on, this method offers a more advantageous limit process.
Finally, Darcy's law governs fluid flow through the porous medium, relating the fluid flux $J_\theta$ to the gradient of pore pressure $p$, that is,
    $J_\theta = -\kappa(\varphi)\nabla p$
where $\kappa(\varphi)$ represents the permeability, accounting for the phase-field's influence on fluid flow. %

Combining the balance equations with the constitutive relations and making appropriate identifications, the resulting Cahn--Hilliard--Biot model comprises a system of partial differential equations, providing insights into complex phenomena within porous media:
\begin{equation} \label{Sys:CHB}\begin{aligned}
\partial_t \varphi &= \div (m(\pf) \nabla \mu) + S_\phi, \\
\mu &=\Psi'(\varphi) - \gamma\Delta \varphi +\delta_\varphi\mathcal{E}_u(\varphi,  u) + \delta_\varphi \mathcal{E}_\theta(\varphi,  u, \theta) ,\\
0 &=\div(\C_v \pt \eps(u))+ \div\big(\mathbb{C}(\varphi)\big({ \varepsilon}( u)-\mathcal{T}(\varphi)\big)\big)   -\nabla(\alpha(\phi) p) + S_u, \\
\partial_t \theta &= \div (\kappa(\phi)\nabla p) + S_\theta,\\
p&=M(\phi)(\theta-\alpha(\phi) \div u). %
\end{aligned}\end{equation}
This system is equipped with the initial conditions $\phi(0)=\phi_0$, $u(0)=u_0$, $\theta(0)=\theta_0$, and the boundary conditions:
$$\begin{aligned} \nabla \phi \cdot n = m(\phi)\nabla \mu \cdot n = u \cdot n =\kappa(\phi)\nabla p \cdot n &=0 &&\text{on } \partial \Omega.
\end{aligned}$$ 
As usual, we consider no-flux boundary conditions for $\phi$ and $\mu$ in the Cahn--Hilliard equation. Similarly, we assume a no-flux boundary condition for the pressure $p$. Following the work of Bociu et al. \cite{bociu2023mathematical} on the Biot model with viscoelastic effects, we postulate a zero Dirichlet condition for the deformation $u$. Analytically, this is the most accessible case that we take, as we want to put the focus on the model. One could consider a boundary that is split in two subregions with a Dirichlet boundary on one for the possible presence of a rigid part of the body such as a bone which prevents variations of the displacement. In the other subregion, one would postulate a homogeneous Neumann condition so that the normal component of the stress is equal to zero. This has been done in other Cahn--Hilliard systems that model deformation, such as in \cite{abels2024existence,garcketumormechanics}, and is also assumed in the analysis of the Cahn--Hilliard--Biot model by Abels et al. \cite{abels2024existence}.

Deriving an a priori estimate under the assumption of a sufficiently regular solution tuple is straight-forward using the structure of the generalized gradient flow and yields
		\begin{equation*}\begin{aligned}
				&\|\phi(t)\|_{H^1}^2 + \|\nabla \mu\|^2_{L^2_t(L^2)} +\|\Psi(\phi(t))\|_{L^1} + \|u(t)\|_{H^1}^2 +\|\partial_t u\|_{L^2_t(H^1)}^2  \\
				&+ \|p(t)\|_{L^2}^2+ \|\nabla p\|^2_{L^2_t(L^2)} 
    \lesssim 1+\|\phi_0\|_{H^1}^2 + \|\theta_0\|_{L^2}^2+\|u_0\|_{H^1}^2.
		\end{aligned}\end{equation*}
  The issue lies in the passage to the limit afterwards when employing a spatial discretization and deriving this inequality on a discrete level. We comment on this in the subsection.

\subsection{Discussion}

In this work, we analyze the Cahn--Hilliard--Biot system as stated in \cref{Sys:CHB} with a spatially dependent Biot--Willis function $x \mapsto \alpha(x)$. This includes the case of the typical Biot--Willis constant $\alpha$. We refer to the discussion below \cref{Def:EnergyTheta} that in most typical applications one may simply select $\alpha=1$.

However, we still want to elaborate on a state dependent Biot--Willis function $\phi \mapsto \alpha(\phi)$ as it was part of the original derivation of the Cahn--Hilliard--Biot system in Storvik et al. \cite{storvik2022cahn}. This makes the analysis significantly more difficult in addition to the obvious point that there is one less term to treat in the equation of the chemical potential; see the second term in $\cref{Eq:EnergyDiff}_3$ with $\alpha'(\phi)$. We intend to make some more comments on the analytical difficulties that arise with a phase-field-dependent Biot--Willis function. 

In other applications of Biot systems of elasticity, it is typically assumed that $M$ and $\alpha$ are of constant value. Then, one can take the time derivative of the pressure equation to get $$\pt p = M (\pt \theta - \alpha \div \pt u).$$
and insert $\pt \theta$ to obtain the following diffusion equation governing the pressure $p$
$$\pt p = M \div (\kappa(\phi) \nabla p) + M S_\theta - M \alpha \div \pt u,$$
see \cite{showalter2000diffusion,bociu2016analysis,bociu2023mathematical}.
We see that in state dependent functions $\alpha$ and $M$ this would invoke by the product rule quadratic nonlinearities that are much more difficult to analyze. This would require linearization techniques and employing a fixed-point approach.

However, we are still able to work things out in a mixed coupled $(\theta,p)$ system. The main difficulty comes from the deformation equation, where we have the term $\nabla(\alpha(\phi)p)$. As usually, we require strong convergence of $\eps(u)$ in $L^2(\OT)$ as we have a quadratic nonlinearity of $\eps(u)$ in $\delta_\phi \E_u$ in the equation for the chemical potential. One way to derive the strong convergence of $\eps(u)$ is to obtain a bound of its time derivative (which is doable due to the presence of the viscoelastic term) and a uniform bound of $\eps(u)$ in $L^2(0,T;W^{1,3/2}(\Omega))$. Then we could apply the Aubin--Lions compactness lemma to infer strong convergence. However, we do not get such high regularity without relaxing the assumptions on $\C$ and $\T$ like assuming Vergard's law. \red{just assume this?} Another way was used in \cite{garcketumormechanics} by testing the discretized deformation equation with $u_k-v_k$ where $v_k$ is a sequence that strongly converges to the limit of $u_k$, that is, $u$. Therefore, we get the term $$(\alpha(\phi_k) p_k,\div(u_k-v_k))_\Omega.$$ However, $\div(u_k-v_k)$ is only weakly converging, and we have at this point no real chance of obtaining strong convergence of $p_k$ as we have no control of its time derivative, and using other (easily accessible) techniques boils down to needing strong convergence of $\theta_k$ and/or $u_k$. However, we do not have some higher spatial regularity of $\theta_k$ to get a better strong convergence than in an inverse Sobolev space. One way is doing integration by parts of the above term to get terms of the form
$(\alpha'(\phi_k) \nabla \phi_k p_k, u_k-v_k)$ and $(\alpha(\phi_k) \nabla p_k, u_k-v_k)$. The second one is actually fine as $u_k-v_k$ is strongly converging, which complements the weak convergence of $\nabla p_k$. However, the first term raises issues since we have no strong convergence of $\nabla \phi_k$. This term is crucial for us and leads us to assume that $\alpha'(\phi_k)=0$ or having strong convergence of $\nabla \phi_k$ available, which we achieve by adding a regularization to the original Cahn--Hilliard--Biot system \cref{Sys:CHB} in the form of $\nu |\Delta \phi|^2$ in the energy $\E$. This leads to $H^2(\Omega)$-regularity in space of $\phi$ which we normally do not obtain easily as in the standard Cahn--Hilliard equation, where we can solve the chemical potential equation for $\Delta \phi$ and using elliptic regularity theory. Here, the $\mu$-equation is more convoluted with additional (quadratic) terms that are not easily in $L^2(\Omega)$.

\subsection{Regularized Cahn--Hilliard--Biot system}
As we discussed, we add the energy 
$$\begin{aligned} \E_\nu(\phi,\theta)&= \int_\Omega \frac{\nu}{2} |\Delta \phi|^2  \dx,
\end{aligned}$$
to the total energy $\E$, see \cref{eq:totenergy}. Taking the functional derivative with respect to $\phi$, we see that only the equation of $\mu$  changes. In fact, one studies the regularized system %
\begin{equation} \label{Sys:CHBreg}\begin{aligned}
\partial_t \varphi &= \div (m(\pf) \nabla \mu) + S_\phi, \\
\mu &=\Psi'(\varphi) - \gamma\Delta \varphi +\nu \Delta^2 \phi +\delta_\varphi\mathcal{E}_u(\varphi,  u) + \delta_\varphi \mathcal{E}_\theta(\varphi,  u, \theta) ,\\
0 &=\div(\C_v \pt \eps(u))+ \div\big(\mathbb{C}(\varphi)\big({ \varepsilon}( u)-\mathcal{T}(\varphi)\big)\big)   -\nabla(\alpha(\phi) p) + S_u, \\
\partial_t \theta &= \div (\kappa(\phi)\nabla p) + S_\theta,\\
p&=M(\phi)(\theta-\alpha(\phi) \div u). %
\end{aligned}\end{equation}
which is equipped with the same initial and boundary conditions as before, but we also demand $\nabla\!\Delta \phi \cdot n = 0$ on $\partial \Omega$ to justify the integration by parts.

\section{Well-posedness of the Cahn--Hilliard--Biot model} \label{sec:well}
In this section, we state and prove the main theorems of this work. We show that the Cahn--Hilliard--Biot system admits a weak solution, and under certain circumstances the solution is unique. As written before, we distinguish the cases whether the Biot--Willis parameter $\alpha$ is taken to be dependent on $x$ or on the phase-field variable $\phi$. In this section, we first study the case of a spatially dependent Biot--Willis function that includes the typical Biot--Willis constant $\alpha(x)=\alpha$. The other case is then studied in \cref{sec:wellreg}. %

Before we prove that a weak solution exists to the Cahn--Hilliard--Biot system, we consider some assumptions on the functions that appear in the model. 

\begin{assumption}\label{Assumption} Let the following assumptions hold: \vspace{.1cm}
	\begin{enumerate}[label=\textup{(A\arabic*)}, ref=A\arabic*, leftmargin=.9cm] \itemsep.2em
   \item  $\phi_0 \in H^1(\Omega)$, $u_0 \in H_0^1(\Omega)^d$, $\theta_0 \in L^2(\Omega)$,\label{Ass:Init}
  \item $\gamma>0$ is fixed.\label{Ass:Con}
  \item \label{Ass:Psi} $\Psi\in C^1(\R;\R_{\geq 0})$ satisfies the growth condition $|\Psi'(x)|\leq C_{\Psi}(1+|x|)$ for any $x \in \R$ with $C_\Psi>0$.
    \item \label{Ass:Fun} $S_\phi$, $S_\theta$ are continuous and bounded, $S_u \in L^2(\Omega)^d$.
    \item $W(x, \F)  = \frac{1}{2} (\F - \T(x)) : \C(x) (\F - \T(x))$ for any $x\in\R$, $\F \in \R^{d\times d}$
with $\C$ bounded, Lipschitz continuous, differentiable and $\T$ Lipschitz continuous, differentiable; further, $\C$ fulfills the usual symmetry
conditions of linear elasticity and that there exists a $C_\C >0$
such that $C_\C|\F|^2\leq \F : \C (x)\F$ for any $x\in\R$, $\F \in \R^{d\times d}_{\rm sym}$; $\C_v$ fulfills the same assumptions as $\C$. \label{Ass:Elastic}
  \item \label{Ass:Mob} $\kappa,m \in C^0(\R)$, $\alpha, M\in W^{1,\infty}(\R)$  s.t. $M(x) \geq M_0>0$, $\kappa(x) \geq \kappa_0>0$, $m(x)\geq m_0>0$ for any $x \in \R$.
\end{enumerate}
\end{assumption}

We note that $\phi$ realistically takes its values in the interval $[0,1]$. Therefore, \cref{Ass:Psi} allows for typical double-well potentials of polynomial or logarithmic type if they are extended outside the relevant interval appropriately. Alternatively to \cref{Ass:Psi}, one can assume growth conditions on $\Psi$ and its derivatives relating to other orders of derivatives; see \cite{garcketumormechanics}. Since $\Psi'$ is assumed to have linear growth, see \cref{Ass:Psi}, this implies that $\Psi$ is semiconvex, that is. 
\begin{equation} \label{Ass:PsiNew}(\Psi'(x)-\Psi'(y))(x-y)\geq -C_{*}|x-y|^2,
\end{equation} for any $x,y \in \R$ for some $C_*\geq 0$.
We note that \cref{Ass:Elastic} implies that it holds
\begin{equation} \label{Ass:Elastic2}
    \begin{aligned}
        (D_{\F_1} W(x,\F_1) - D_{\F_2} W(x,\F_2)) : (\F_1-\F_2) &\geq C |\F_1-\F_2|^2, \\
        |W(x,\F)|+|D_x W(x,\F)| &\leq C(1+|x|^2+|\F|^2), \\
        |D_{\F} W(x,\F) | &\leq C(1+|x|+|\F|).
    \end{aligned}
\end{equation}
for any $s \in \R$, $\F_1,\F_2 \in \R^{d \times d}$, $\F \in \R^{d \times d}_{\text{sym}}$. \smallskip

Then, the existence theorem  reads as follows.

\begin{theorem}[Existence of a weak solution] \label{Theorem}
Let Assumption \ref{Assumption} hold.
    Then there exists a global-in-time weak solution
    $$\begin{aligned}
    \phi &\in C_w([0,T];H^1(\Omega)) \cap C^0([0,T];L^r(\Omega)) \cap H^1(0,T;(H^1(\Omega))'), \\
    \mu &\in L^2(0,T;H^1(\Omega)), \\
    u &\in H^1(0,T;H_0^1(\Omega)^d), \\
    p &\in L^\infty(0,T;L^2(\Omega)) \cap L^2(0,T;H^1(\Omega)), \\
    \theta &\in C_w([0,T];L^2(\Omega)) \cap H^1([0,T];(H^1(\Omega))'),
\end{aligned}$$
with $r<6$ for $d=3$ or $r<\infty$ for $d\leq 2$,
to the Cahn--Hilliard--Biot system with the Biot--Willis function $x \mapsto \alpha(x)$, see \cref{Sys:CHB}, in the sense that the variational equations
    \begin{equation} \begin{aligned}
(\pt \phi,\zeta_\phi)_\Omega + (m(\phi) \nabla \mu,\nabla \zeta_\phi)_\Omega&=(S_\phi,\zeta_\phi)_\Omega, \\
 - (\mu,\zeta_\mu)_\Omega+\gamma  (\nabla\phi,\nabla \zeta_\mu)_\Omega + (\Psi'(\phi),\zeta_\mu)_\Omega &= (\delta_\phi \E_u + \delta_\phi \E_\theta,\zeta_\mu)_\Omega,    \\
(\C_v \p_t \eps(u),\eps(\zeta_u))_\Omega + \big(\mathbb{C}(\varphi)\left({ \varepsilon}( u)-\mathcal{T}(\varphi)\right),\eps(\zeta_u)\big)_\Omega &= (\alpha  p,\div(\zeta_u))_\Omega+(S_u,\zeta_u)_\Omega,  \\
(\pt \theta,\zeta_\theta)_\Omega + (\kappa(\phi)\nabla p,\nabla\zeta_\theta)_\Omega   &= (S_\theta,\zeta_\theta)_\Omega, \\
 (p,\zeta_p)_\Omega  &= (\theta-\alpha \div u,M(\phi)\zeta_p)_\Omega,
\end{aligned} \label{Sys:Variational} \end{equation}
are satisfied for any test functions $\zeta_\phi,\zeta_\mu,\zeta_\theta,\zeta_p  \in H^1(\Omega)$, $\zeta_u \in H_0^1(\Omega)^d$, and the initials $\phi(0)=\phi_0$, $u(0)=u_0$, $\theta(0)=\theta_0$. 
    Moreover, the solution fulfills for any $t \in (0,T]$ the energy inequality
\begin{equation} \label{Est:FinalCont}\begin{aligned}
&\|\phi(t)\|_{H^1}^2 + \|\mu\|^2_{L^2_t(H^1)} + \|p\|^2_{L^2_t(H^1)} + \|p(t)\|_{L^2}^2+\|\Psi(\phi(t))\|_{L^1} + \|u\|_{H^1(H^1_0)}^2   \\
&\lesssim 1+\|\phi_0\|_{H^1}^2 +\|u_0\|_{L^2}^2+ \|\theta_0\|_{L^2}^2.
\end{aligned}\end{equation}
\end{theorem}

Under additional assumptions, we are able to prove the uniqueness of the weak solution. Assuming an affine linear eigenstrain is quite standard in the literature, and this case is referred to as Vergard's law, see \cite{garcketumormechanics}. The appropriate testing and cancelation require constants $m$ and $\kappa$. Moreover, we assume that $M$ and $\alpha$ are constant, which drastically simplifies the model, since it implies that $\delta_\phi \E_\theta=0$, that is, it eliminates quadratic nonlinearities.

\begin{theorem}[Continuous dependence and uniqueness] \label{Thm:Uniqueness}
    Let \cref{Ass:Init}--\cref{Ass:Psi} of Assumption \ref{Assumption} hold. Additionally, we assume:
    \begin{enumerate}[start=4,label=\textup{(A\arabic*$^*$)}, ref=A\arabic*$^*$, leftmargin=1.05cm] \itemsep.2em
    \item \label{Ass:FunNew} In addition to \cref{Ass:Fun}, $S_\theta$ and $S_\phi$ are Lipschitz continuous. %
    \item In addition to \cref{Ass:Elastic}, $\C(\phi)=\C$ constant, $\T(\phi)$ affine linear, that is, $\T(\phi)=\T_1+\T_2\phi$ for $\T_1,\T_2 \in \R^{d \times d}$.\label{Ass:ElasticNew}
  \item \label{Ass:MobNew} $\kappa,m,\alpha,M>0$ are constant.
\end{enumerate}
Then the weak solution $(\phi,\mu,u,p,\theta)$ to the Cahn--Hilliard--Biot system is unique. Moreover, any two weak solutions $(\phi_i,\mu_i,u_i,p_i,\theta_i)$, $i \in \{1,2\}$, continuously depend on the data $(\phi_{0,i},\theta_{0,i},u_{0,i})$, $i \in \{1,2\}$, in the sense that it holds
		\begin{equation} \label{ContDep}
			\begin{aligned}
				&\|\phi_1-\phi_2\|_{L^\infty(H^1)'\cap L^2(L^2)}^2+\|\mu_1-\mu_2\|_{L^2(H^1)'}^2+\|\pt u_1-\pt u_2\|_{L^2(H^1_0)}^2\\ &\quad +\|u_1-u_2\|_{L^\infty(H^1_0)}^2 +  \|\theta_1-\theta_2\|_{L^2(L^2)}^2 +\|p_1-p_2\|_{L^\infty(L^2)\cap L^2(H^1)}^2\\ &\lesssim  \|\phi_{0,1}-\phi_{0,2} \|_{(H^1)'}^2+ \|u_{0,1} - u_{0,2} \|_{H^1}^2 +\|\theta_{0,1}-\theta_{0,2} \|_{L^2}^2 .
			\end{aligned}
		\end{equation}
\end{theorem}

To prove the existence theorem, we follow the Galerkin procedure by deriving energy estimates on a discrete level before returning to the continuous problem. In fact, the steps of the procedure read as follows:

\begin{enumerate} \itemsep0.5em
\item \textbf{Approximate problem}  We approximate the system by a problem in a finite-dimensional space. This reduces the problem to ODEs, and we can apply standard theory to ensure the existence of a solution $(\phi_k,\mu_k,u_k,\theta_k,p_k)$ to this finite-dimensional problem. 

	\item \textbf{Energy estimates.} In this step, we show that $(\phi_k,\mu_k,u_k,\theta_k,p_k)$ is uniformly bounded in norms of reflexive/separable Banach spaces. Thus, there is a subsequence (denoted by the same index) $(\phi_k,\mu_k,u_k,\theta_k,p_k)$ that converges weakly/weakly-$*$ to some limit element $(\phi,\mu,u,\theta,p)$.
	\item \textbf{Strong convergence.} We prove that the derivative of $\phi_{k}$ is bounded in another Bochner space and thus we can apply the Aubin--Lions lemma to conclude that $(\phi_{k})_k$ converges strongly in some space. This strong convergence is essential for the limit process later on. Moreover, we show the strong convergence of $(\theta_{k})_k$, $(u_{k})_k$, and $(p_{k})_k$. %
    \item \textbf{Energy inequality and initial conditions.} We show that the limit functions fulfill an energy inequality. Further, we show that the limit functions $\phi$, $u$ and $\theta$ also fulfill the imposed initial condition $\phi(0)= \phi_0$, $u(0)=u_0$ and $\theta(0)=\theta_0$ in some sense. This is performed using the strong convergence at $t=0$ and the uniqueness of limits. 
	\item \textbf{Limit process.} We are at the point where we have already proved the existence of functions $(\phi_k,\mu_k,u_k,\theta_k,p_k)$ that fulfill the $k$-th Galerkin equations, respectively. In this step, we take the limit $k \to \infty$ of the $k$-th Galerkin equations to obtain the variational Cahn--Hilliard equation. %
 This finishes the proof of \cref{Theorem}. %
 \item \textbf{Continuous dependence and uniqueness.} In \cref{Thm:Uniqueness}, we consider several simplifications, e.g., by assuming a constant mobility.  We prove the solution's continuous dependence on the data. From this we obtain the uniqueness of the solution.
\end{enumerate}

\begin{proof}[Proof of \cref{Theorem}] As described, we follow the Galerkin procedure to prove the existence of a weak solution. \smallskip
	
	\noindent\textbf{Step 1 (Approximate problem).}
	We choose $\{z_i\}_{i \in \N}$ as the set of eigenfunctions of the Neumann--Laplacian operator that is orthonormal in $L^2(\Omega)$ and orthogonal in $H^1(\Omega)$ with $z_1$ being the constant function $|\Omega|^{-1/2}$ and $(z_i,1)_\Omega = 0$ for $i \geq 2$.  In \cite[\S 3]{garcke2016global} it is shown that $\{z_i\}_{i \in \N}$ forms a basis in $H^2_n(\Omega)$,  the space of $H^2(\Omega)$ functions with vanishing normal component on $\partial \Omega$.
	Further, we choose $\{y_i\}_{i \in \N}$ as the eigenfunctions of for the vector-valued Dirichlet--Laplacian operator. These eigenfunctions are orthonormal in $L^2(\Omega)^d$ and orthogonal in $H_0^1(\Omega)^d$. We define the finite-dimensional spaces $Z_k$ and $Y_k$
	as the linear span of the first $k$ eigenfunctions $\{z_i\}_{i \in \N}$
	and $\{y_i\}_{i \in \N}$, respectively. Moreover, we denote by $\Pi_{Z_k}$ and $\Pi_{Y_k}$ the $L^2$-projection onto  $Z_k$ and $Y_k$, respectively.  Then, the Galerkin
	aapproximation of \cref{Sys:Variational} reads as: for any $k \in \N$
	find $(\phi_k, \mu_k, u_k,\theta_k,p_k)$ of the form
	$$\begin{aligned}
		\varphi_k(t,x) &= \sum_{i=1}^k a_{i}^k(t) z_i(x), \quad \mu_k(t,x) = \sum_{i=1}^k b_{i}^k(t) z_i(x),  \quad u_k(t,x) = \sum_{i=1}^k c_{i}^k(t) y_i(x), \\ \theta_k(t,x) &= \sum_{i=1}^k d_{i}^k(t) z_i(x), \quad p_k(t,x) = \sum_{i=1}^k e_{i}^k(t) z_i(x),
	\end{aligned}$$
	satisfying for a.e. $t \in (0,T)$  the system
	\begin{equation} \label{Sys:GalerkinVar} \begin{aligned}
			(\pt \phi_k,z)_\Omega + (m(\phi_k)\nabla \mu_k,\nabla z)_\Omega&=  (S_\phi^k,z)_\Omega, \\
			- (\mu_k,z)_\Omega+\gamma  (\nabla\phi_k,\nabla z)_\Omega + (\Psi'(\phi_k),z)_\Omega &= -(\delta_{\phi_k} \E_\theta^k + \delta_{\phi_k} \E_u^k,z)_\Omega,    \\
			(\C_v\pt \eps(u_k),\eps(y))_\Omega + \big(\mathbb{C}(\varphi_k)\left({ \varepsilon}( u_k)-\mathcal{T}(\varphi_k)\right),\eps( y)\big)_\Omega &= ( \alpha  p_k,\div y)_\Omega + (S_u,y)_\Omega, \\
			(\pt \theta_k,z)_\Omega + (\kappa(\phi_k)\nabla p_k,\nabla z)_\Omega  &= (S_\theta^k,z)_\Omega,  \\
 (p_k,z)_\Omega&= (M(\phi_k)(\theta_k-\alpha \div u_k),z)_\Omega, \end{aligned}\end{equation}
	for all $z \in Z_k$, $y \in Y_k$. Here, we have set $\E_u^k=\E_u(\phi_k,u_k)$ and in the same way for $S_\phi^k$, $S_\theta^k$ and $\E_\theta^k$ where $\E_u$ and $\E_\theta$ are defined as in \cref{Eq:EnergyDiff}. We equip the system with the initial conditions $\phi_k(0)=\phi_{k,0}:=\Pi_{Z_k} \phi_0$, $\theta_k(0)=\theta_{k,0}:=\Pi_{Z_k} \theta_0$ and $u_k(0)=u_{k,0}:=\Pi_{Y_k} u_0$. We solve \cref{Sys:GalerkinVar}$_2$ for $\mu_k$ and plug it into \cref{Sys:GalerkinVar}$_1$. Similarly, we put $p_k$ from \cref{Sys:GalerkinVar}$_5$ into \cref{Sys:GalerkinVar}$_4$.
	Then, the orthogonality of $\{z_i\}_{i \in \N}$ with respect to the $L^2(\Omega)$-inner product allows us to
express the Galerkin approximation as a system of ordinary differential equations
in the coefficient vectors $\bm{a} := (a_1^k, \dots, a_k^k)$,
	$\bm{c} := (c_1^k, \dots, c_k^k)$, and
	$\bm{d}: = (d_1^k, \dots, d_k^k)$. 
    Since it holds $\phi_{k}(0,x)=\Pi_{Z_k} \phi_0(x) = \sum_{i=1}^k (\phi_0,z_i)_\Omega z_i(x)$, we set the initial condition as $a_i^k(0)=(\phi_0,z_i)_\Omega$ for any $i \in \N$ and in the same way for the other variables. Since all involved functions $m$, $\alpha$, $M$, $\kappa$, $\Psi'$, $S_\phi$, $S_\theta$, $\C$, $\T$ are continuous with respect to their arguments, the differential-algebraic system contains only contributions that are continuous in $\bm{a}$, $\bm{c}$, and $\bm{d}$. 
		We invoke the Cauchy--Peano theorem to obtain the existence of $T_k \in (0,T]$ and local solutions $\bm{a}, \bm{c}, \bm{d} \in C^1([0,T_k];\R^k)$ solving the Galerkin system. From this we further obtain $\bm{b}$,  $\bm{e} \in C^0([0,T_k];\R^k)$. For an interested reader who is interested in more details on the structure of the ordinary differential equations, we refer to \cite[Lemma 3]{riethmuller2023well} and \cite[Lemma 18]{abels2024existence} where the existence of a discretized solution to very similar Cahn--Hillaird--Biot systems were studied in full detail. \medskip
		
		\noindent\textbf{Step 2 (Energy estimates).} In this step, we derive the important energy estimate that will guarantee the existence of weakly converging subsequences.  \smallskip
		
		\noindent\textbf{Step 2a.} First, we consider the test functions $\mu_k$ and $K\phi_k$ in \cref{Sys:GalerkinVar}$_1$, $K>0$ to be determined later on, and $\pt \phi_k$ in \cref{Sys:GalerkinVar}$_2$ to get
		$$\begin{aligned}
			(\pt \phi_k,\mu_k)_\Omega  + (m(\phi_k),|\nabla \mu_k|^2)_\Omega&= (S_\phi^k,\mu_k)_\Omega, \\
			K (\pt \phi_k,\phi_k)_\Omega +K(m(\phi_k)\nabla \mu_k,\nabla \phi_k)_\Omega &= K(S_\phi^k,\phi_k)_\Omega, \\
			- (\mu_k,\pt \phi_k)_\Omega+\gamma  (\nabla\phi_k,\pt\nabla \phi_k)_\Omega + (\Psi'(\phi_k),\pt \phi_k)_\Omega &= -(\delta_{\phi_k} \E_\theta^k + \delta_{\phi_k} \E_u^k,\pt \phi_k)_\Omega.
		\end{aligned}$$
		We add the three equations, which cancels the mixed term $(\pt \phi_k,\mu_k)_\Omega$ and we obtain
		\begin{equation} \label{Est:CH1}\begin{aligned}
				&K(\pt \phi_k,\phi_k)_\Omega+(m(\phi_k),|\nabla \mu_k|^2)_\Omega+\gamma  (\nabla\phi_k,\pt\nabla \phi_k)_\Omega \\
				&\quad + (\Psi'(\phi_k),\pt \phi_k)_\Omega+ (\delta_{\phi_k} \E_u^k +\delta_{\phi_k} \E_\theta^k,\pt \phi_k)_\Omega  \\ &=(S_\phi^k,\mu_k+K\phi_k)_\Omega - K(m(\phi_k)\nabla\mu_k,\nabla\phi_k)_\Omega.
		\end{aligned}\end{equation}
		The left-hand sides of \cref{Est:CH1} can be rewritten using the chain rule and the term involving the mobility function $m$ can be estimated from below by the lower bound of $m$, see \cref{Ass:Mob}. The second term on right-hand side of \cref{Est:CH1} can be estimated by the upper bound of the mobility $m$ (which we denote by $m_\infty$ and in the same manner for the other upper bounds) and the Young inequality as follows:
		\begin{equation} 
  \label{Est1:RHS1} -K(m(\phi_k)\nabla\mu_k,\nabla\phi_k)_\Omega \leq \frac{m_0}{4}\|\nabla\mu_k\|_{L^2}^2 + \frac{m_\infty^2K^2}{m_0}\|\nabla \phi_k\|_{L^2}^2. \end{equation}
		We note that we can absorb the term involving $\|\nabla \mu_k\|_{L^2}^2$ by the left-hand side of \cref{Est:CH1}.
		Further, the source function $S_\phi$ is bounded by \cref{Ass:Fun} and for one of the two terms, we simply have by the Young inequality
	\begin{equation} 
  \label{Est1:RHS2}(S_\phi^k,K\phi_k)_\Omega \leq CK + C \|\phi_k\|_{L^2}^2.\end{equation}
		
  For the other term, we add and subtract the mean $\langle \mu_k \rangle$ and then apply the Poincar\'e and Young inequality, that is,
		\begin{equation} 
  \label{Est1:RHS3}\begin{aligned} (S_\phi^k,\mu_k)_\Omega  &= (S_\phi^k,\mu_k- \langle \mu_k \rangle)_\Omega+\langle \mu_k \rangle (S_\phi^k,1)_\Omega  \\ &\leq C\|\mu_k-\langle \mu_k \rangle\|_{L^1} + C |\langle \mu_k \rangle| 
			\\&\leq C+ \frac{m_0}{4} \|\nabla \mu_k\|_{L^2}^2 + C |\langle \mu_k \rangle|.
		\end{aligned}
		\end{equation}
		We observe that we still need to estimate the mean of $\mu_k$ on the right-hand side of the estimate. We test \cref{Sys:GalerkinVar}$_2$ by $z=1 \in Z_k$ to obtain
		\begin{equation} \label{Eq:MeanMu} \langle \mu_k \rangle = \int_\Omega \Big[ \Psi'(\phi_k) + D_{\phi_k} W(\phi_k,\eps(u_k)) + \delta_{\phi_k} \E_\theta^k \Big] \dx
  \end{equation}
        We can estimate the first term $\Psi'(\phi_k)$ by $C_\Psi(1+|\phi_k|)$ as  assumed in \cref{Ass:Psi}. The second term can be estimated using the growth estimate \cref{Ass:Elastic2}. That is,
        $$\int_\Omega D_{\phi_k} W(\phi_k,\eps(u_k)) \lesssim 1+\|\phi_k\|_{L^2}^2+\|\eps(u_k)\|_{L^2}^2.$$
		Regarding the last term, we use the upper bounds of $M'$ of $\alpha$, see \cref{Ass:Mob}, to get
		$$\begin{aligned}\int_\Omega \delta_{\phi_k} \E_\theta^k \dx &=\int_\Omega  \frac{M'(\phi_k)}{2} (\theta_k-\alpha \div u_k)^2 \dx \lesssim \|\theta_k\|_{L^2}^2 + \|\div u_k\|_{L^2}^2.
		\end{aligned}$$
		We can eliminate the $\theta_k$-norm on the right-hand side by relating it to $p_k$ and $\div u_k$. In fact, considering the test function $\theta_k$ in \eqref{Sys:GalerkinVar}$_5$  gives
		\begin{equation*}  M_0 \|\theta_k\|_{L^2}^2 \leq  \|p_k\|_{L^2} \|\theta_k\|_{L^2} + M_\infty \alpha_\infty \|\div u_k\|_{L^2} \|\theta_k\|_{L^2}. 
		\end{equation*}
			Thus, the Young inequality yields
			\begin{equation}\label{Est:Theta}\|\theta_k\|_{L^2}^2 \lesssim  \|p_k\|_{L^2}^2+ \|\div u_k\|_{L^2}^2.
			\end{equation}
		We put these estimates back into \cref{Eq:MeanMu} to get the following bound of the mean of $\mu_k$
		\begin{equation}\label{Est:MuMean}\begin{aligned} \langle \mu_k \rangle 
				\lesssim 1+ \|\phi_k\|_{L^2}^2 + \|u_k\|_{H^1}^2 + \|p_k\|_{L^2}^2.
		\end{aligned}\end{equation}
		Plugging the estimates \cref{Est1:RHS1}--\cref{Est1:RHS3}, \cref{Est:MuMean} back into the right-hand side of the inequality \cref{Est:CH1}, we obtain the following:
		\begin{equation}\label{Est:CH2}\begin{aligned}
				&\ddt \bigg[ \frac{K}{2}  \|\phi_k\|_{L^2}^2 +\frac{\gamma}{2}  \|\nabla \phi_k\|_{L^2}^2 +\|\Psi(\phi_k)\|_{L^1} \bigg]    
				+\frac{m_0}{2} \|\nabla \mu_k\|^2_{L^2} \\ &\quad +(\delta_{\phi_k} \E_u^k +\delta_{\phi_k} \E_\theta^k,\pt \phi_k)_\Omega \\ &\lesssim 1+\|\phi_k\|_{H^1}^2  + \|u_k\|_{H^1}^2 + \|p_k\|_{L^2}^2.
		\end{aligned}\end{equation}
		We observe that we require bounds on $u_k$ and $p_k$ to complete the energy estimate. \medskip
		
		\noindent\textbf{Step 2b.} Secondly, we consider the test function $\pt u_k$ in \cref{Sys:GalerkinVar}$_3$,  $p_k$ in \cref{Sys:GalerkinVar}$_4$,  and $\pt \theta_k$ in \cref{Sys:GalerkinVar}$_5$ to obtain the following:
		\begin{equation}\begin{aligned}
			(\C_v \pt \eps(u_k),\pt\eps(u_k))_\Omega + (\delta_{\eps(u_k)} \E_u^k,\pt \eps(u_k) )_\Omega- (\alpha p_k I,\pt \eps(u_k))_\Omega &= (S_u,\pt u_k)_\Omega\\ 
   (\pt \theta_k,p_k)_\Omega  + (\kappa(\phi_k),|\nabla p_k|^2)_\Omega  &= (S_\theta^k,p_k)_\Omega, \\
			(M(\phi_k)(\theta_k-\alpha\div u_k), \pt \theta_k)_\Omega&=(p_k,\pt \theta_k)_\Omega. 
		\end{aligned}\label{Eq:2b}\end{equation}
  We use the assumption of $\C_v$, see \cref{Ass:Elastic}, to derive a lower estimate of the first term in the first equation in the form of $C\|\pt \eps(u_k)\|_{L^2}^2$ for some constant $C>0$. In the first equation, we note that we cannot deduce that $-(\alpha p_k I, \pt \eps(u_k))_\Omega$ is the same as $(\delta_{\eps(u_k)} \E_\theta^k,\pt \eps(u_k))_\Omega$ since $p_k$ is not the same as $M(\phi_k)(\theta_k-\alpha \div u_k)$. However, we add and subtract the latter term to get on the left-hand side of the first equation
  $$\begin{aligned}&-(\alpha p_k I,\pt \eps(u_k))_\Omega + (\delta_{\eps(u_k)} \E_\theta^k,\pt \eps(u_k))_\Omega + (\alpha M(\phi_k)(\theta_k-\alpha \div u_k) I,\pt \eps(u_k))_\Omega \\
  &=(\delta_{\eps(u_k)} \E_\theta^k,\pt \eps(u_k))_\Omega + (\alpha (\text{Id}-\Pi_{Z_k})(M(\phi_k)(\theta_k-\alpha \div u_k) I),\pt \eps(u_k))_\Omega.
  \end{aligned}$$
  When estimating the second inner product on the right-hand side, we use a combination of the H\"older and Young inequalities and give $\|\pt \eps(u_k)\|_{L^2}^2$ a sufficiently small prefactor, which allows us to absorb this term. Moreover, we estimate the resulting term $\|\theta_k\|_{L^2}$ on the right-hand side using \cref{Est:Theta}. Thus, the first equation becomes the following
  		$$\begin{aligned}
			&C \|\pt \eps(u_k)\|_{L^2}^2+(\delta_{\eps(u_k)} (\E_u^k + \E_\theta^k),\pt \eps(u_k) )_\Omega \lesssim 1+ \|u_k\|_{H^1}^2 +\|\pt u_k\|_{L^2}^2+ \|p_k\|_{L^2}^2.
		\end{aligned}$$
  When adding this inequality to \cref{Eq:2b}$_2$ and \cref{Eq:2b}$_3$, we notice that the mixed term $(\pt \theta_k,p_k)_\Omega$ cancels. Then, we add the result to the previous estimate \cref{Est:CH2} to obtain
		$$\begin{aligned}
			&\ddt \bigg[ \frac{K}{2}  \|\phi_k\|_{L^2}^2 +\frac{\gamma}{2}  \|\nabla \phi_k\|_{L^2}^2 +\|\Psi(\phi_k)\|_{L^1} +\E_u^k+\E_\theta^k \bigg] \\ & \quad +C \|\pt \eps(u_k)\|_{L^2}^2+\frac{m_0}{2} \|\nabla \mu_k\|^2_{L^2} + \kappa_0 \|\nabla p_k\|_{L^2}^2
			\\
   &\lesssim 1+\|\phi_k\|_{H^1}^2  + \|u_k\|_{H^1}^2 + \|\pt u_k\|_{L^2}^2+ \|p_k\|_{L^2}^2,
		\end{aligned}$$
		where we already estimated $\kappa$ by $\kappa_0$ from below according to \cref{Ass:Mob} and $S_\theta^k$ by its upper bound according to \cref{Ass:Fun}. We integrate the estimate over the time interval $(0,t)$ for $t\leq T_k$ to get
		\begin{equation} \label{Est:Integrated}\begin{aligned}
				&\frac{K}{2} \|\phi_k(t)\|_{L^2}^2+\frac{\gamma}{2} \|\nabla\phi_k(t)\|_{L^2}^2 + \|\Psi(\phi_k(t))\|_{L^1}+ \E_u(\phi_k(t),u_k(t)) \\ &\quad + \E_\theta(\phi_k(t),u_k(t),\theta_k(t))+C \|\pt \eps(u_k)\|_{L^2_t(L^2)}^2 +\frac{m_0}{2} \|\nabla \mu_k\|^2_{L^2_t(L^2)}    + \kappa_0 \|\nabla p_k\|_{L^2_t(L^2)}^2 
				\\ &\lesssim 1+\|\phi_{k,0}\|_{H^1}^2+\|\Psi(\phi_{k,0})\|_{L^1}+\E_u(\phi_{k,0},u_{k,0})+ \E_\theta(\phi_{k,0},u_{k,0},\theta_{k,0})  \\&\quad +\|\phi_k\|_{L^2_t(H^1)}^2  + \|u_k\|_{L^2_t(H^1)}^2 + \|\pt u_k\|_{L^2_t(L^2)}^2+ \|p_k\|_{L^2_t(L^2)}^2.
			\end{aligned}
		\end{equation}
		It remains to treat the energies $\E_u$ and $\E_\theta$ in the inequality. This is done in the next substep before we return to this energy estimate. \medskip 
		
		\noindent{\bf Step 2c.} First, we study the elastic energy $\E_u$ that appears on both sides in \cref{Est:Integrated}. Using \cref{Ass:Elastic} and \cref{Ass:Elastic2}, that is, the strict monotonicity of $D_{\eps(u_k)} W$ with respect to its second argument, we find
		$$\begin{aligned} W(s,\F) &=W(s,0)+\int_0^1 D_2 W(s,t\F):t\F\frac{1}{t} \dt \geq C|\F|^2-C(1+|s|^2),
		\end{aligned}$$
		for any $s \in \R$ and $\F \in \R^{d\times d}$. Thus,
		it gives 
		\begin{equation} \label{Aux:1} \E_u(\phi_k,u_k)=\int_\Omega W(\phi_k,\eps(u_k)) \dx \geq C \|\eps(u_k)\|^2_{L^2} - C(1+\|\phi_k\|^2_{L^2}).
		\end{equation}
		On the other hand, the initial energy $\E_u(\phi_{k,0},u_{k,0})$ on the right-hand side of \cref{Est:Integrated} can be estimated from above by using the upper bound of $W$, see \cref{Ass:Elastic2}$_2$, as follows
		\begin{equation} \label{Aux:2} \E_u(\phi_{k,0},u_{k,0})  \leq C(1+\|\phi_{k,0}\|_{L^2}^2 + \|\eps(u_{k,0})\|_{L^2}^2).%
		\end{equation}

		Next, we investigate fluid energy on both sides of the estimate \cref{Est:Integrated}. First, we use the definition of $\E_\theta$ to obtain
		\begin{equation} \label{Aux:3} \begin{aligned} \E_\theta(\phi_k,u_k,\theta_k) &= \int_\Omega \frac{M(\phi_k)}{2}(\theta_k-\alpha \div u_k)^2  \dx. 
  \end{aligned}\end{equation}
  Further, by testing \cref{Sys:GalerkinVar}$_5$ with $p_k$ and using Young's inequality, we have
  $$\|p_k\|_{L^2}^2 = (M(\phi_k) (\theta_k-\alpha \div u_k),p_k)_\Omega \leq \frac12\|M(\phi_k) (\theta_k-\alpha \div u_k)\|_{L^2}^2 +\frac12 \|p_k\|_{L^2}^2.$$
  We split $M(\phi_k)$ into $M^{1/2}(\phi_k) M^{1/2}(\phi_k)$ and bound one of them with $M_\infty^{1/2}$. Then we obtain
  $$\frac12 \|p_k\|_{L^2}^2 \leq \frac{M^{1/2}_\infty}{2} \|M^{1/2}(\phi_k) (\theta_k-\alpha \div u_k)\|_{L^2}^2.$$
  Using the definition of $\E_\theta$, see \cref{Aux:3}, we see that
  $$\E_\theta(\phi_k,u_k,\theta_k) \geq \frac{1}{2M_\infty^{1/2}} \|p_k\|_{L^2}^2.$$
		Moreover, we can bound the initial fluid energy as follows:
		\begin{equation} \label{Aux:4}\begin{aligned}\E_\theta(\phi_{k,0},u_{k,0},\theta_{k,0}) &= \frac12 (M(\phi_{k,0}), (\theta_{k,0} - \alpha \div u_{k,0} )^2)_\Omega \\ &\leq M_\infty \big( \|\theta_{k,0}\|^2_{L^2} +\alpha_\infty^2 \|\div u_{k,0}\|_{L^2}^2\big).\end{aligned}\end{equation} \smallskip
		
		\noindent\textbf{Step 2d.} Now, we are in a position to return to  the integrated estimate \cref{Est:Integrated}. We insert the upper and lower bounds of the elastic and fluid energy, see \cref{Aux:1}--\cref{Aux:4}, back into  \cref{Est:Integrated} to obtain
		$$\begin{aligned}
			&(\tfrac{K}{2}-C) \|\phi_k(t)\|_{L^2}^2+\frac{\gamma}{2} \|\nabla\phi_k(t)\|^2_{L^2} +\|\Psi(\phi_k(t))\|_{L^1}+C \|\partial_t u_k\|_{L^2_t(H^1)}^2    \\ &\quad+ C\|u_k(t)\|_{H^1}^2+ \frac{m_0}{2}\|\nabla \mu_k\|^2_{L^2_t(L^2)}+ \kappa_0 \|\nabla p_k\|^2_{L^2_t(L^2)}  + \frac{1}{2M_\infty^{1/2}}\|p_k(t)\|_{L^2}^2  \\ 
			&\lesssim 1+\|\phi_{k,0}\|_{H^1}^2+ \|\theta_{k,0}\|_{L^2}^2+ \|\Psi(\phi_{k,0})\|_{L^1}+\|u_{k,0}\|_{H^1}^2 \\ &\quad +\|\phi_k\|_{L^2_t(H^1)}^2 + \|u_k\|_{L^2_t(H^1)}^2 + \|p_k\|_{L^2_t(L^2)}^2, 
		\end{aligned}$$
		where we used Korn's inequality on the left-hand side of the inequality to estimate $\|\eps(u_k)\|_{L^2}$ and $\|\eps(\pt u_k)\|_{L^2}$ by their full $H^1(\Omega)^d$-norm.
		At this point, we choose $K$ sufficiently large to get the prefactor in front of $\|\phi_k(t)\|_{L^2}^2$ positive. We can estimate the initials on the right-hand side, using the properties of the orthogonal projection, by $\|\phi_{k,0}\|_{H^1}^2 \leq \|\phi_0\|_{H^1}^2$ and similarly for $\|u_{k,0}\|_{H^1}^2$ and $\|\theta_{k,0}\|_{L^2}^2$. Moreover, we can integrate the growth condition in \cref{Ass:Psi} to get $$\|\Psi(\phi_{k,0})\|_{L^1} \leq C_\Psi (1+\|\phi_{k,0}\|_{L^2}^2) \leq C_\Psi (1+\|\phi_{0}\|_{L^2}^2).$$ 
		Then, by an application of the Gronwall inequality, we get
		\begin{equation} \label{Est:Final}\begin{aligned}
				&\|\phi_k(t)\|_{H^1}^2 + \|\nabla \mu_k\|^2_{L^2_t(L^2)} +\|\Psi(\phi_k(t))\|_{L^1} + \|u_k(t)\|_{H^1}^2 +\|\partial_t u_k\|_{L^2_t(H^1)}^2  \\
				&+ \|p_k(t)\|_{L^2}^2+ \|\nabla p_k\|^2_{L^2_t(L^2)} %
    \lesssim 1+\|\phi_0\|_{H^1}^2 + \|\theta_0\|_{L^2}^2+\|u_0\|_{H^1}^2.
		\end{aligned}\end{equation}
		Since the right-hand side is uniform in $k$, we can argue by a no-blow-up criterion to extend the existence interval by setting $T_k=T$ for any $k$.
		Moreover, we already proved that $\langle \mu_k \rangle(t)$ is bounded in $L^\infty(0,T)$, see \cref{Est:MuMean}, and thus we obtain a $k$-uniform bound of $\mu_k$ in the $L^2(0,T;H^1(\Omega))$-norm.   Moreover, we see that $\theta_k$ is uniformly bounded in the $L^\infty(0,T;L^2(\Omega))$-norm by \cref{Est:Theta} and the derived bounds of $p_k$ and $u_k$. 
		By the energy estimate \cref{Est:Final} and the Eberlein--Smulian theorem, we can already extract weakly converging subsequences (which we denote by the same index by a typical abuse of notation).  In fact, we get the existence of limit functions $(\phi,\mu,u,\theta,p)$ such that
		\begin{equation} \label{Convergence:Weak}\begin{aligned}
				\phi_k &\to \phi_k &&\text{weakly* in } L^\infty(0,T;H^1(\Omega)), \\
				\mu_k &\to \mu_k &&\text{weakly in } L^2(0,T;H^1(\Omega)), \\
				u_k &\to u &&\text{weakly* in } L^\infty(0,T;H^1_0(\Omega)^d), \\
                \pt u_k &\to \pt u &&\text{weakly* in } L^2(0,T;H^1_0(\Omega)^d), \\
				p_k &\to p &&\text{weakly* in } L^\infty(0,T;L^2(\Omega)) \cap L^2(0,T;H^1(\Omega)), \\
				\theta_k &\to \theta &&\text{weakly* in } L^\infty(0,T;L^2(\Omega)).
		\end{aligned}\end{equation} 
		\medskip

		\noindent\textbf{Step 3 (Strong convergence).} Since the system is nonlinear, we require strong convergence. To do so, we want to apply the Aubin--Lions compactness lemma, see \cite[Section 7.3]{roubivcek2013nonlinear}, which still requires a uniform bound of a time derivative. We consider an arbitrary element $\zeta \in L^2(0,T;H^1(\Omega))$. Then we test \cref{Sys:GalerkinVar}$_1$ with $\Pi_{Z_k} \zeta(t) \in Z_k$ which gives
		$$
		\begin{aligned} \langle \pt \phi_k, \zeta \rangle_{L^2(H^1)}&=\langle \pt \phi_k, \Pi_{Z_k} \zeta  \rangle_{L^2(H^1)} \\ &= -(m(\phi_k)\nabla \mu_k,\nabla \Pi_{Z_k}\zeta)_{\Omega_T} + (S_\phi^k,\Pi_{Z_k}\zeta)_{\Omega_T} \\
			& \leq m_\infty \|\nabla \mu_k\|_{L^2(L^2)} \|\Pi_{Z_k}\zeta\|_{L^2(H^1)} + C \|\Pi_{Z_k}\zeta\|_{L^2(L^2)} \\
			&\leq C \|\zeta\|_{L^2(H^1)},
		\end{aligned}$$
		where we used the uniform bound of $\nabla \mu_k$ as shown in \cref{Est:Final}. Since $\zeta$ was arbitrarily chosen, we obtain the uniform bound of $\pt \phi_k$ in $L^2(0,T;(H^1(\Omega))')$. Inferring the Aubin--Lions compactness lemma, we obtain the compact embedding
		$$L^\infty(0,T;H^1(\Omega)) \cap H^1(0,T;(H^1(\Omega))') \hookrightarrow \hookrightarrow C^0([0,T];L^r(\Omega)),$$
		where $r<6$ for $d=3$ and $r<\infty$ for $d\leq 2$ to ensure $H^1(\Omega) \hookrightarrow \hookrightarrow L^r(\Omega)$.
		Thus, we obtain the convergences
		\begin{equation}
			\label{Convergence:Strong}
			\begin{aligned}\phi_k &\to \phi &&\text{ strongly in } C^0([0,T];L^r(\Omega)), \\
				\pt \phi_k &\to \pt \phi &&\text{ weakly in } L^2(0,T;(H^1(\Omega))'),
		\end{aligned}\end{equation}

		We note that $W^{1,r'}(\Omega)$ is compactly embedded in $L^2(\Omega)$ for $r'$ being the H\"older conjugate of $r$, that is, it holds $r'>6/5$ for $d=3$ and $r'>1$ for $d=2$. Thus, we find %
		\begin{equation}
			\label{Convergence:Strong2}
			\begin{aligned}\theta_k &\to \theta &&\text{ strongly in } C^0([0,T];(W^{1,r'}(\Omega))'), \\
				\pt \theta_k &\to \pt \theta &&\text{ weakly in } L^2(0,T;(H^1(\Omega))').
		\end{aligned}\end{equation}
		Here we used that testing \cref{Sys:GalerkinVar}$_4$ by $\Pi_{Z_k} \zeta(t) \in Z_k$ yields in a straightforward manner
		$$\begin{aligned} \langle \pt \theta_k, \zeta \rangle_{L^2(H^1)} &=\langle \pt \theta_k, \Pi_{Z_k}\zeta \rangle_{L^2(H^1)} \\ &= -(\kappa(\phi_k)\nabla p_k,\nabla \Pi_{Z_k}\zeta)_{\Omega_T} + (S_\theta^k,\Pi_{Z_k}\zeta)_{\Omega_T} \\
			&\leq C \|\zeta\|_{L^2(H^1)}.
		\end{aligned}$$ 
	Specifically, as it holds $4/3 > 6/5$, $\theta_k$ strongly convergence in the Banach space $C^0([0,T];(W^{1,4/3}(\Omega))')$.

  Next, we deduce the strong convergence of
$u_k$ to $u$ in $L^2(0,T;H^1_0(\Omega)^d)$, which is required to pass the limit in the quadratic nonlinearity (in terms of $\eps(u_k)$) in $\delta_{\phi_k} \E_u^k$ appearing in the equation of the chemical potential $\mu_k$. First off, we apply the Aubin--Lions compactness lemma to infer the comapct emebdding
$$L^\infty(0,T;H_0^1(\Omega)^d) \cap H^1(0,T;H_0^1(\Omega)^d) \hookrightarrow \hookrightarrow C^0([0,T];L^r(\Omega)^d),$$
 with $r$ as above, and thus we obtain the following strong convergence:
\begin{equation}
		\label{Convergence:Strong4}
			\begin{aligned}u_k &\to u &&\text{ strongly in } C^0([0,T];L^r(\Omega)^d).
		\end{aligned}\end{equation}
Since 
$\cup_k Y_k$ is dense in $H^1_0(\Omega)^d$, we can choose a sequence $\{
v_k\}_{k \in \N}$ with $v_k(t) \in Y_k$ for a.e. $t \in (0,T)$ and 
\begin{equation} \label{Conv:VK}\begin{aligned}v_k &\to u &&\text{ strongly in } H^1(0,T;H^1_0(\Omega)^d), 
\end{aligned}\end{equation}
Thus, by the convergence properties of $u_k$, see \cref{Convergence:Weak} and \cref{Convergence:Strong4}, it yields 
\begin{equation} \label{Conv:UKVK}\begin{aligned} u_k - v_k &\to 0 &&\text{ weakly* in } L^\infty(0,T;H^1_0(\Omega)^d), \\
u_k - v_k &\to 0 &&\text{ strongly in } C^0([0,T];L^r(\Omega)^d).
\end{aligned}\end{equation}  Let us consider the test function $\zeta_u = u_k - v_k$
in \eqref{Sys:GalerkinVar}$_3$, which gives 
\begin{equation*}
\begin{aligned}
&(S_u,\eps(u_k-v_k))_\Omega +(\alpha p_k,\div(u_k-v_k))_\Omega \\ & =(\C_v\pt \eps(u_k),\eps(u_k-v_k))_\Omega+(\C (\phi_k)(\eps(u_k)-  \T(\phi_k)) , \eps(u_k - v_k))_\Omega 
\end{aligned}
\end{equation*}
and together with the coercivity property
\eqref{Ass:Elastic2}$_1$ it yields
\begin{equation} \label{Eq:DeriveStrongU}
\begin{aligned}&(S_u,\eps(u_k-v_k))_\Omega+(\alpha p_k,\div(u_k-v_k))_\Omega  - (\C(\phi_k)(\eps(v_k) - \T(\phi_k)), \eps(u_k - v_k))_\Omega \\
	&= (\C_v\pt \eps(u_k-v_k),\eps(u_k-v_k))_\Omega +(\C_v\pt \eps(v_k),\eps(u_k-v_k))_\Omega \\&\quad +(D_2 W(\varphi_k,\eps(u_k)) - D_2 W(\varphi_k,\eps(v_k)), \eps(u_k - v_k))_\Omega \\
& \geq C \ddt \|\eps(u_k-v_k)\|_{L^2}^2 + C \|\eps(u_k - v_k)\|_{L^2}^2  + (\C_v\pt \eps(v_k),\eps(u_k-v_k))_\Omega
\end{aligned}
\end{equation}
We have used that we can write $(\C_v \F,\mathcal{G})_\Omega=(\C_v^{1/2} \F,\C_v^{1/2} \mathcal{G})_\Omega$ for any $\F,\mathcal{G} \in \mathbb{R}^{d \times d}$ where $\C_v^{1/2}$ exists according to the square-root theorem (interpreting $\C_v$ as a linear operator from the Hilbert space $\mathbb{R}^{d \times d}$ to $\mathbb{R}^{d \times d}$) as $\C_v$ is symmetric and positive definite.
We integrate the inequality over $(0,T)$ and pass to the limit $k \to \infty$. 
The third term on the right-hand side converges to zero as $\eps(u_k-v_k)$ converges weakly to zero in $L^2(\OT)^{d \times d}$ by \cref{Conv:UKVK} and $\pt \eps(v_k)$ converges strongly to $\pt \eps(u)$ in $L^2(\OT)^{d \times d}$ by \cref{Conv:VK}. Similarly, the first and third terms on the left-hand side converge to zero.
The $\|\eps(u_{k,0})-\eps(v_{k,0})\|_{L^2}^2$ appears after the time integration and goes to zero as $\eps(u_{k,0}) \to \eps(u_0)$ in $L^2(\Omega)^{d \times d}$ by construction of the discrete initial condition, and $\eps(v_k(0)) \to \eps(u(0))$  by the strong convergence of $v_k$ in $C([0,T];H_0^1(\Omega)^d)$, see \cref{Conv:VK}. It remains to treat the second term on the left-hand side of the inequality. First, we integrate by parts and exploit the homogeneous Dirichlet boundary condition of $u_k-v_k$ to get
$$(\alpha p_k,\div(u_k-v_k))_\OT =-(\alpha' p_k,u_k-v_k)_\OT - (\alpha \nabla p_k,u_k-v_k)_\OT.
$$
Here, we use the strong convergence of $u_k-v_k$ to zero in $L^2(\OT)^d$ together with the weak convergence of $p_k$ and $\nabla p_k$ in $L^2(\OT)$. We note that in the case of a state-dependent Biot--Willis function $\phi_k \mapsto \alpha(\phi_k)$, an additional term $\nabla \phi_k$ would appear after integration by parts and we would require strong convergence of this term to pass to the limit.
Collecting all the derived convergences, we see from \cref{Eq:DeriveStrongU} that
\begin{align*}
\|\eps(u_k - v_k)\|_{L^2(\Omega_T)}^2 \to 0 \quad \text{ as } k \to \infty.
\end{align*}
By Korn's inequality this shows that $u_k-v_k$ converges strongly to zero in the space $L^2(0,T;H^1_0(\Omega))$ and hence 
\begin{equation} \label{Eq:StrongU2} u_k \to u \quad \text{ strongly in }  L^2(0,T;H^1_0(\Omega)).
\end{equation} %

Next, we consider the discretized pressure equation and, having already derived the strong convergence of $u_k$, we are able to conclude a strong convergence result for $p_k$. As $\cup_k Z_k$ is dense in $H^1(\Omega)$, we can choose a sequence $\{q_k\}$ with $q_k(t) \in Z_k$ for a.e. $t \in (0,T)$ and $q_k \to p$ strongly in $L^2(0,T;H^1(\Omega))$. Then, it also yields $p_k-q_k \to 0$ weakly in $L^2(0,T;H^1(\Omega))$. We consider the test function $\xi_p=p_k-q_k$ in \eqref{Sys:GalerkinVar}$_5$, which gives after integration
\begin{equation} \label{Eq:PkQk} \begin{aligned} \|p_k-q_k\|^2_{L^2(\OT)} &= (q_k,p_k-q_k)_\OT  -(\alpha M(\phi_k) \div u_k, p_k-q_k)_\OT 
	\\&\quad +(M(\phi_k) \theta_k,p_k-q_k)_\OT .\end{aligned}\end{equation}
The terms on the right-hand side all converge to zero as $k \to \infty$.
In fact, $q_k$ strongly converges in $L^2(\OT)$ and $p_k-q_k$ weakly converges to $0$ in $L^2(\OT)$. Furthermore, we have just derived that $\div u_k$ strongly converges in $L^2(\OT)$ and together with the boundedness of $M$ and the strong convergence of $\phi_k$, we conclude that $M(\phi_k) \div u_k$ converges strongly in $L^2(\OT)$. 
 The third term remains to be treated. %
 We have
\begin{equation} \label{Eq:PuzzlePiece} \begin{aligned} &\left( M(\phi_k) \theta_k,   p_k-q_k \right)_\OT\\
&=\left( M(\phi_k)\theta_k-M(\phi)\theta,  p_k-q_k  \right)_\OT+ \left( M(\phi)\theta,  p_k-q_k \right)_\OT \\
&\leq \|M(\phi_k)\theta_k-M(\phi)\theta\|_{L^2(H^1)'} \|p_k-q_k\|_{L^2(H^1)}+ \left( M(\phi)\theta,  p_k-q_k \right)_\OT.
\end{aligned} \end{equation}
We observe that the second term converges to zero as $M(\phi) \theta \in L^2(\OT)$ and $p_k-q_k$ converges weakly to zero in $L^2(\OT)$. From the weak convergence of $p_k-q_k$ in $L^2(0,T;H^1(\Omega))$ we are also able to conclude  that $p_k-q_k$ is bounded in said space. We prove that the first norm in the first term on the right-hand side goes to zero. %
We first note that
$$M(\phi_k) \theta_k - M(\phi) \theta = M(\phi_k) (\theta_k-\theta) - (M(\phi_k)-M(\phi)) \theta,$$
and here the second term converges to zero in $L^2(\OT)$ by the dominated convergence theorem as $M(\phi_k)\theta \to M(\phi) \theta$ a.e. in $\OT$ and the sequence is uniformly bounded in $L^2(\Omega)$. The first term goes to zero as
$$\begin{aligned} \|M(\phi_k) (\theta_k-\theta)\|_{L^2(H^1)'} &=\sup_{\|\xi\|_{L^2(H^1)} \leq 1} |(\theta_k-\theta,M(\phi_k)\xi)_\OT |  \\ 
	&\leq \sup_{\|\xi\|_{L^2(H^1)}\leq 1} \|\theta_k-\theta\|_{L^2(W^{1,4/3})'} \|M(\phi_k) \xi\|_{L^2(W^{1,4/3})} 
	\end{aligned}$$
and we have for any $\xi \in L^2(0,T;H^1(\Omega))$
$$\begin{aligned}
	\|M(\phi_k) \xi\|_{L^2(L^{4/3})}  &\leq M_\infty \|\xi\|_{L^2(L^{4/3})} \leq C  \|\xi\|_{L^2(L^2)}, \\
	\|M(\phi_k) \nabla \xi\|_{L^2(L^{4/3})}  &\leq M_\infty \|\nabla \xi\|_{L^2(L^{4/3})} \leq C  \|\nabla \xi\|_{L^2(L^2)}, \\
	\|M'(\phi_k) \nabla \phi_k \xi\|_{L^2(L^{4/3})}  &\leq M'_\infty \|\nabla \phi_k\|_{L^\infty(L^2)} \|\xi\|_{L^2(L^{4})} \leq C  \|\xi\|_{L^2(H^1)}, 
\end{aligned}$$
where we used H\"older's inequality with $\frac{1}{2}+\frac{1}{4}=\frac{3}{4}$ and the $k$-uniform bound of $\nabla \phi_k$ in $L^\infty(0,T;L^2(\Omega))$. Thus, $M(\phi_k) (\theta_k-\theta)$ converges to zero in $L^2(0,T;(H^1(\Omega))')$.  Hence, the term in \eqref{Eq:PuzzlePiece} and likewise in \eqref{Eq:PkQk} goes to zero.  We conclude
$p_k-q_k \to 0$ in $L^2(\OT)$ and 
$$p_k \to p \quad \text{ strongly in}L^2(\OT).$$

Again, we consider the discretized pressure equation and, having already derived the strong convergence of $u_k$ and $p_k$, we are able to conclude a strong convergence result for $\theta_k$. As $\cup_k Z_k$ is dense in $L^2(\Omega)$, we can choose a sequence $\{\Theta_k\}$ with $\Theta_k(t) \in Z_k$ for a.e. $t \in (0,T)$ and $\Theta_k \to \theta$ strongly in $L^2(\OT)$. Then, it also yields $\theta_k-\Theta_k \to 0$ weakly in $L^2(\OT)$. We consider the test function $\xi_p=\theta_k-\Theta_k$ in \eqref{Sys:GalerkinVar}$_5$, which gives after integration
\begin{equation} \label{Eq:PkQk2} \begin{aligned} M_0^2 \| \theta_k-\Theta_k\|^2_{L^2(\OT)}  &\leq \|M(\phi_k)(\theta_k-\Theta_k)\|^2_{L^2(\OT)} \\ &= (M(\phi_k) \Theta_k ,\theta_k-\Theta_k)_\OT + (p_k,\theta_k-\Theta_k)_\OT  \\ &\quad  +(\alpha M(\phi_k) \div u_k, \theta_k-\Theta_k)_\OT.
\end{aligned} \end{equation}
The right-hand side is converging to zero as $k \to \infty$ by the weak convergence of $\theta_k-\Theta_k$ and the strong convergences of $\Theta_k$, $p_k$, $\div u_k$. We conclude
$$\theta_k \to \theta \quad \text{ strongly in} L^2(\OT).$$

		\noindent\textbf{Step 4 (Energy inequality and initial conditions).} We note that the energy inequality \cref{Est:Final} holds in a continuous setting in the sense that we may replace $\phi_k$ by $\phi$ and so on. This is achieved by taking the limit inferior as $k \to \infty$ and using that norm are weakly/weakly* lower semicontinuous. Moreover, we apply the Fatou lemma on the non-negative continuous function $\Psi$ to achieve
		$$\int_\Omega \Psi(\phi(x)) \dx \leq \liminf_{k \to \infty}\int_\Omega \Psi(\phi_k(x)) \dx.$$
		Thus, the quintuple $(\phi,\mu,u,\theta,p)$ satisfies \cref{Est:FinalCont} as stated in \cref{Theorem}.
  
		Regarding the initial conditions, we note that it holds $\phi_k(0) \to \phi(0)$ in $L^r(\Omega)$ as $k \to \infty$ according to the strong convergence \cref{Convergence:Strong} by setting $t=0$. However, it holds $\phi_k(0)=\phi_{k,0} \to \phi_0$ in $H^1(\Omega)$ by the definition of $\phi_{k,0}=\Pi_{Z_k} \phi_0$ and the properties of the orthogonal projection. Thus, we obtain $\phi(0)=\phi_0$ in $L^r(\Omega)$ by the uniqueness of limits. Moreover, we make use of the embedding, see \cite[Lemma II.5.9]{boyer2012mathematical}
		$$C([0,T];L^r(\Omega)) \cap L^\infty(0,T;H^1(\Omega)) \hookrightarrow C_w([0,T];H^1(\Omega)),$$
		to infer that $\phi$ is weakly continuous with values in $H^1(\Omega)$ and thus, the initial condition $\phi_0$ is satisfied in the sense
		$$\langle w,\phi(t) \rangle_{H^1} \to \langle w,\phi_0 \rangle_{H^1} \quad \forall w \in (H^1(\Omega))'.$$ In the same manner, we obtain $u(0)=u_0$ in $L^r(\Omega)$ and $\theta(0)=\theta_0$ in $L^2(\Omega)$.
  \smallskip
		
		\noindent\textbf{Step 5 (Limit process).} In this step, we pass to the limit $k \to \infty$ in the $k$-th Galerkin system using the convergences that we have derived in \cref{Convergence:Weak}--\cref{Convergence:Strong2}. First, we multiply each of the Galerkin equations \cref{Sys:GalerkinVar} by an arbitrary function $\eta \in C_0^\infty(0,T)$ and integrate over $(0,T)$, giving
		\begin{equation} \label{Sys:GalerkinVarInt} \begin{aligned}
				\langle \pt \phi_k,\eta z\rangle_{L^2(H^1)} + (m(\phi_k)\nabla \mu_k,\eta\nabla z)_\OT&=  (S_\phi^k,\eta z)_\OT, \\
				- (\mu_k,\eta z)_\OT+\gamma  (\nabla\phi_k,\eta \nabla z)_\OT + (\Psi'(\phi_k),\eta z)_\OT &= -(\delta_{\phi_k} \E_\theta^k + \delta_{\phi_k} \E_u^k,\eta z)_\OT,    \\
				(\C_v\pt \eps(u_k),\eta \eps(y))_\OT +\big(\mathbb{C}(\varphi_k)\left({ \varepsilon}( u_k)-\mathcal{T}(\varphi_k)\right),\eta \nabla y\big)_\OT &=  (\alpha p_k,\eta \div y)_\OT, \\
				\langle\pt \theta_k,\eta z\rangle_{L^2(H^1)} + (\kappa(\phi_k)\nabla p_k,\eta \nabla z)_\OT  &= (S_\theta^k,\eta z)_\OT, \\
				(M(\phi_k)(\theta_k-\alpha \div u_k), \eta z)_\OT   &= (p_k,\eta z)_\OT,   
		\end{aligned}\end{equation}
		for any $z \in Z_k$, $y \in Y_k$ and $\eta \in C_0^\infty(0,T)$.
		
		Taking the limit $k \to \infty$ in the linear terms such as $(p_k,\eta z)_{\OT}$ follows directly due to the weak convergences \cref{Convergence:Weak}. Therefore, we only study the nonlinearities. We notice that the nonlinear functions terms depending on $\phi_k$ (such as $\kappa(\phi_k)$) are bounded. The convergence of these terms can be treated using the strong convergence of $\phi_k$, see \cref{Convergence:Strong}, and the Lebesgue dominated convergence theorem. By the weak-strong convergence lemma, we can take the limit in all product-structured terms such as $(m(\phi_k)\nabla \mu_k,\eta \nabla z)_{\OT}$. 
  
  Due to the strong convergence $u_k \to u$ in $L^2(0,T;H_0^1(\Omega))$ and the growth estimates of $W$ and its derivatives, see \cref{Ass:Elastic2}, we can pass to the limit in $\delta_{\phi_k} \E_u^k$. Further, we use the strong convergence of $\theta_k$ and $u_k$ to pass the limit in the quadratic nonlinearities in $\delta_{\phi_k} \E_\theta^k$.
  
  Then we use that $\cup_k Z_k$ is dense in $H^1(\Omega)$ and $\cup_{k} Y_k$ is dense in $H_0^1(\Omega)^d$. Together with the fundamental lemma of calculus of variations, it yields that $(\phi,\mu,u,\theta,p)$ is a weak solution to the Cahn--Hilliard--Biot system in the sense as stated in \cref{Theorem}. \end{proof}

  We have proved the existence of a weak solution. Next, we prove the uniqueness and continuous dependence on the data of the said solution in the case of stricter assumptions, as stated in \cref{Thm:Uniqueness}.

  \begin{proof}[Proof of \cref{Thm:Uniqueness}]
Now, we are in the setting that $m(\phi)=m$, $\C(\phi)=\C$, $\kappa(\phi)=\kappa$, $\alpha(\phi)=\alpha$, $M(\phi)=M$ are constant and $\T(\phi)$ is affine linear, that is, $\T(\phi)=\T_1+\T_2\phi$. In particular, this implies $$D_\phi W(\phi,\eps(u))=-\C(\eps(u)-\T_1-\T_2\phi):\T_2.$$
Moreover, due to $\alpha'=M'=0$, we have $\delta_\phi \E_\theta=0$.
In this case, we can easily prove that it holds $\Psi'(\phi) \in L^2(\OT)$ and $\phi \in L^2(0,T;H^2_n(\Omega))$ by solving the equation of the chemical potential $\mu$ for $\Delta \phi$, see also \cite[Section 5.2]{garcketumormechanics}.
Therefore, the inverse Neumann--Laplacian is well-defined when applied to $\phi$.

		We consider two weak solutions $(\phi_1,\mu_1,u_1,\theta_1,p_1)$ and $(\phi_2,\mu_2,u_2,\theta_2,p_2)$, and we denote their difference by $\phi=\phi_1-\phi_2$ and in the same way for the other variables. Subtracting their weak forms, we obtain
			\begin{equation} \label{Eq:UniqueTest1}\begin{aligned}
			\langle\pt \phi,\zeta_\phi\rangle_{H^1} + m(\nabla \mu,\nabla \zeta_\phi)_\Omega&=(S_\phi(\phi_1,\eps(u_1),\theta_1),\zeta_\theta) \\ &\quad -(S_\phi(\phi_2,\eps(u_2),\theta_2),\zeta_\phi)_\Omega, \\
			- (\mu,\zeta_\mu)_\Omega+\gamma  (\nabla\phi,\nabla \zeta_\mu)_\Omega + (\Psi'(\phi_1)\!-\!\Psi'(\phi_2),\zeta_\mu)_\Omega &= -(\C(\eps(u)\!-\!\T_2\phi),\T_2\zeta_\mu)_\Omega,  \\
			(\C_v \eps(\pt u),\eps(\zeta_u))_\Omega+\big(\C\big({ \varepsilon}( u)-\T_2\varphi\big),\eps(\zeta_u)\big)_\Omega &= \alpha(p,\div \zeta_u)_\Omega,  \\
			\langle\pt \theta,\zeta_\theta\rangle_{H^1} + \kappa(\nabla p,\nabla\zeta_\theta)_\Omega  &= (S_\theta(\phi_1,\eps(u_1),\theta_1),\zeta_\theta) \\ &\quad -(S_\theta(\phi_2,\eps(u_2),\theta_2),\zeta_\theta)_\Omega, \\
			(p,\zeta_p)_\Omega   &= M(\theta-\alpha \div u,\zeta_p)_\Omega,  
		\end{aligned} \end{equation}
		for any $\zeta_\phi, \zeta_\mu,\zeta_\theta,\zeta_p \in H^1(\Omega)$, $\zeta_u \in H^1_0(\Omega)$. 
		We have $p=M(\theta-\alpha \div u)$ in $L^2(\OT)$. Moreover, as the right-hand side is in $L^2(0,T;(H^1(\Omega))')$, we obtain $\pt p = M \pt \theta- \alpha \div \pt u$ in $L^2(0,T;(H^1(\Omega))')$. From here, we see that
		\begin{equation} \label{Eq:UniqueTest2} \langle\pt p,\zeta_\theta\rangle_{H^1} + M\kappa(\nabla p,\nabla\zeta_\theta)_\Omega  = (S_\theta(\phi_1)-S_\theta(\phi_2),\zeta_\theta)_\Omega - \alpha (\div \pt u,\zeta_\theta)_\Omega,
			\end{equation} 
		for any $\zeta_\theta \in H^1(\Omega)$ and the initial condition $p(0)=p_0:=M \theta_0 - \alpha \div u_0 \in L^2(\Omega)$.
				\medskip
		
		\noindent\textbf{First testing.} Taking the test functions $\zeta_\phi=(-\Delta)^{-1} \phi$ in \eqref{Eq:UniqueTest1}$_1$ and $\zeta_\mu=m\phi$ in \eqref{Eq:UniqueTest1}$_2$, it yields
		\begin{equation*} 
			\begin{aligned}
				\langle \p_t \phi,(-\Delta)^{-1} \phi\rangle_{H^1} + m (\nabla \mu,\nabla (-\Delta)^{-1} \phi)_\Omega &= (S_\phi(\phi_1)-S_\phi(\phi_2),(-\Delta)^{-1}\phi)_\Omega, \\
				-m(\mu,\phi)_\Omega+ m\gamma (\nabla \phi,\nabla \phi)_\Omega &= - m(\Psi'(\phi_1)-\Psi'(\phi_2),\phi)_\Omega \\ &\quad\, -m(\C(\eps(u)\!-\!\T_2\phi),\T_2\varphi)_\Omega  
			\end{aligned}
		\end{equation*} 
		Exploiting the property $(\nabla \mu, \nabla (-\Delta)^{-1} \phi)_\Omega = (\mu,\phi)_\Omega$ of the Neumann--Laplace operator,  after adding the equations and canceling, it yields
		\begin{equation} \label{Eq:UniqueDifference}
			\begin{aligned}
				&\langle \p_t \phi, (-\Delta)^{-1} \phi\rangle_{H^1} + m \gamma \|\nabla \phi\|_{L^2}^2 + m(\Psi'(\phi_1)-\Psi'(\phi_2),\phi)_\Omega \\
				&=(S_\phi(\phi_1)-S_\phi(\phi_2),(-\Delta)^{-1} \phi)_\Omega - m(\C(\eps(u)\!-\!\T_2\phi),\T_2\varphi)_\Omega. 
			\end{aligned}
		\end{equation}
  The graph norm $\|\nabla (-\Delta)^{-1} \cdot\|_{L^2}$ is equivalent to the usual norm of $(H^1(\Omega))'$. Thus, we set $\|\cdot\|_{(H^1)'}=\|\nabla (-\Delta)^{-1} \cdot\|_{L^2}$. First, we note that we may rewrite the first term on the left-hand side of \cref{Eq:UniqueDifference} as %
		\begin{equation}  
			\begin{aligned} 
				\langle \p_t \phi, (-\Delta)^{-1} \phi\rangle_{H^1} &= \langle -\Delta (-\Delta)^{-1} \p_t \phi, (-\Delta)^{-1} \phi\rangle_{H^1}  \\&= \langle\p_t \nabla (-\Delta)^{-1} \phi, \nabla (-\Delta)^{-1} \phi\rangle_{H^1} \\
				&= \frac12 \ddt \|\phi\|_{(H^1)'}^2.
			\end{aligned} \label{Eq:UniqueChain} 
		\end{equation}
		Using the semiconvexity of $\Psi$, see \cref{Ass:PsiNew}, it yields
		$$m(\Psi'(\phi_1)-\Psi'(\phi_2),\phi)_\Omega \geq - C_* m\|\phi\|_{L^2}^2,
		$$
		and consequently, we obtain by the Young inequality 
	\begin{equation} \label{Eq:L2EstimateInvLap}\begin{aligned}m(\Psi'(\phi_2)-\Psi'(\phi_1),\phi)_\Omega &\leq C_*m \|\phi\|_{L^2}^2 \\ &= C_*m (\nabla (-\Delta)^{-1} \phi,\nabla \phi)_\Omega \\ &\leq \frac{m\gamma}{4} \|\nabla \phi\|^2_{L^2} + \frac{C_*^2m}{\gamma} \|\nabla (-\Delta)^{-1} \phi\|_{L^2}^2.\end{aligned}\end{equation}
Here, we have seen that we can estimate any $L^2(\Omega)$-norm of $\phi$ using the Laplacian and inverse Laplacian operators. This property is also used in the next step.
		Using the Lipschitz continuity of $S_\phi$, see \cref{Ass:FunNew}, the right-hand side of \cref{Eq:UniqueDifference} can be treated  as follows:
		$$\begin{aligned}
			&(S_\phi(\phi_1,\eps(u_1),\theta_1)-S_\phi(\phi_2,\eps(u_2),\theta_2),(-\Delta)^{-1} \phi)_\Omega - m(\C(\eps(u)\!-\!\T_2\phi),\T_2\varphi)_\Omega  \\
			&\leq C \big( \|\phi\|_{L^2} + \|\eps(u)\|_{L^2}+\|\theta\|_{L^2} \big) \|(-\Delta)^{-1} \phi\|_{L^2} + C \|\eps(u)\|_{L^2} \|\phi\|_{L^2} + C\|\phi\|_{L^2}^2   \\
			&\leq \frac{m\gamma}{4} \|\nabla \phi\|_{L^2}^2  + C\|\phi\|_{(H^1)'}^2 +  C \|\eps(u)\|_{L^2}^2+ \frac{M}{4} \|\theta\|_{L^2}^2.
		\end{aligned}$$
		Therefore, applying this estimate and \cref{Eq:UniqueChain} to \cref{Eq:UniqueDifference}, it yields
		\begin{equation}  \label{Est:UniqueFinal1}
			\begin{aligned}
				\frac12 \ddt \|\phi\|_{(H^1)'}^2 + \frac{m \gamma}{2} \|\nabla \phi\|_{L^2}^2 \leq C\|\phi\|_{(H^1)'}^2 +C \|\eps(u)\|_{L^2}^2+ \frac{M}{4} \|\theta\|_{L^2}^2.
			\end{aligned}
		\end{equation} \smallskip
	
\noindent\textbf{Second testing.}		Next, we consider the test functions $\zeta_u=\pt u$ in \eqref{Eq:UniqueTest1}$_3$, $z_p=\theta$ in  \eqref{Eq:UniqueTest1}$_5$, $z_\theta=p$ in \eqref{Eq:UniqueTest2} to obtain
		$$\begin{aligned} C_\C \|\pt \eps(u)\|_{L^2}^2 + \frac{1}{2} \ddt \|\C^{1/2}\eps(u)\|_{L^2}^2  &= \big(\C\T_2\varphi,\eps(\pt u)\big)_\Omega+\alpha  (p,\div \pt u)_\Omega,  \\
			  M \|\theta\|_{L^2}^2  &= (p,\theta)_\Omega + \alpha M (\div u,\theta)_\Omega, \\
			\frac12 \ddt \|p\|_{L^2}^2 + M \kappa \|\nabla p\|_{L^2}^2 &= (S_\theta(\phi_1,\eps(u_1),\theta_1)-S_\theta(\phi_2,\eps(u_2),\theta_2),p)_\Omega \\ &\quad - \alpha (\div \pt u,p)_\Omega.
		\end{aligned}$$
		We notice that the mixed term $\alpha(p,\div \pt u)_\Omega)$ cancels. Adding the tested equations gives
		$$\begin{aligned}
			&C_\C \|\pt\eps(u)\|_{L^2}^2+\frac{1}{2}\ddt \|C_\C^{1/2}\eps(u)\|_{L^2}^2+ M \|\theta\|_{L^2}^2+\frac12 \ddt \|p\|_{L^2}^2+ M\kappa \|\nabla p\|_{L^2}^2 \\
   &= (\C(\T_2\phi),\eps(\pt u))_\Omega  +(p,\theta)_\Omega + \alpha M (\div u,\theta)_\Omega \\&\quad +(S_\theta(\phi_1,\eps(u_1),\theta_1)-S_\theta(\phi_2,\eps(u_2),\theta_2),p)_\Omega 
		\end{aligned}$$
  We apply the H\"older inequality on the terms on the right-hand side of this equality and with the boundedness of $S_\theta$, see \cref{Ass:Fun}, we obtain
  		\begin{equation}\label{Est:Unique2}\begin{aligned}
			&C_\C \|\pt\eps(u)\|_{L^2}^2+\frac12 \ddt \|\C^{1/2}\eps(u)\|_{L^2}^2+M \|\theta\|_{L^2}^2+\frac12 \ddt \|p\|_{L^2}^2+ M\kappa \|\nabla p\|_{L^2}^2 \\
			&\leq  C \|\phi\|_{L^2} \|\pt \eps(u)\|_{L^2}+C \|p\|_{L^2} \|\theta\|_{L^2} + C \|\div u\|_{L^2} \|\theta\|_{L^2} \\ &\quad + C \big( \|\phi\|_{L^2}  
+ \|\eps(u)\|_{L^2} + \|\theta\|_{L^2} \big) \|p\|_{L^2}
		\end{aligned}\end{equation}
  We apply the Young inequality on the right-hand side of the inequality. In particular, we have for the terms where $\phi$ appears
  $$\begin{aligned}&C \|\phi\|_{L^2} \|\pt \eps(u)\|_{L^2}+ C \big( \|\phi\|_{L^2}  
+ \|\eps(u)\|_{L^2} + \|\theta\|_{L^2} \big) \|p\|_{L^2} \\ &\leq 
  	C\|\phi\|_{L^2}^2 + \frac{C_\C}{2} \|\pt \eps(u)\|_{L^2}^2 + C\|\eps(u)\|_{L^2}^2 + \frac{M}{4} \|\theta\|_{L^2}^2 + C \|p\|_{L^2}^2  \\
  	&\leq C \|\phi\|_{(H^1)'}^2 + \frac{m\gamma}{4} \|\nabla \phi\|_{L^2}^2+ \frac{C_\C}{2} \|\pt \eps(u)\|_{L^2}^2 + C\|\eps(u)\|_{L^2}^2 + \frac{M}{4} \|\theta\|_{L^2}^2+ C \|p\|_{L^2}^2 ,
  \end{aligned}$$
  where we applied \cref{Eq:L2EstimateInvLap} to estimate $\phi$ in the $L^2(\Omega)$-norm in the last step.
  Inserting this estimate into \cref{Est:Unique2}, we obtain
\begin{equation}  \label{Est:UniqueFinal2}	\begin{aligned}
		&\frac{C_\C}{2} \|\pt\eps(u)\|_{L^2}^2+\frac12 \ddt \|\C^{1/2}\eps(u)\|_{L^2}^2+\frac{3M}{4} \|\theta\|_{L^2}^2+\frac12 \ddt \|p\|_{L^2}^2+ M\kappa \|\nabla p\|_{L^2}^2 \\
		&\leq C \|\phi\|_{(H^1)'}^2 + \frac{m\gamma}{4} \|\nabla \phi\|_{L^2}^2+  C \|p\|_{L^2}^2+C\|\eps(u)\|_{L^2}^2.
    \end{aligned}\end{equation}
\smallskip

		\noindent\textbf{Together.}
  We add \cref{Est:UniqueFinal1} to   \cref{Est:UniqueFinal2}, which gives
  \begin{equation*}  \begin{aligned}
			&\frac12 \ddt \|\phi\|_{(H^1)'}^2 +\frac{m \gamma}{4} \|\nabla \phi\|_{L^2}^2+\frac{C_\C}{2} \|\pt \eps(u)\|_{L^2}^2 +\frac12 \ddt \|\C^{1/2}\eps(u)\|_{L^2}^2 + \frac{M}{2} \|\theta\|_{L^2}^2  \\ &\quad  +\frac{1}{2} \ddt \|p\|_{L^2}^2 + M\kappa \|\nabla p\|_{L^2}^2 \\
  		&\lesssim   \|\phi\|_{(H^1)'}^2 + \|p\|_{L^2}^2 + \|\eps(u)\|_{L^2}^2.
    \end{aligned}\end{equation*}
		We integrate and apply the Gronwall inequality, which gives, after taking the essential supremum of $t\in(0,T)$,
		\begin{equation*} %
			\begin{aligned}
				& \|\phi\|_{L^\infty(H^1)' \cap L^2(L^2)}^2+ \|\pt u\|_{L^2(H_0^1) }^2+\| u\|_{L^\infty(H_0^1) }^2 + \|\theta\|_{L^2(L^2)}^2 +   \|p\|_{L^\infty(L^2) \cap L^2(H^1)}^2\\ &\lesssim  \|\phi_{0,1}-\phi_{0,2} \|_{(H^1)'}^2 +\|u_{0,1}-u_{0,2} \|_{H^1}^2 + \|p_{0,1}-p_{0,2}\|_{L^2}^2.
			\end{aligned}
		\end{equation*}
	As noted before, we may use the representation of the initial pressure to estimate the last term on the right-hand side as follows:
	$$\|p_{0,1}-p_{0,2}\|_{L^2}^2 \leq C \|\theta_{0,1}-\theta_{0,2}\|_{L^2}^2 + C \|u_{0,1}-u_{0,2}\|_{H^1}^2.$$
    Moreover, from \cref{Eq:UniqueTest1}$_2$ we have that
    $$\|\mu\|_{L^2(H^1)'}^2 \lesssim \|\phi\|_{L^2(H^1)} + \|\eps(u)\|_{L^2(L^2)}^2.$$
  This completes the continuous dependence on the data as stated in \cref{ContDep}.
Moreover, in the case of $\phi_{0,1}=\phi_{0,2}$, $u_{0,1}=u_{0,2}$ and $\theta_{0,1}=\theta_{0,2}$, the weak solutions coincide. 
	\end{proof}
\section{Well-posedness of the regularized Cahn--Hilliard--Biot model} \label{sec:wellreg}
As written before, we distinguish the cases in which the Biot--Willis parameter $\alpha$ is taken to be spatially or state dependent. After studying the first case in the previous section, we now study the case of a Biot--Willis function $\phi \mapsto \alpha(\phi)$, that is, we study the $\nu$-regularized Cahn--Hilliard--Biot system \cref{Sys:CHBreg} as discussed in \cref{sec:model}.

The existence theorem to the Cahn--Hilliard--Biot system with the Biot--Willis function $\phi \mapsto \alpha(\phi)$ reads as follows.

\begin{theorem}[Existence of a weak solution] \label{RegTheorem}
Let Assumption \ref{Assumption} hold. In addition, we assume $\phi_0 \in H^2_n(\Omega)$ and $\nu>0$ arbitrary but fixed.
    Then there exists a global-in-time solution 
    $$\begin{aligned}
    \phi &\in C_w([0,T];H^2_n(\Omega)) \cap C^0([0,T];W^{1,r}(\Omega)) \cap H^1(0,T;(H^1(\Omega))'), \\
    \mu &\in L^2(0,T;H^1(\Omega)), \\
    u &\in H^1(0,T;H_0^1(\Omega)^d), \\
    p &\in L^2(0,T;H^1(\Omega)) \cap L^\infty(0,T;L^2(\Omega)), \\
    \theta &\in C_w([0,T];L^2(\Omega)) \cap H^1([0,T];(H^1(\Omega))'),
\end{aligned}$$
with $r<6$ for $d=3$ or $r<\infty$ for $d\leq 2$, to the $\nu$-regularized Cahn--Hilliard--Biot system \cref{Sys:CHB} in the sense that the variational equations
    \begin{equation} \begin{aligned}
(\pt \phi,\zeta_\phi)_\Omega + (m(\phi) \nabla \mu,\nabla \zeta_\phi)_\Omega&=(S_\phi,\zeta_\phi)_\Omega, \\
 - (\mu,\zeta_\mu)_\Omega\!+\!\gamma  (\nabla\phi,\nabla \zeta_\mu)_\Omega\!-\!\nu (\Delta \phi,\Delta \zeta_\mu)_\Omega \!+\!(\Psi'(\phi),\zeta_\mu)_\Omega &= (\delta_\phi \E_u + \delta_\phi \E_\theta,\zeta_\mu)_\Omega,    \\
(\C_v \p_t \eps(u),\eps(\zeta_u))_\Omega + \big(\mathbb{C}(\varphi)\left({ \varepsilon}( u)-\mathcal{T}(\varphi)\right),\eps(\zeta_u)\big)_\Omega &= (\alpha(\varphi)p,\div \zeta_u)_\Omega+(S_u,\zeta_u)_\Omega,  \\
(\pt \theta,\zeta_\theta)_\Omega + (\kappa(\phi)\nabla p,\nabla\zeta_\theta)_\Omega   &= (S_\theta,\zeta_\theta)_\Omega, \\
 (p,\zeta_p)_\Omega  &= (\theta-\alpha(\phi) \div u,M(\phi)\zeta_p)_\Omega,
\end{aligned} \label{Sys:VariationalReg} \end{equation}
are satisfied for any test functions $\zeta_\phi \in H^2_n(\Omega)$, $\zeta_\mu,\zeta_\theta,\zeta_p  \in H^1(\Omega)$, $\zeta_u \in H_0^1(\Omega)^d$, and the initials $\phi(0)=\phi_0$, $u(0)=u_0$, $\theta(0)=\theta_0$.
Moreover, the solution fulfills for any $t \in (0,T]$ the energy inequality
\begin{equation} \label{RegEst:FinalCont}\begin{aligned}
&\|\phi(t)\|_{H^2}^2 + \|\nabla \mu\|^2_{L^2_t(L^2)} + \|\nabla p\|^2_{L^2_t(L^2)} +\|\Psi(\phi(t))\|_{L^1} + \|u\|_{H^1(H^1_0)}^2  + \|p(t)\|_{L^2}^2 \\
&\lesssim 1+\|\phi_0\|_{H^1}^2 +\|u_0\|_{L^2}^2+ \|\theta_0\|_{L^2}^2.
\end{aligned}\end{equation}
\end{theorem}

To prove the existence theorem, we follow the Galerkin procedure by deriving energy estimates on a discrete level before returning to the continuous problem. As the proof structure is the same as before, we simply point out the differences and the influence of the new term with the regularization parameter. We also see the influence of the Biot--Willis function $\alpha$. However, this will only make a significant difference later in the limit passage.

\begin{proof}[Proof of \cref{RegTheorem}] As described, we follow the Galerkin procedure to prove the existence of a weak solution. The existence of a discrete solution follows as before. \smallskip
		
		\noindent\textbf{Energy estimate.} We consider the same test functions as before. We obtain almost the same time integrated estimate as before in \cref{Est:Integrated}. 
		In contrast to the estimate from before, the only difference is in the term $\|\Delta \phi_k\|_{L^2}^2$, giving higher regularity than before due to elliptic regularity. Moreover, the energies $\E_u$ and $\E_\theta$ in the estimate can be treated as before.
		We observe that in the estimate of the mean of $\mu_k$ the term  $\int_\Omega \delta_{\phi_k} \E_\theta^k  \dx$ has now an additional term due to $\alpha$ depending on $\phi_k$. However, the new term can be estimated by the H\"older and Young inequalities as usual, not changing the right-hand side for us. 
		Then, using $\phi_0 \in H^2_n(\Omega)$ and the results in \cite{garcketumormechanics,abels2024existence}, by an application of the Gronwall inequality, we get
\begin{equation} \label{RegEst:Final}\begin{aligned}
		&\|\phi_k(t)\|_{H^2}^2 + \|\nabla \mu_k\|^2_{L^2_t(L^2)} +\|\Psi(\phi_k(t))\|_{L^1} + \|u_k(t)\|_{H^1}^2 +\|\partial_t u_k\|_{L^2_t(H^1)}^2  \\
		&+ \|p_k(t)\|_{L^2}^2+ \|\nabla p_k\|^2_{L^2_t(L^2)} + \|\theta(t)\|_{L^2}^2\lesssim 1+\|\phi_0\|_{H^2}^2 + \|\theta_0\|_{H^1}^2+\|u_0\|_{H^1}^2.
\end{aligned}\end{equation}
Since the right-hand side is uniform in $k$, we can extend the existence interval by setting $T_k=T$ for any $k$.
We can extract weakly converging subsequences (that we denote by the same index by a typical abuse of notation).  In fact, we get the existence of limit functions $(\phi,\mu,u,\theta,p)$ such that \cref{Convergence:Weak} holds again and in addition, we have
\begin{equation} \label{RegConvergence:Weak}\begin{aligned}
		\phi_k &\to \phi &&\text{weakly* in } L^\infty(0,T;H^2(\Omega)). 
\end{aligned}\end{equation} 
Again, bounding the time derivatives $\pt \phi_k$, we obtain with the Aubin--Lions lemma the improved strong convergences
\begin{equation} \label{RegConvergence:Strong}\begin{aligned}
		\phi_k &\to \phi &&\text{strongly in } C([0,T];W^{1,r}(\Omega)), 
\end{aligned}\end{equation} 
wirh $r$ defined as before.
We proceed similarly as before and consider the test function $u_k-v_k$ with $v_k$ strongly converging to $u$ in $H^1(0,T;H_0^1(\Omega))$. Then we end up with the same estimate as in \cref{Eq:DeriveStrongU}, that is,
\begin{equation} \label{RegEq:DeriveStrongU}
	\begin{aligned}&(S_u,\eps(u_k-v_k))_\Omega-(\C_v\pt \eps(v_k),\eps(u_k-v_k))_\Omega+(\alpha(\phi_k) p_k,\div(u_k-v_k))_\Omega \\ 
		&\quad - (\C(\phi_k)(\eps(v_k) - \T(\phi_k)), \eps(u_k - v_k))_\Omega \\
		& \geq C \ddt \|\eps(u_k-v_k)\|_{L^2}^2 + C \|\eps(u_k - v_k)\|_{L^2}^2 .
	\end{aligned}
\end{equation}
We have already investigated all the terms besides the term involving $\alpha(\phi_k)$. Again, we apply integration by parts to obtain
$$(\alpha(\phi_k) p_k,\div(u_k-v_k))_\Omega=-(\alpha'(\phi_k) \nabla \phi_k p_k,u_k-v_k)_\Omega-(\alpha(\phi_k) \nabla p_k,u_k-v_k)_\Omega.$$
In both terms, $u_k-v_k$ converges strongly in $C([0,T];L^4(\Omega))$. Moreover, $\nabla \phi_k$ converges strongly in $C([0,T];L^4(\Omega)^d)$. Together with the weak convergence of $p_k$ and $\nabla p_k$ in $L^2(\OT)$, we observe that both terms converge to $0$. Hence, we obtain the strong convergence of $u_k$ to $u$ in $L^2(0,T;H_0^1(\Omega))$. Again, we can show that $p_k \to p$ in $L^2(\OT)$ and $\theta_k \to \theta$ in $L^2(\OT)$ and thus we can pass to the limit in the discretized system, recovering the regularized Cahn--Hilliard--Biot system. 
	\end{proof}

\section{Numerical simulations} \label{Sec:Numerics}
In this section, we present some numerical simulations to highlight the differences of the Cahn--Hilliard--Biot model with the established Cahn--Hilliard and Cahn--Larch\'e equations. We assume an endless supply of nutrients by assuming a growth function $S_\phi=\lambda \phi(1-\phi)$ with $\lambda>0$ being the tumor proliferation factor. Other approaches include a diffusion equation for nutrients; see \cite{colli2017asymptotic}. However, to ensure the boundedness of $S_\phi$, we replace $\phi$ by $\mathcal{C}(\phi)$ where $\mathcal{C}$ is the cutoff operator defined by $\mathcal{C}(\phi)=\max\{0,\min\{1,\phi\}\}$. Moreover, it would be straightforward to model the nutrients by their own reaction-diffusion equation, as done in \cite{fritz2019localnonlocal,fritz2019unsteady,garcke2017analysis}. However, this is not the focus of this work. We are interested in the effects of the Biot model on the evolution of tumor mass.

\subsection{Setup}
We consider the unit square domain $\Omega=[0,1]^2$ and the time domain $[0,T]$ with $T=1.5$, which are discretized with $\Delta x = 2^{-8}$ and $\Delta t = 2^{-7}$. As initials, we choose $\theta_0(x)=\frac12$ and $\phi_0(x)=\text{exp}(1-\frac{1}{1-h(x)})$ with $$h(x)=(\sin(14.4 x_1 + 11.2 x_2-12.8) + 1)  (8x_1-4.2)^2 + (\sin(16x_1-8) + 1)  (16x_2-8)^2,$$
while we select zero initial data for the other variables. The variational system \cref{Sys:Variational} is discretized in time using a semi-implicit Euler method by using the classical convex-concave splitting of the nonlinear functional $\Psi'=\Psi_e'+\Psi_c'$ into its expansive $\Psi_e'$ and contractive part $\Psi_c'$. In the case of $\Psi(\phi)=\frac14\phi^2(1-\phi^2)$, we set
$\Psi_e'(\phi)=\phi^3-\frac32\phi^2-\frac14\phi$ and $\Psi_c'(\phi)=\frac34\phi$.
We treat the expansive part explicitly and the contractive part implicitly. The three-way coupled nonlinear system is then solved by an iterative decoupling scheme, starting with the coupled elasticity-flow equations governing $(u,\theta,p)$ and afterward the Cahn--Hilliard model governing $(\phi,\mu)$. 
The nonlinear equations are solved by the Newton method in each iterative decoupling iteration. For each variable, we select the bilinear rectangular finite element space $Q_1$. The system is implemented in the finite element library FEniCS \cite{alnaes2015fenics}.

In the model, we choose the functions $S_\phi=5\phi(1-\phi)$, $S_\theta=S_u=0$,  $m(\phi)=10^{-16}+\frac12\phi^2(1-\phi)^2$, $\T(\phi)=\frac{3}{10} \phi I$, $\C_v=10^{-16}I$. We are interested in the influence of the Biot model on tumor growth and not viscoelastic effects.  In addition, the interfacial parameter $\gamma$ is selected as $\gamma=10^{-4}$. Following \cite{storvik2022cahn}, the  permeability $\kappa(\pf)$, compressibility $M(\pf)$, Biot--Willis coefficient $\alpha(\pf)$ and elasticity tensor $\mathbb{C}(\pf)$ are dependent on the phase-field variable $\phi$ through the interpolation function $\pi(\pf)$ in the following way
$$\kappa(\pf) = \kappa_0 + \pi(\pf)(\kappa_1-\kappa_0),$$
and in the same manner for the other material functions
where the interpolation function $\pi$ is defined by
$\pi(\pf) = -2\phi^3+3\phi^2$ for $\pf \in [0,1]$ and it is set to $0$ for $\pf<0$ and to $1$ for $\pf>1$. However, we choose the material parameters differently from \cite{storvik2022cahn} because tumor tissue is more compressible and permeable than normal tissue due to its loose and disorganized structure, which facilitates fluid leakage. Additionally, tumors tend to have a higher Biot--Willis coefficient because of their high fluid content and soft structure, meaning that fluid pressure strongly affects the solid tumor matrix. Therefore, we choose the limit values of the material parameters as written in \cref{table}. Generally, the compressibility and permeability should be a few magnitudes lower in a real-world scenario, but we increase their value to highlight the effects of the new couplings in the Cahn--Hilliard--Biot model. 
Moreover, the elasticity tensors $\mathbb{C}_0$ and $\mathbb{C}_1$ are written in the Voigt notation 
$$\mathbb{C}_i =\begin{pmatrix} 
    \lambda_i+2G_i & \lambda_i & 0  \\
    \lambda_i & \lambda_i+2G_i & 0 \\ 
    0 & 0 & G_i\end{pmatrix} $$
where $G_i$ and $\lambda_i$ denote the shear modulus and Lam\'e parameter, respectively, for $i=0$ (healthy tissue) and $i=1$ (tumor region). They are given by the formulas $G_i=\frac{E_i}{2+2\nu_i}$ and $\lambda_i = \frac{E_i\nu_i}{(1+\nu_i)(1-2\nu_i)}$, where $E_i$ and $\nu_i$ are the Young modulus and Poisson ratio, respectively; see \cref{table} for selected chosen values.

\begin{table}[htbp!]
\centering
\caption{Table of material parameters.}
\begin{tabular}{c|c||c|c}
 Symbol & Value &  Symbol & Value\\
\hline
 $\kappa_0$ & 0.5 &   $\kappa_1$ & 5  \\
 $M_0$ & 0.5 &  $M_1$ & 1  \\
 $\alpha_0$ & 0.5 &   $\alpha_1$ & 1  \\
 $E_0$ & 2.8 &   $E_1$ & 1.4  \\
  $\nu_0$ & 0.4 &   $\nu_1$ & 0.2  \\
\end{tabular}
\label{table}
\end{table}

\subsection{Results}
 As our first test, we compare the Cahn--Hilliard (CH), Cahn--Larch\'e (CL) and Cahn--Hilliard--Biot (CHB) systems by investigating the evolution of tumor mass $t \mapsto \int_\Omega \phi(t,x) \dx$. Here, the Cahn--Hilliard equation reads
\begin{equation*} \label{Sys:CH}\begin{aligned}
\partial_t \varphi - \div (m(\pf) \nabla \mu) &= S_\phi, \\
\mu +\gamma\Delta \varphi - \Psi'(\varphi)  &= 0,%
\end{aligned}\end{equation*}
and the Cahn--Larch\'e equation
\begin{equation*} \label{Sys:CHL}\begin{aligned}
\partial_t \varphi - \div (m(\pf) \nabla \mu) &= S_\phi, \\
\mu +\gamma\Delta \varphi - \Psi'(\varphi)  &= \frac{1}{2}\left({ \varepsilon}( u) - \mathcal{T}(\varphi)\right)\!:\!\mathbb{C}'(\varphi)\left({ \varepsilon}( u)  - \mathcal{T}(\varphi)\right) \\&\quad- \mathcal{T}'(\varphi)\!:\!\mathbb{C}(\pf)\left({ \varepsilon}( u) - \mathcal{T}(\varphi)\right),\\
\div\big(\mathbb{C}(\varphi)\left({ \varepsilon}( u)-\mathcal{T}(\varphi)\right)\big)  &= 0. %
\end{aligned}\end{equation*}

\begin{figure}[H]
    \centering
    \begin{tikzpicture}
	\definecolor{color0}{HTML}{4E79A7}
	\definecolor{color1}{HTML}{F28E2B}
	\definecolor{color2}{HTML}{E15759}
	\definecolor{color3}{HTML}{76B7B2}
	\definecolor{color4}{HTML}{59A14F}
	\definecolor{color5}{HTML}{EDC948}
	\definecolor{color6}{HTML}{B07AA1}
	\definecolor{color7}{HTML}{FF9DA7}
	\definecolor{color8}{HTML}{9C755F}
	\definecolor{color9}{HTML}{BAB0AC}
	\begin{axis}[
		height=.33\textheight,
        width=.96\textwidth,
		legend pos=north west,
		legend cell align={left},
        ymin=0.02, ymax=0.145,
		xmin=0, xmax=1.5,
		ytick={0.14, 0.1, 0.06,0.02},
		xtick={0,0.5, 1.0,1.5},
		]
  		\addplot[ultra thick,color0,smooth] table [x expr=\coordindex/100, y index=0]{Figures/cahn2.txt};
		\addlegendentry{Cahn--Hilliard}
  		\addplot[ultra thick,color1,smooth] table [x expr=\coordindex/100, y index=0]{Figures/larche2.txt};
		\addlegendentry{Cahn--Larch\'e}
		\addplot[ultra thick,color2,smooth] table [x expr=\coordindex/100, y index=0]{Figures/biot2.txt};
		\addlegendentry{Cahn--Hilliard--Biot}
	\end{axis}
\end{tikzpicture}
    \caption{Evolution of the tumor mass $t \mapsto \int_\Omega \phi(t,x) \dx$ for the different models.}
    \label{Fig:Volume}
\end{figure}
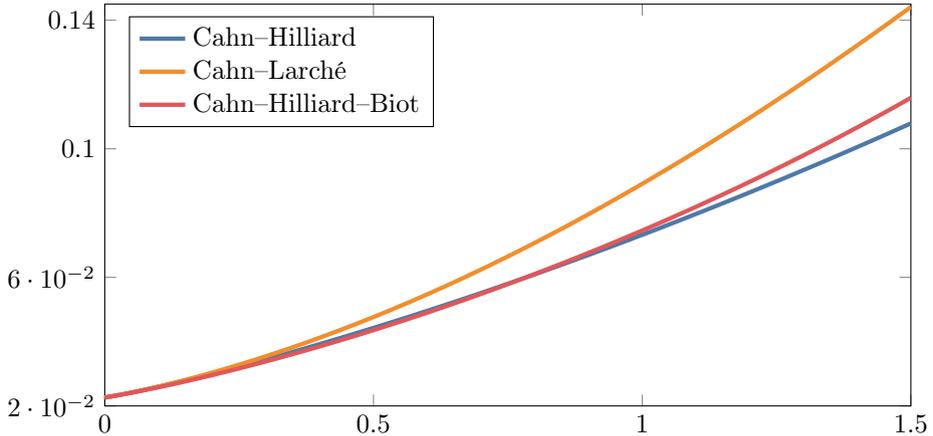

The simulated mass evolution for the different models is shown in \cref{Fig:Volume}. The first curve, represented in blue, corresponds to the Cahn--Hilliard equation. This basic model serves as a reference point for understanding tumor growth dynamics. As expected, the tumor mass exhibits an almost linear increase over time, aligning with conventional growth expectations; see \cite[Section 5.3]{fritz2022time}. The slower initial growth is due to the unstructured initial condition. The second curve, depicted in orange, corresponds to the Cahn--Larch\'e equation, which introduces elastic effects into the model. This alteration significantly affects the growth behavior of the tumor. In particular, we observe a substantial increase in tumor mass over the simulated time frame and the growth does not follow a linear trajectory. Instead, it exhibits behavior similar to a quadratic growth pattern. This observation underscores the influence of elastic forces on the dynamics of tumor growth. 
The third curve, shown in red, extends the model by incorporating flow effects through the Biot model, in addition to elasticity. This combined approach results in an intermediate growth pattern. Although not as pronounced as the Cahn--Larch\'e model, the tumor's growth still deviates from linearity. In particular, the inclusion of flow effects appears to moderate tumor growth, suggesting that flow phenomena have a restraining influence on tumor mass expansion.

\begin{figure}[H] \begin{center}
\begin{tabular}{cM{.195\textwidth}M{.195\textwidth}M{.195\textwidth}M{.195\textwidth}}
&
$t=0.1$&$t=0.5$&$t=1.0$&$t=1.5$ \\
CH&
\includegraphics[clip,trim={4cm 4cm 4cm 4cm},width=.22\textwidth]{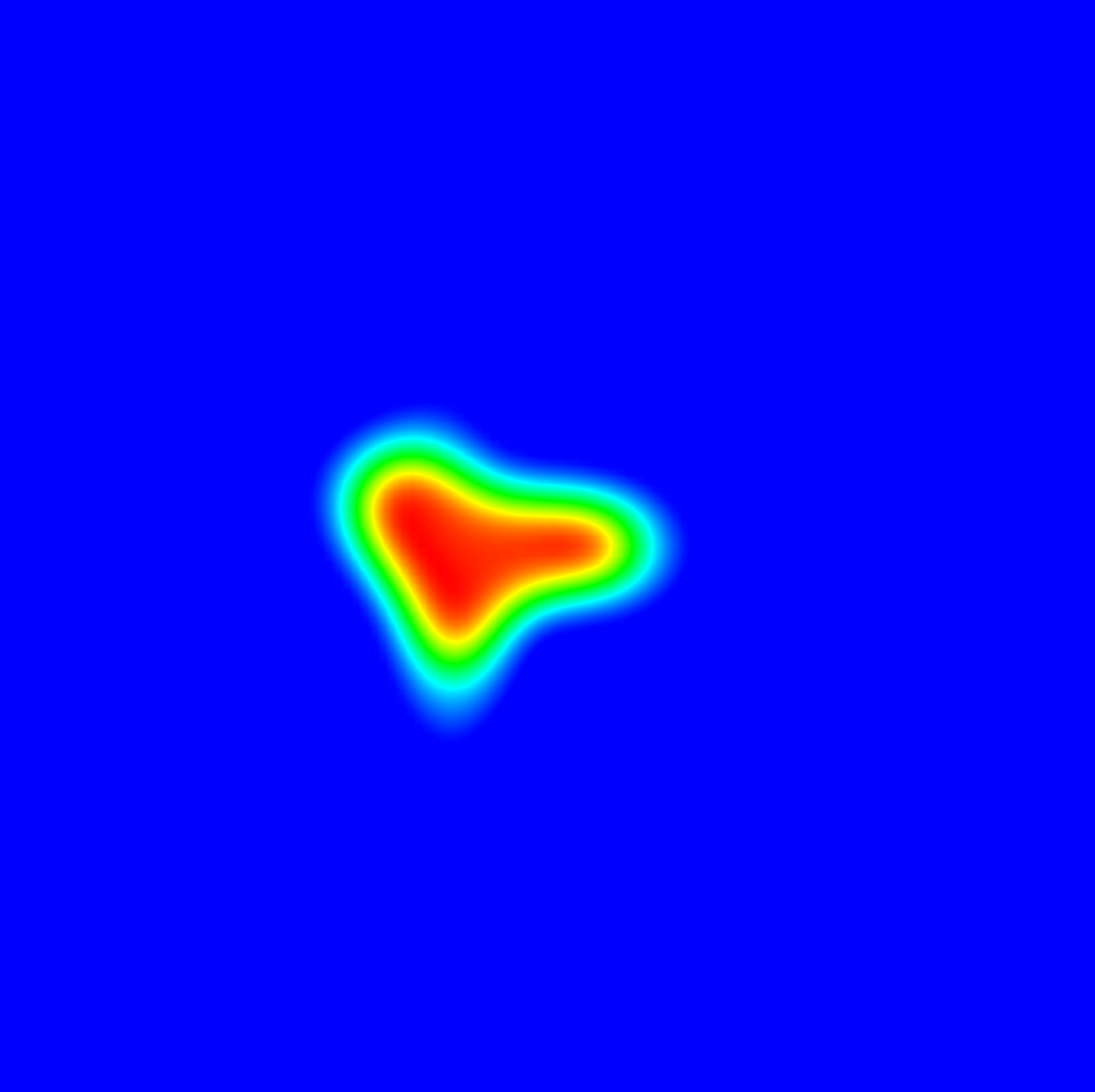}&
\includegraphics[clip,trim={4cm 4cm 4cm 4cm},width=.22\textwidth]{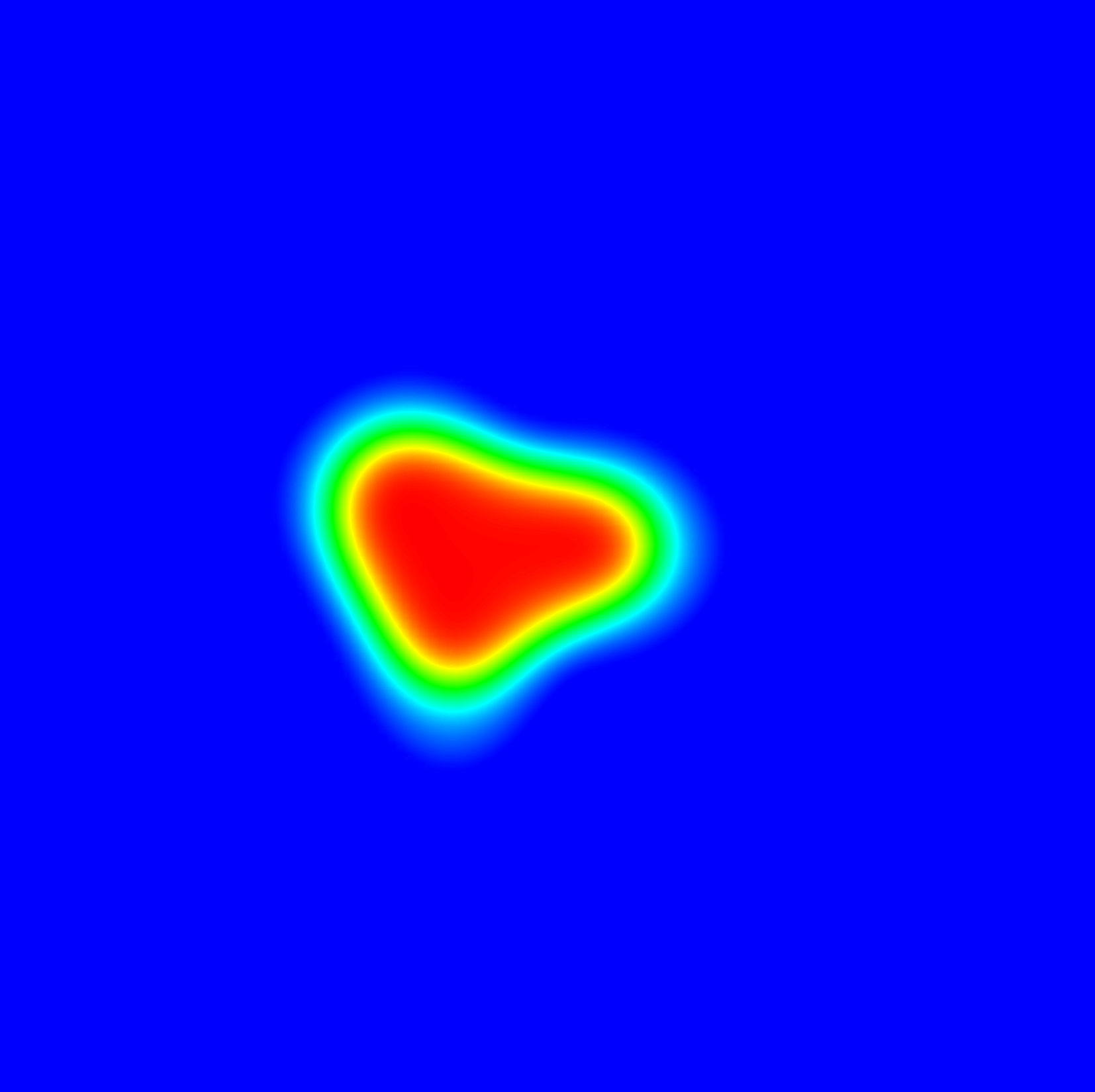}&
\includegraphics[clip,trim={4cm 4cm 4cm 4cm},width=.22\textwidth]{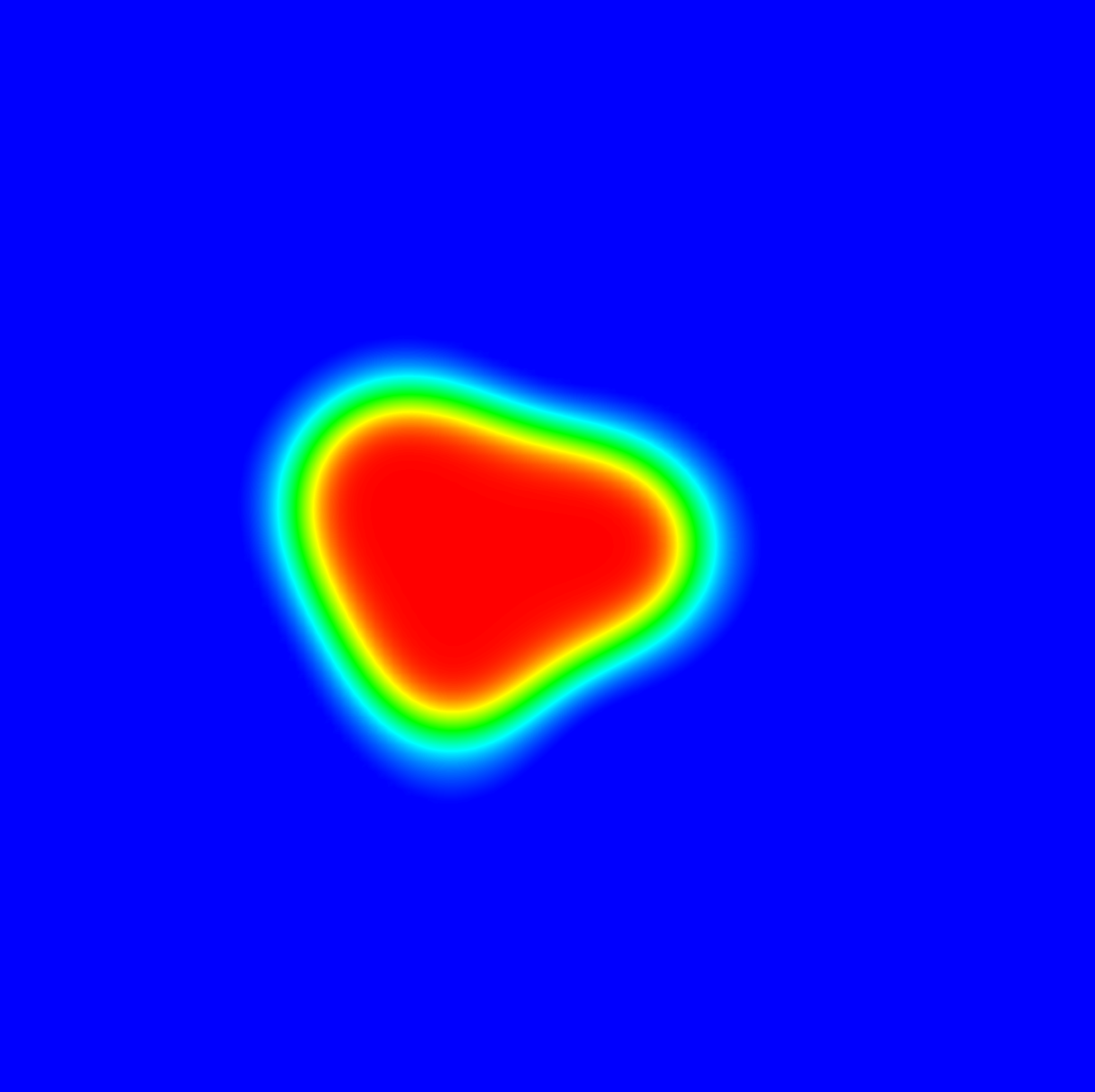}&
\includegraphics[clip,trim={4cm 4cm 4cm 4cm},width=.22\textwidth]{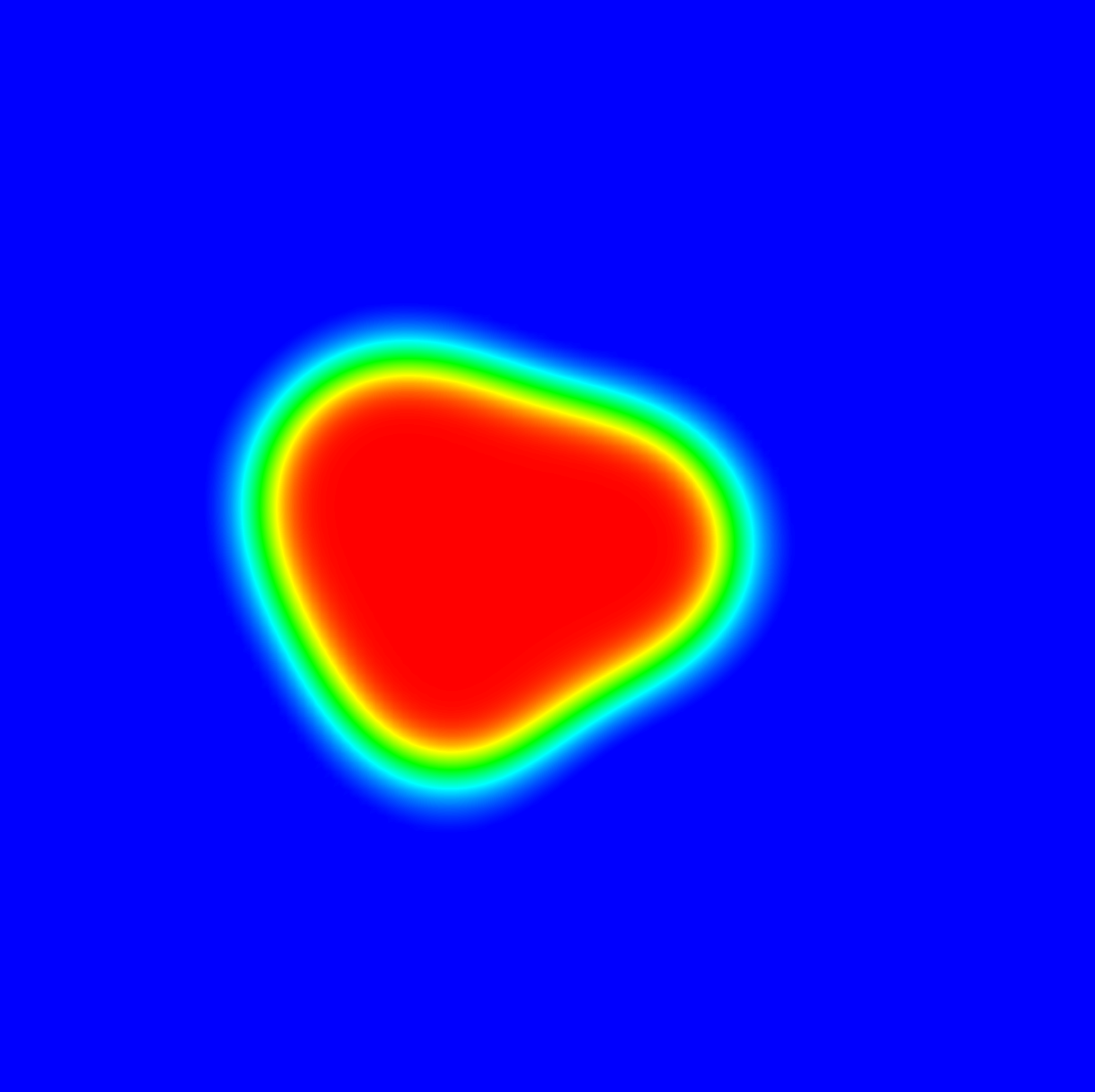}\\[-.1cm]
CL&
\includegraphics[clip,trim={4cm 4cm 4cm 4cm},width=.22\textwidth]{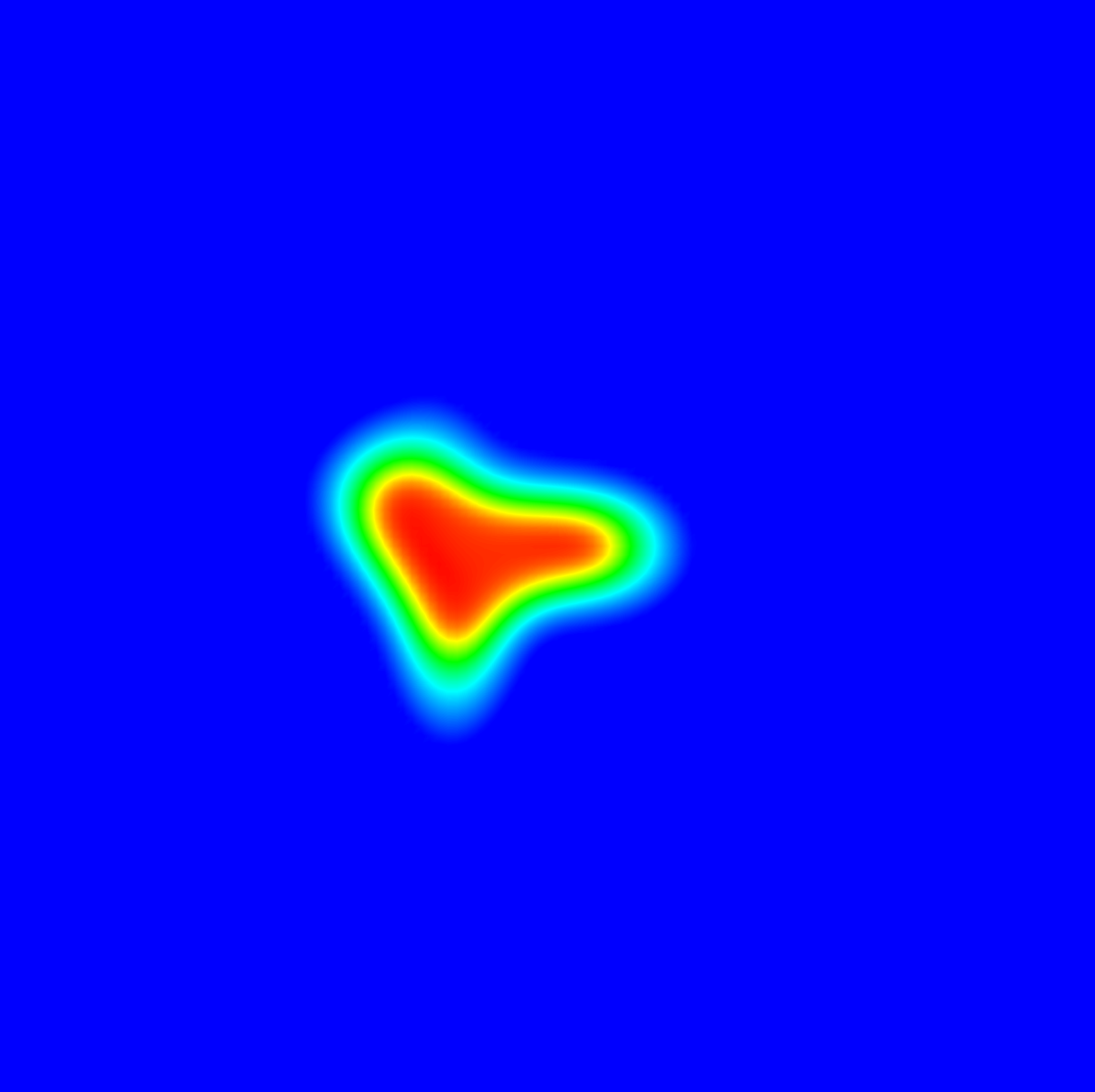}&
\includegraphics[clip,trim={4cm 4cm 4cm 4cm},width=.22\textwidth]{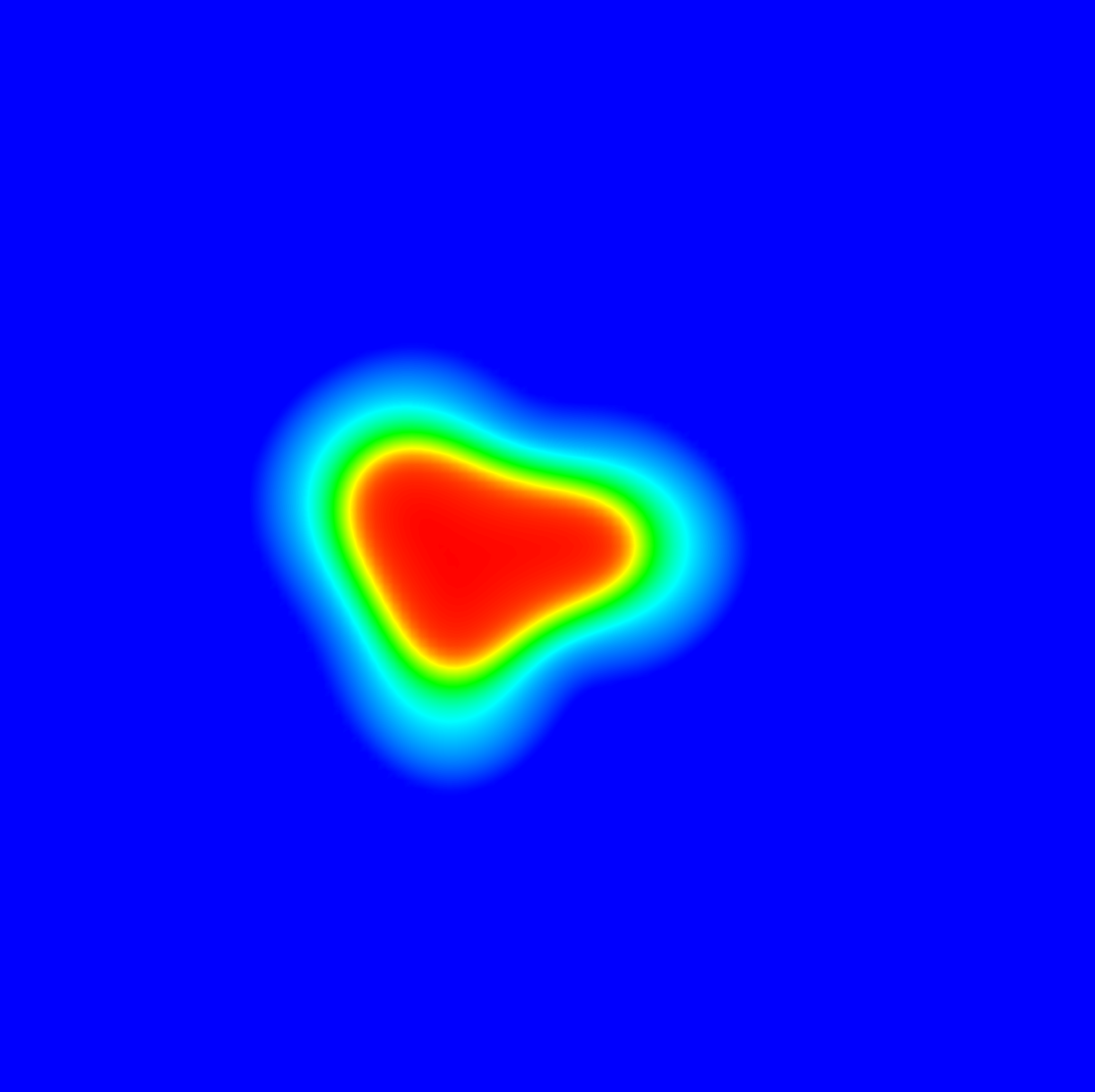}&
\includegraphics[clip,trim={4cm 4cm 4cm 4cm},width=.22\textwidth]{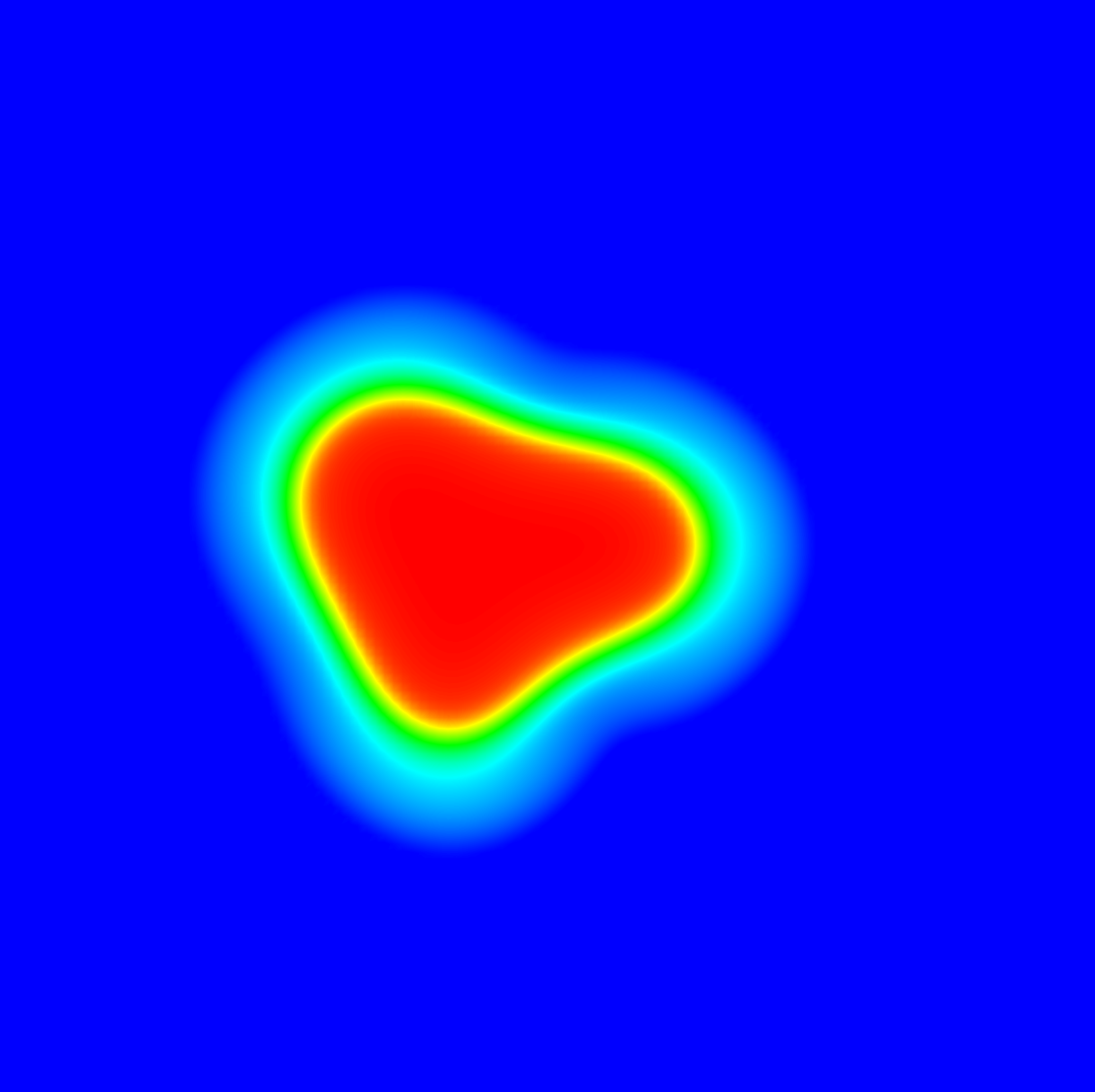}&
\includegraphics[clip,trim={4cm 4cm 4cm 4cm},width=.22\textwidth]{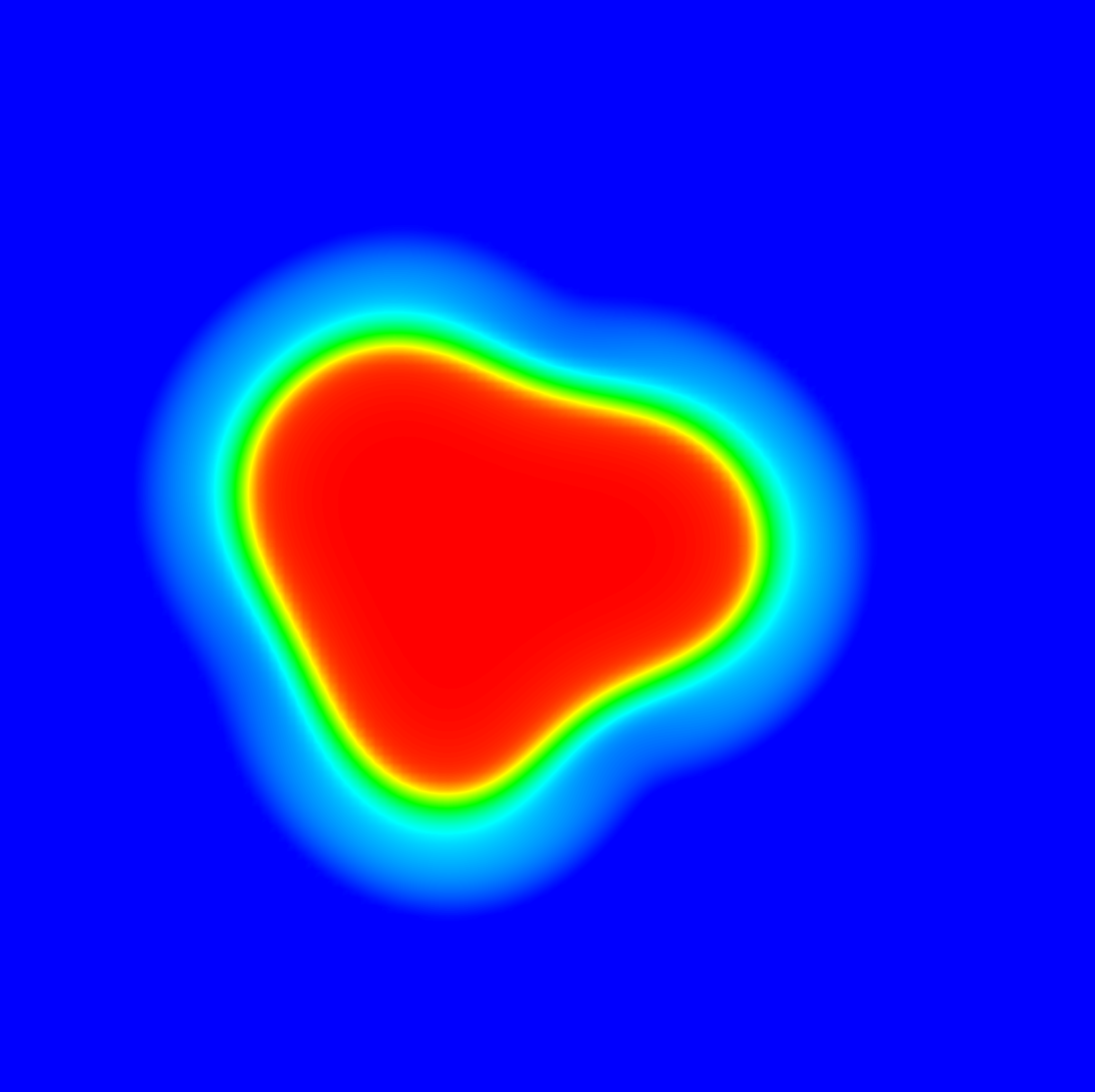}\\[-.1cm]
\!\!\!\!\!CHB\!\!\!\!\!&
\includegraphics[clip,trim={4cm 4cm 4cm 4cm},width=.22\textwidth]{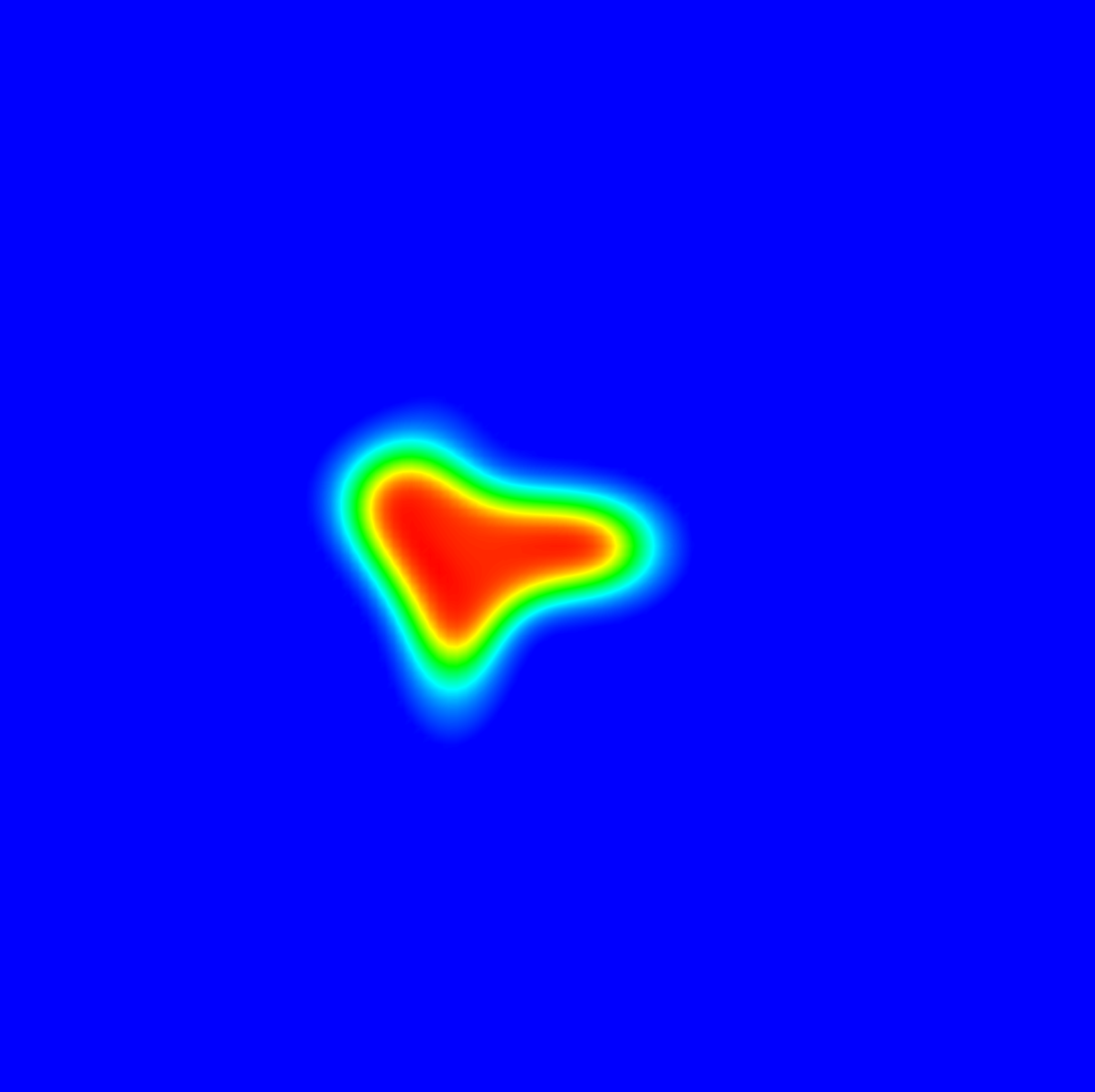}&
\includegraphics[clip,trim={4cm 4cm 4cm 4cm},width=.22\textwidth]{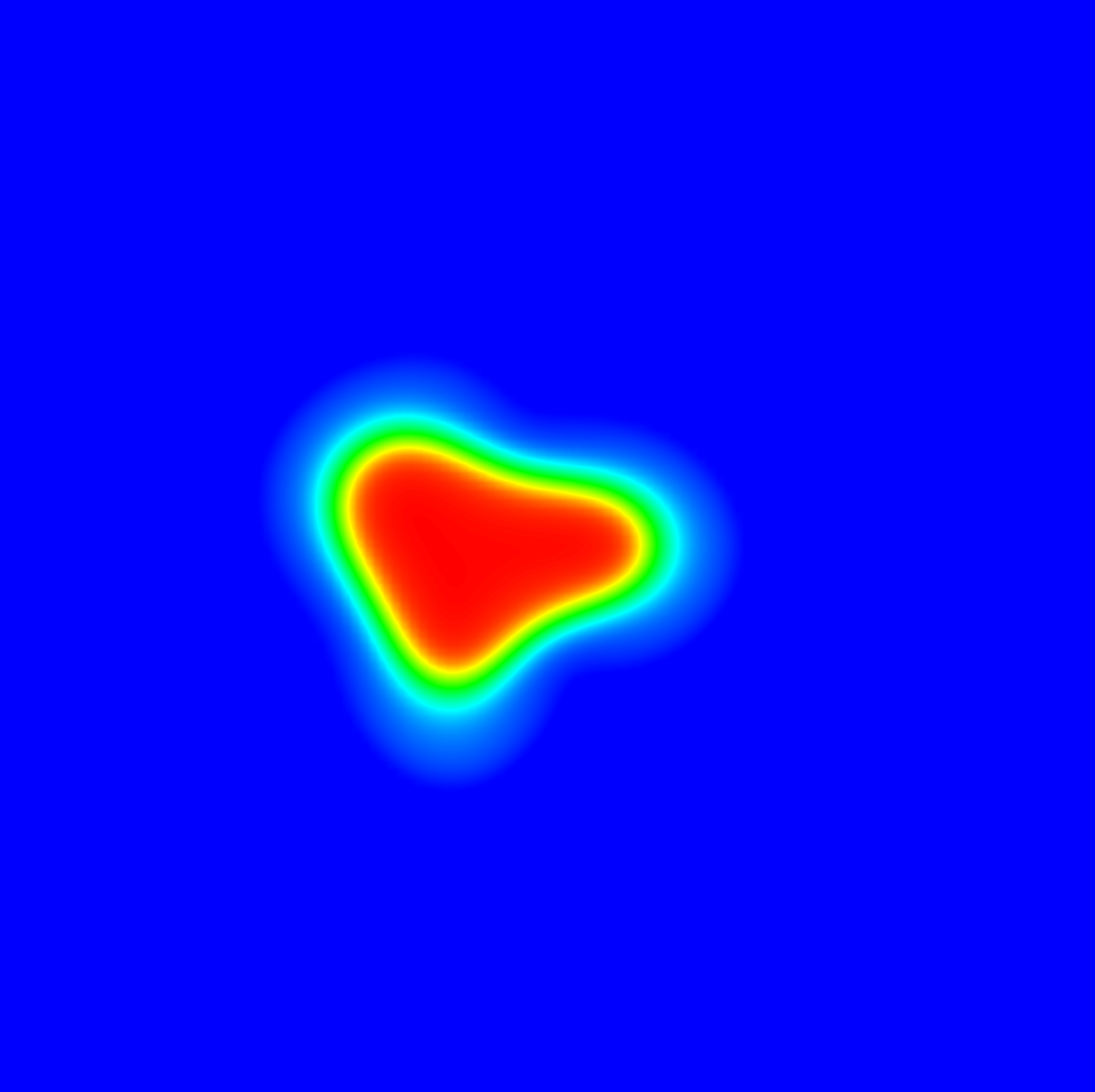}&
\includegraphics[clip,trim={4cm 4cm 4cm 4cm},width=.22\textwidth]{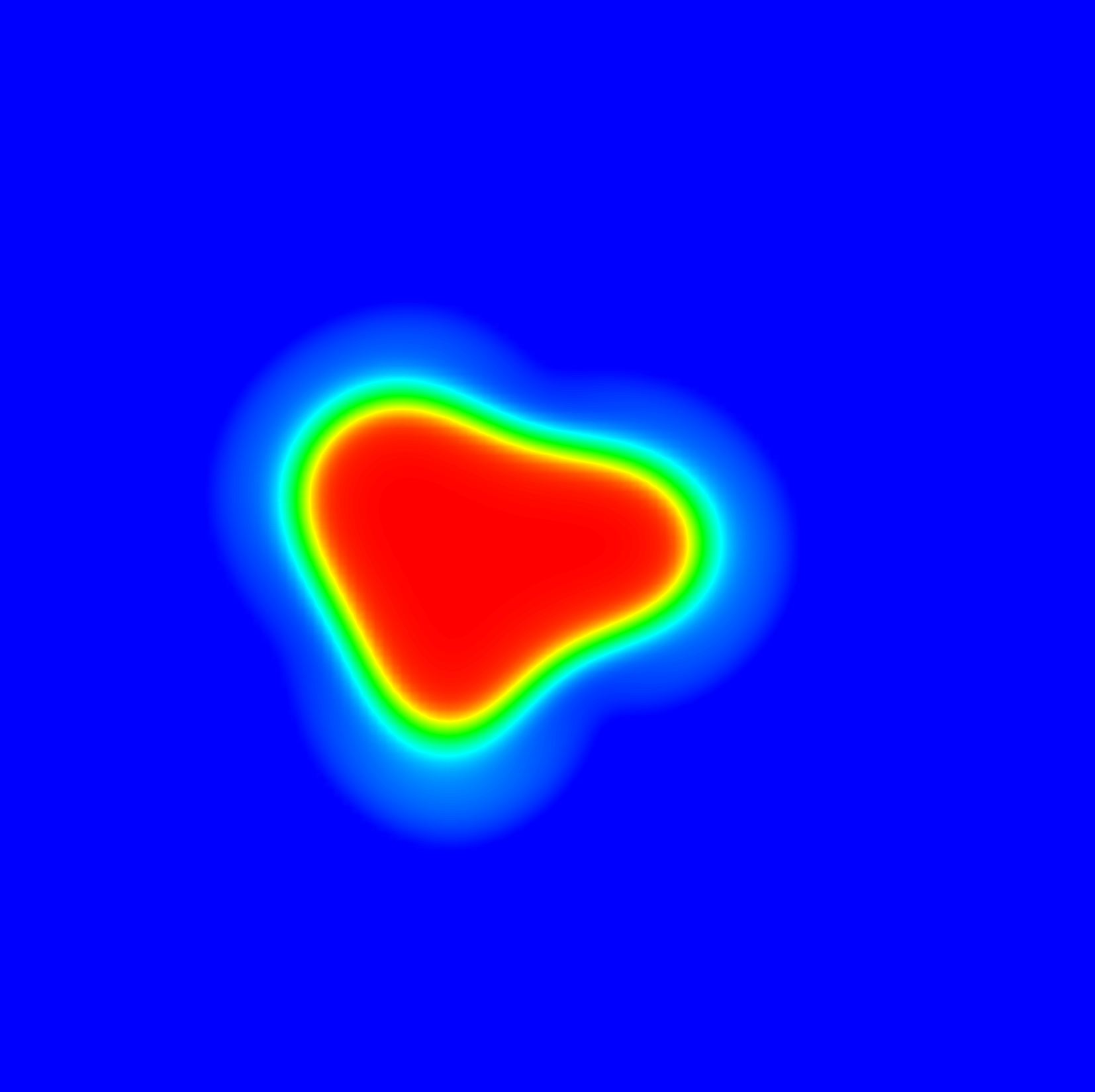}&
\includegraphics[clip,trim={4cm 4cm 4cm 4cm},width=.22\textwidth]{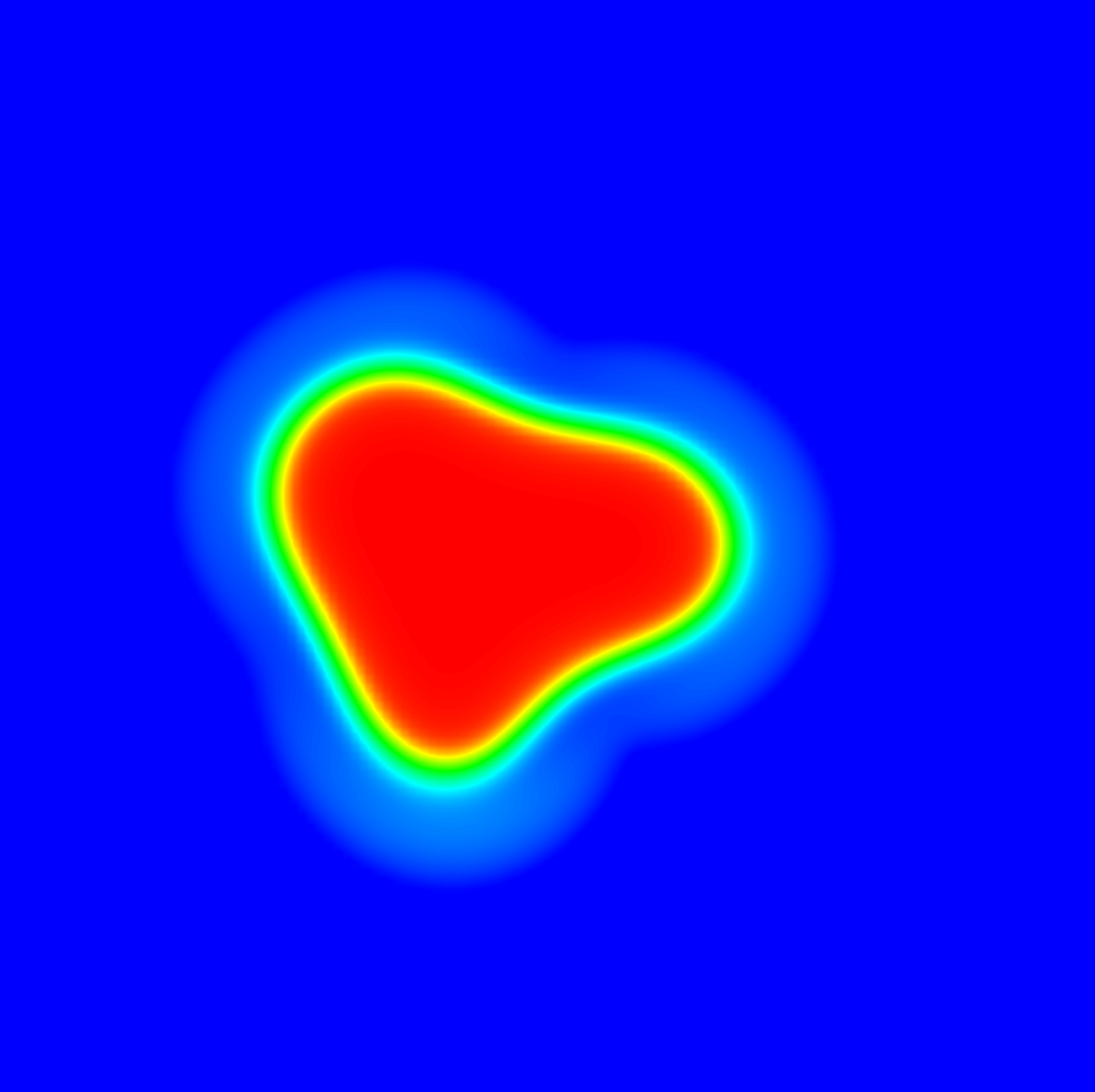}
\end{tabular} \\
\hfill\includegraphics[width=.88\textwidth]{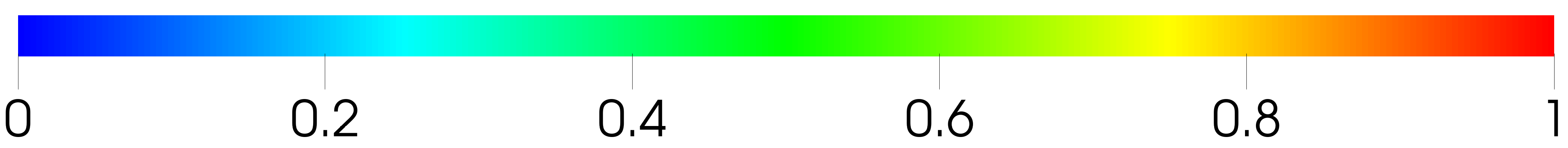}
\caption{\label{Fig:TumorMass}Evolution of the tumor volume fraction $\phi(t,x)$ over time in the domain $\Omega$}
\end{center} 
\end{figure}

\begin{figure}[H] \begin{center}
\begin{tabular}{cM{.2\textwidth}M{.2\textwidth}M{.2\textwidth}M{.2\textwidth}}
&$t=0.1$&$t=0.5$&$t=1$&$t=1.5$ \\
$\phi=0.5$&\includegraphics[width=.21\textwidth]{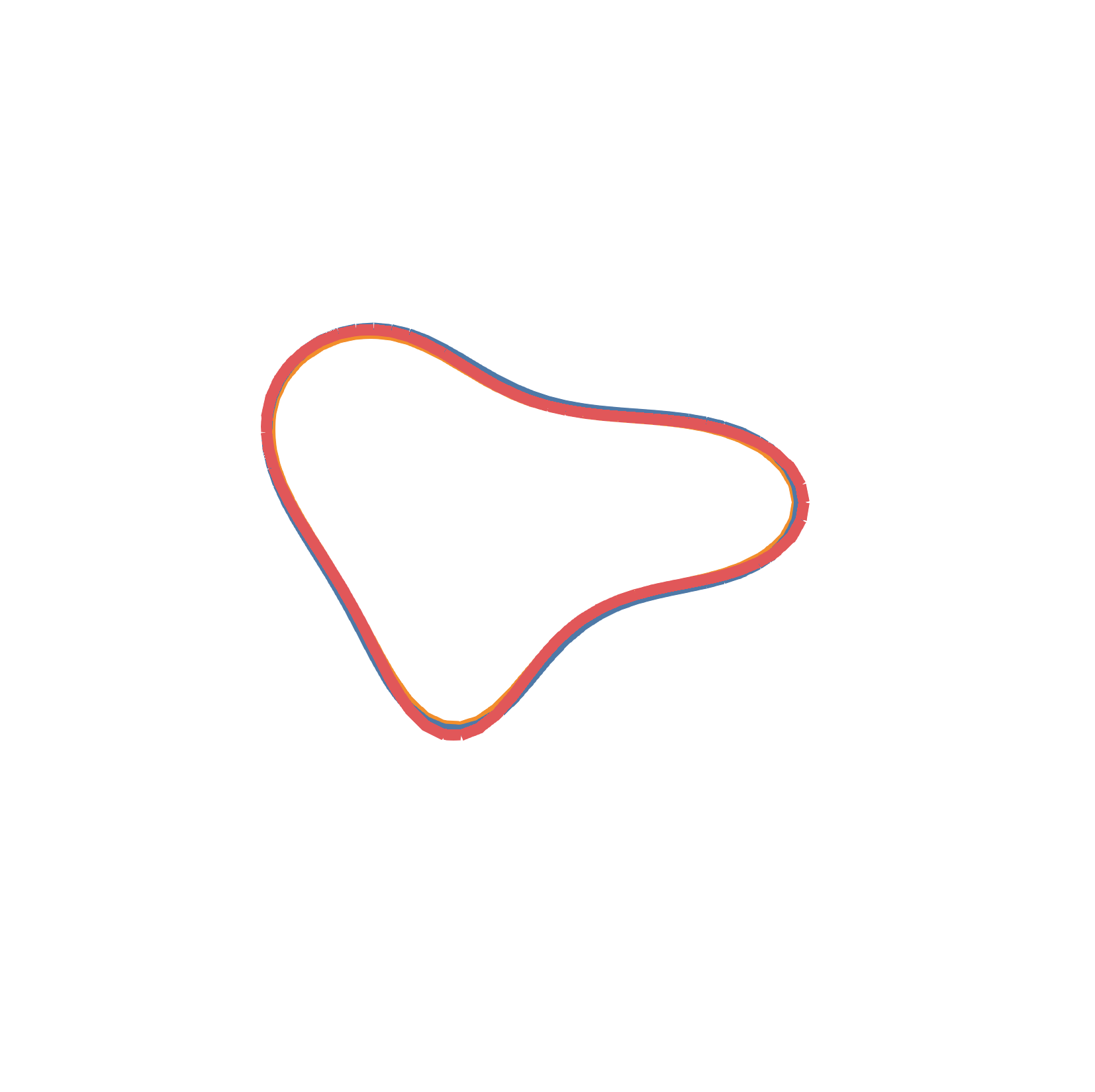}&
\includegraphics[width=.21\textwidth]{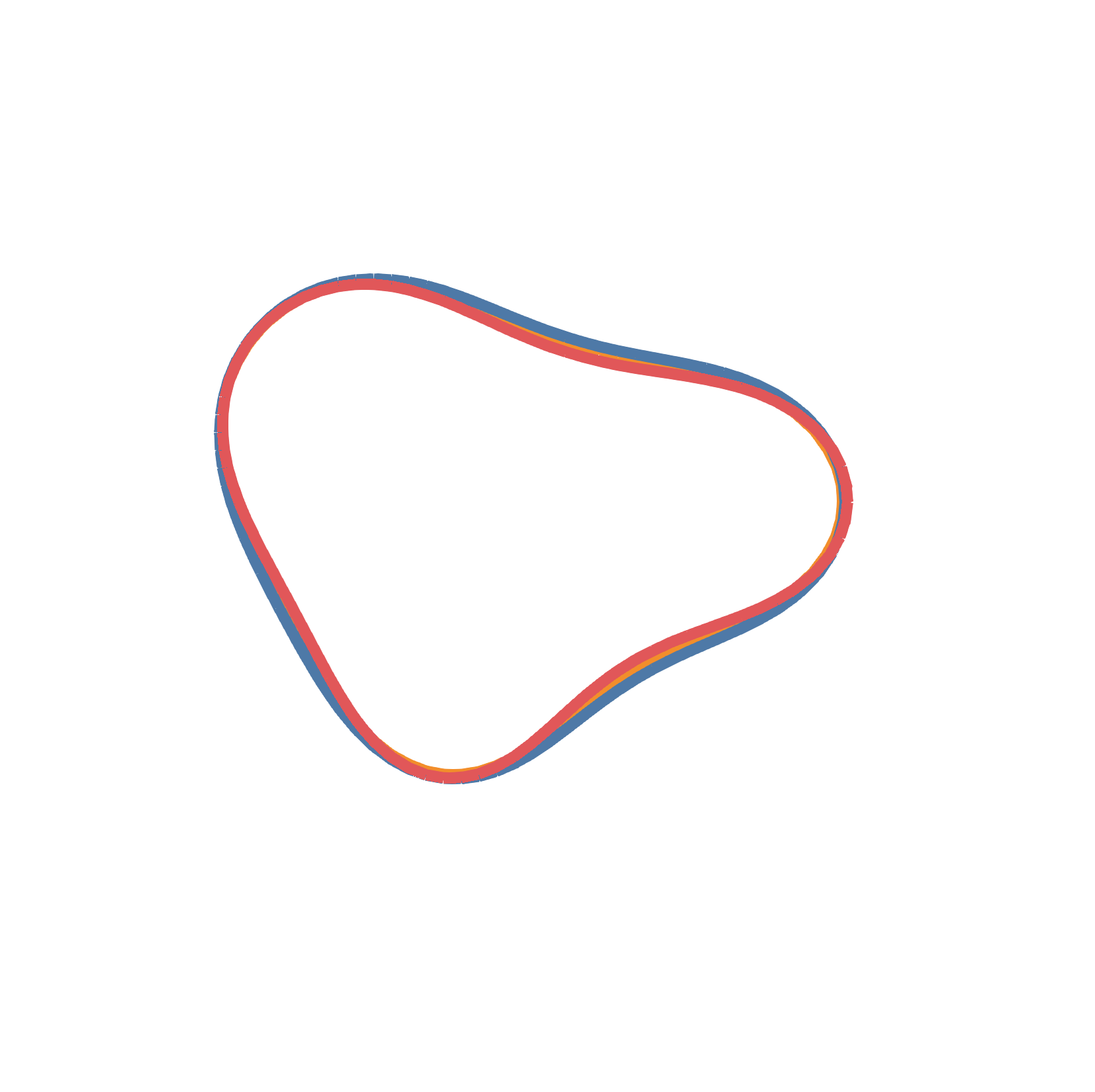}&
\includegraphics[width=.21\textwidth]{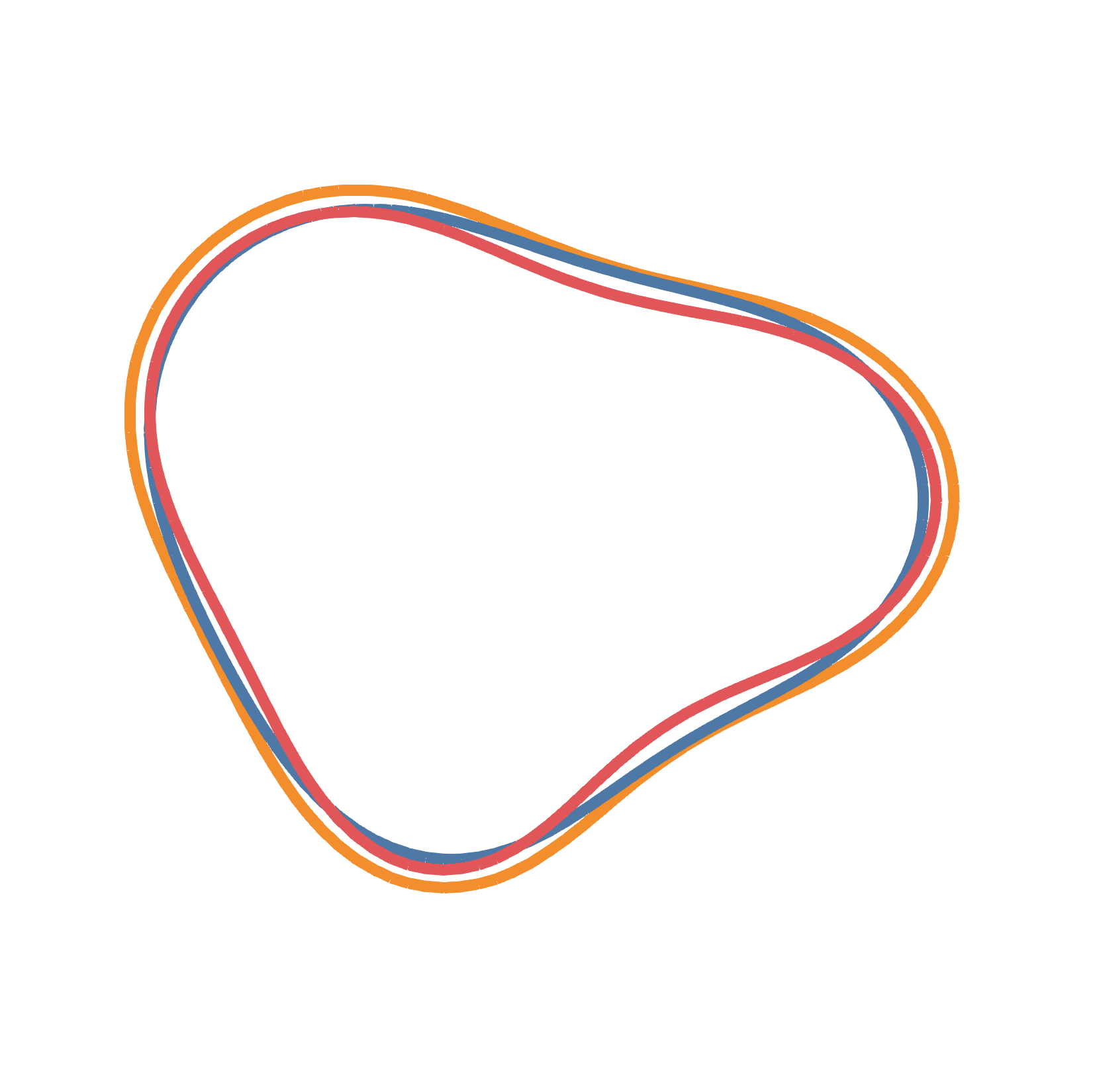}&
\includegraphics[width=.21\textwidth]{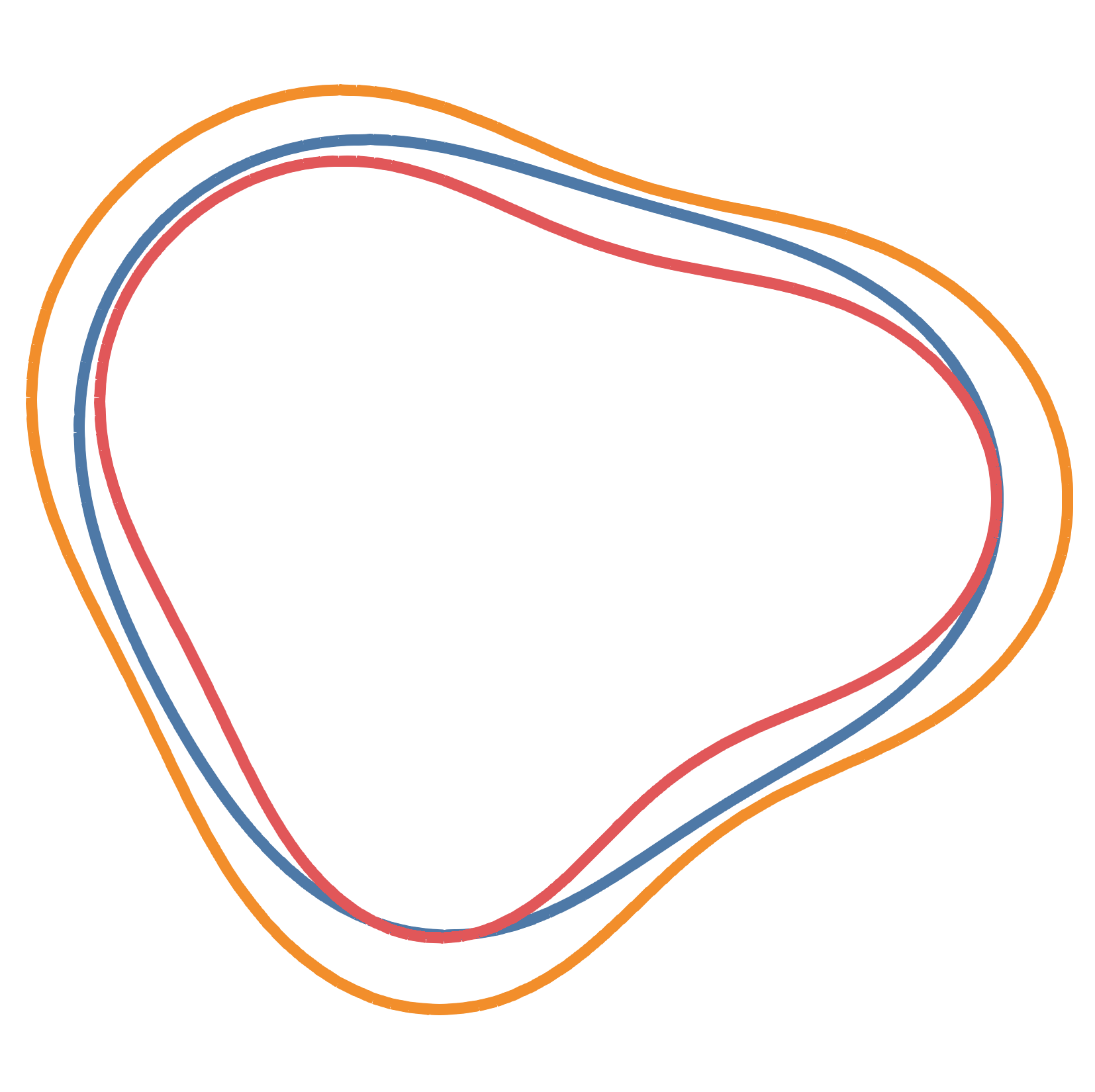} \\
$\phi=0.9$&\includegraphics[width=.21\textwidth]{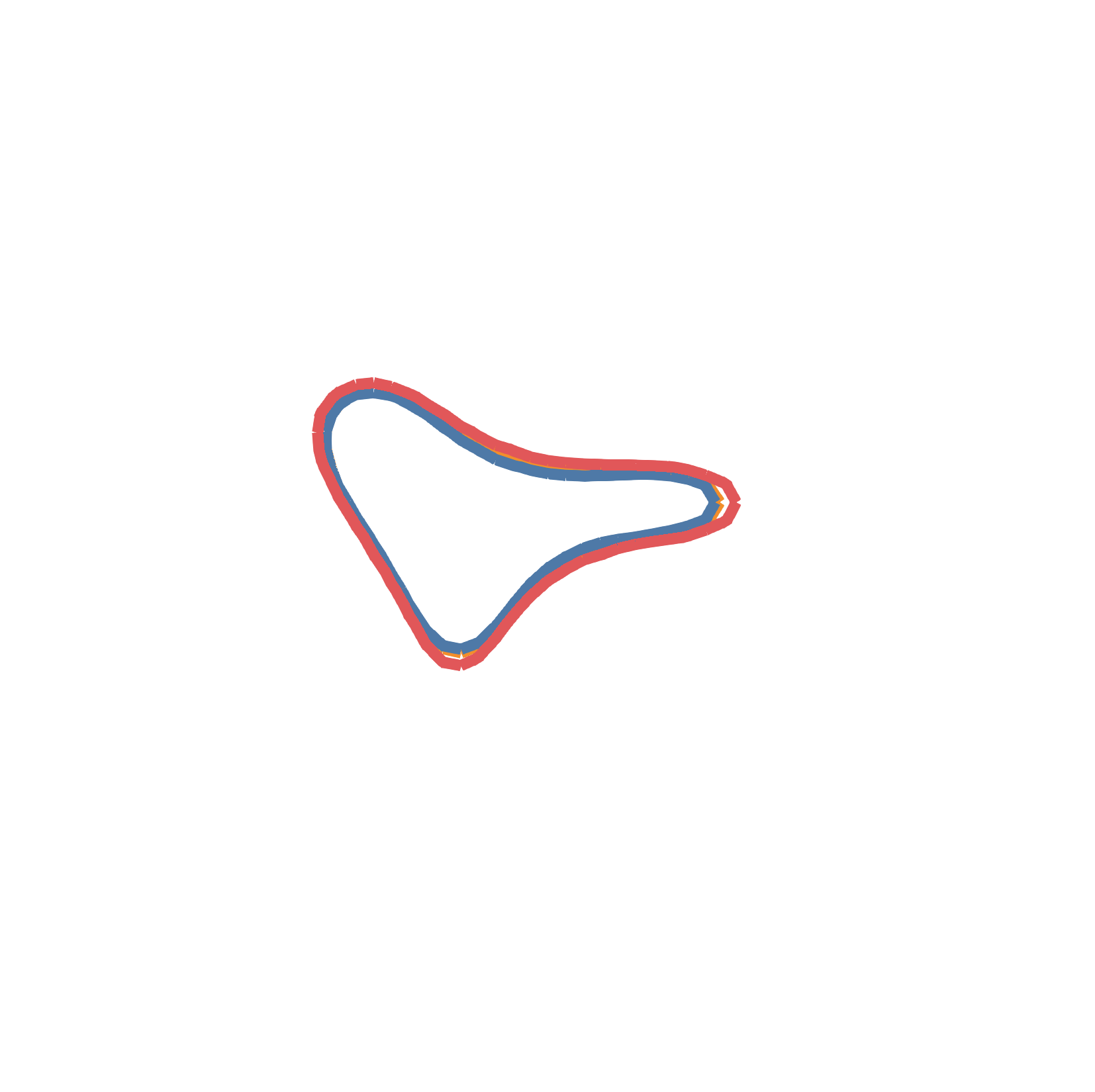}&
\includegraphics[width=.21\textwidth]{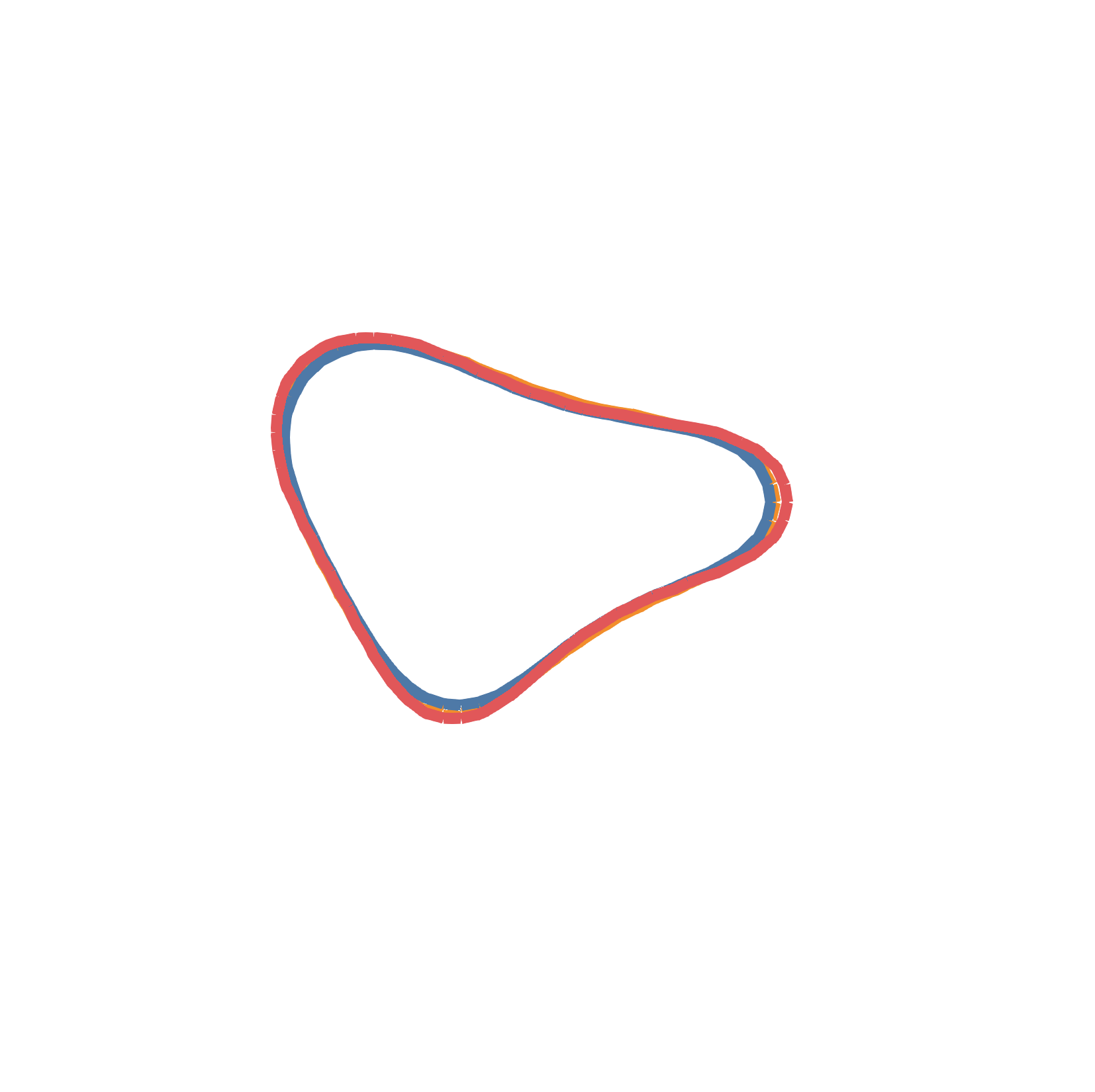}&
\includegraphics[width=.21\textwidth]{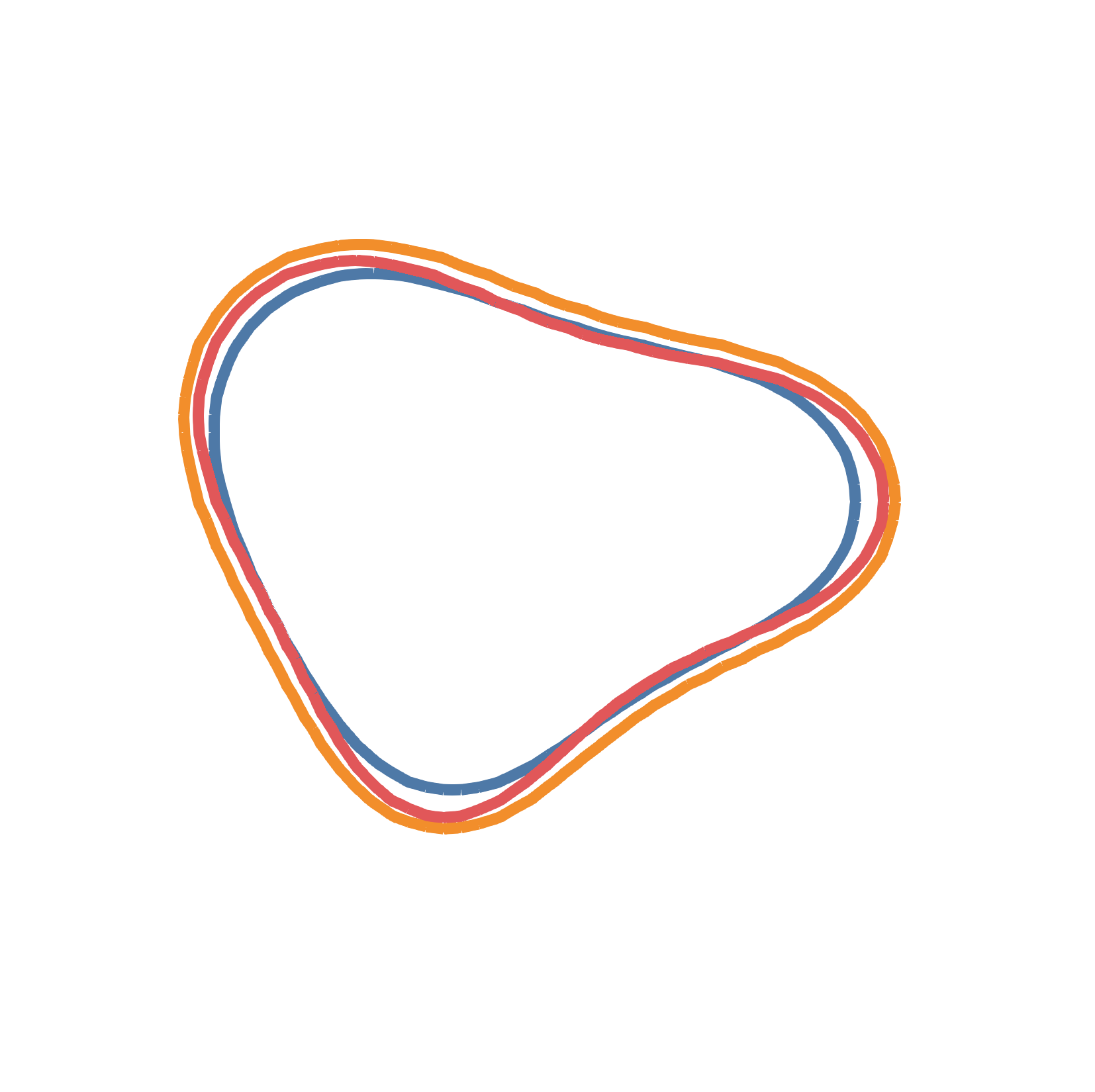}&
\includegraphics[width=.21\textwidth]{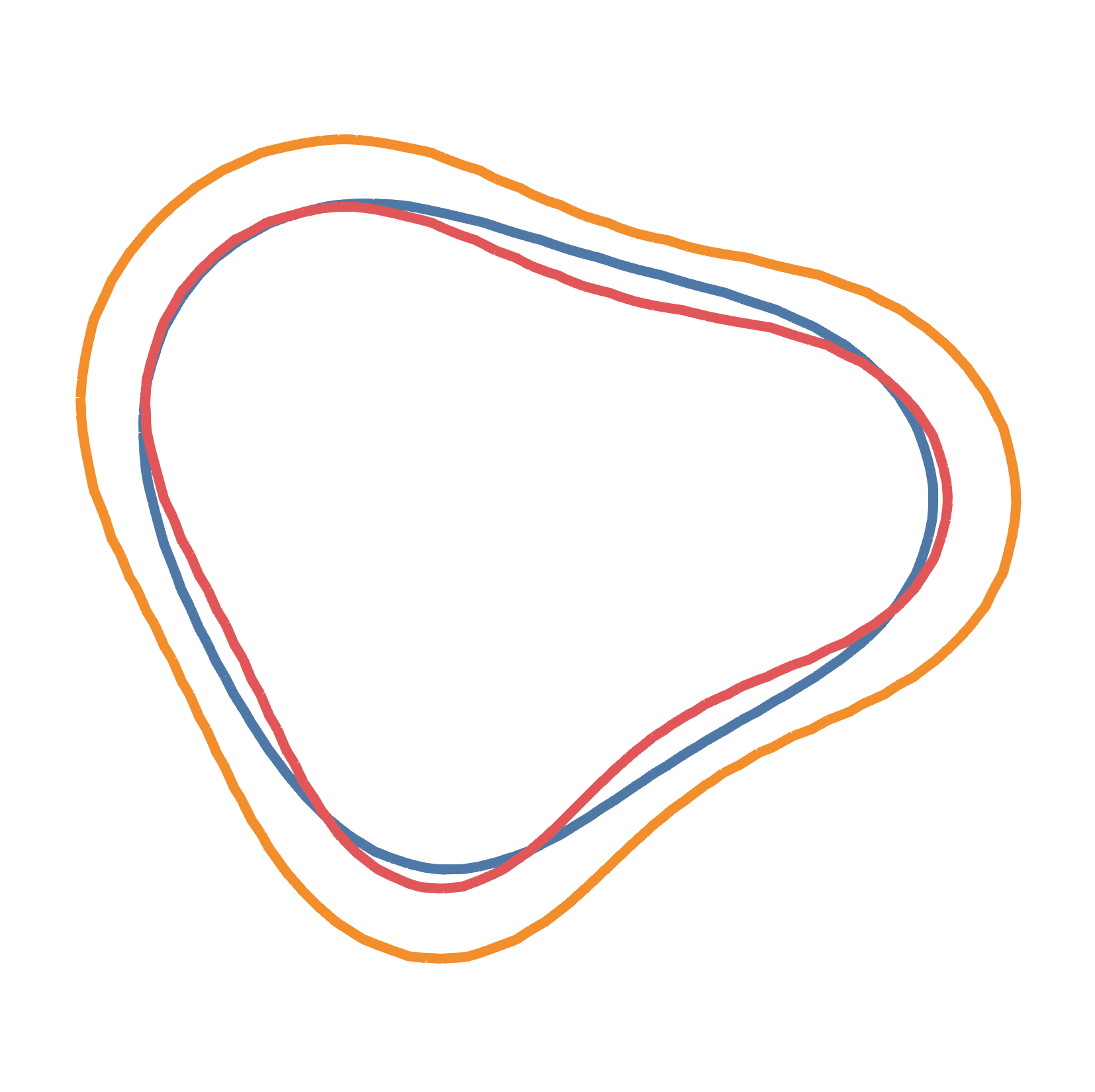} \\[-.5cm]
\end{tabular} 
\caption{\label{Fig:Contour}Comparison of the contour lines $\phi(t,x)=0.5$ (upper) and $\phi(t,x)=0.9$ (lower) at different times for the Cahn--Hilliard (\textcolor{color00}{$\boldsymbol{-}$}), Cahn--Larch\'e equation (\textcolor{color01}{$\boldsymbol{-}$}) and Cahn--Hilliard--Biot equations (\textcolor{color02}{$\boldsymbol{-}$}).}
\end{center}
\end{figure}

\cref{Fig:TumorMass} presents a visual representation of the dynamic evolution of the tumor over time, simulated using three mathematical models: the Cahn--Hilliard (CH), Cahn--Larch\'e (CL), and Cahn--Hilliard--Biot (CHB) equations. Each model provides distinct insights into how various factors influence the tumor's shape and symmetry.
In the upper row of \cref{Fig:TumorMass}, the evolution of the tumor is shown under the CH equation. As expected, the tumor tends to form a circular shape, which is consistent with the lack of additional physical factors in the standard model that might alter the geometry. This behavior is consistent with observations from previous studies, such as \cite{fritz2019unsteady}, where different flow models influence the irregularities in the shape of the tumor.
The middle row of \cref{Fig:TumorMass} displays the simulation results for the CL model. In this case, the tumor maintains its irregular structure, with visible dents that persist throughout the evolution. Notably, this model produces the largest interface width among the three models, indicating a broader transition zone between tumor and non-tumor tissue. This can be attributed to the choice of the interpolation function, $\pi$, which acts on the interface.
The bottom row of \cref{Fig:TumorMass} shows the evolution of the tumor in the more complex CHB model. The tumor's behavior is similar to that of the CL model, though with a smaller mass. Like the elasticity model, the CHB model also generates a wider interface between tumor and non-tumor tissue than the standard CH model, although the interface is narrower than in the CL model.
Further analysis of tumor shapes is provided by comparing the contour lines at $\phi = 0.5$ and $\phi = 0.9$ for the different models and at various time points in \cref{Fig:Contour}. The most significant differences are observed between the contour lines of the standard CH model and the CHB system, particularly in the irregular dented regions. While the standard CH model attempts to smooth out the dents by forming a more circular shape, the CHB model preserves these irregularities.
We observe that the different models affect the evolving shape and symmetry of the tumor. The inclusion of elasticity and other complex factors leads to notable variations in tumor morphology.

\begin{figure}[H] \begin{center}
\begin{tabular}{cM{.265\textwidth}M{.265\textwidth}M{.265\textwidth}}
&
$t=0.5$&$t=1.0$&$t=1.5$ \\
CL&
\includegraphics[width=.29\textwidth]{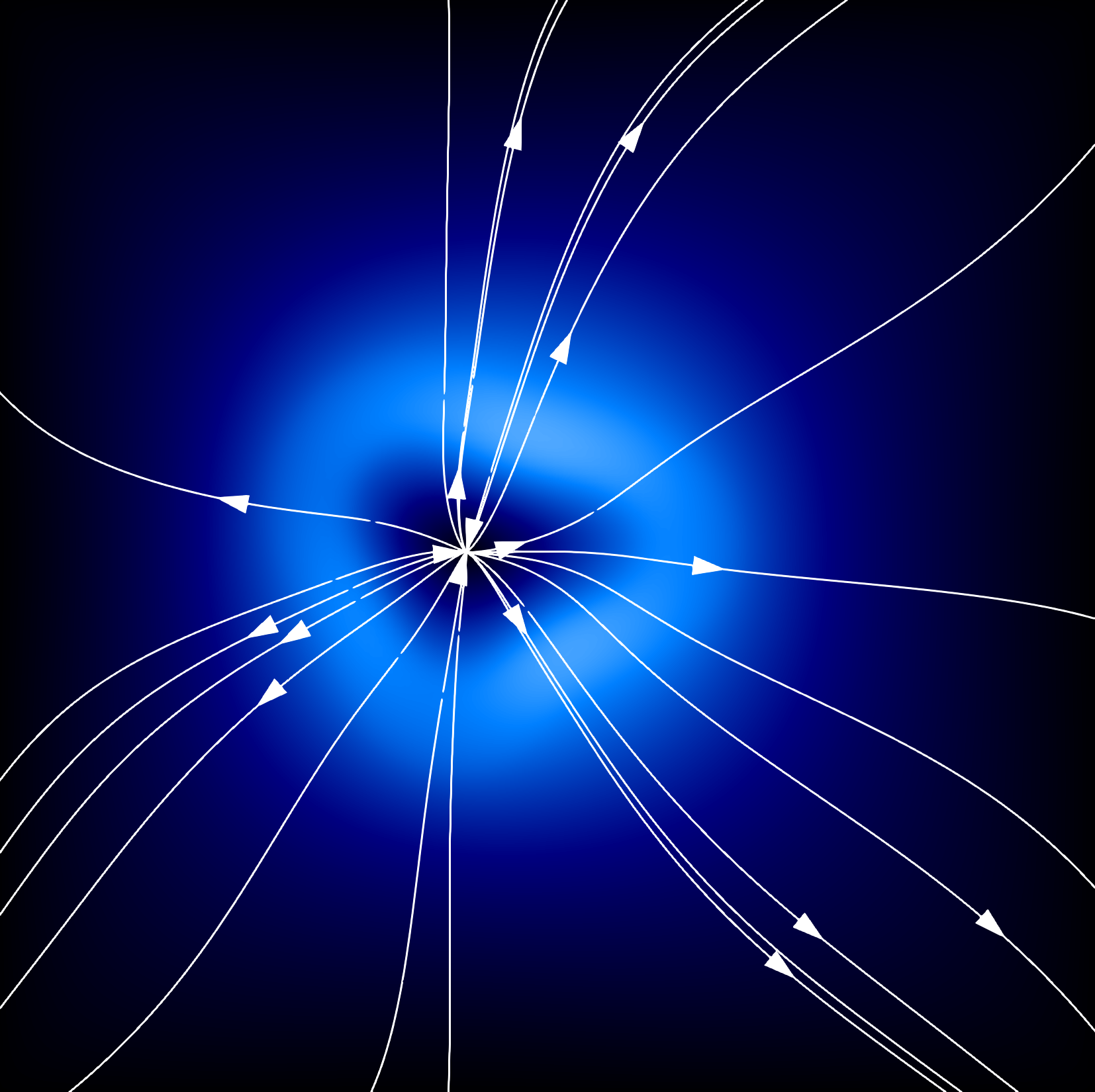}&
\includegraphics[width=.29\textwidth]{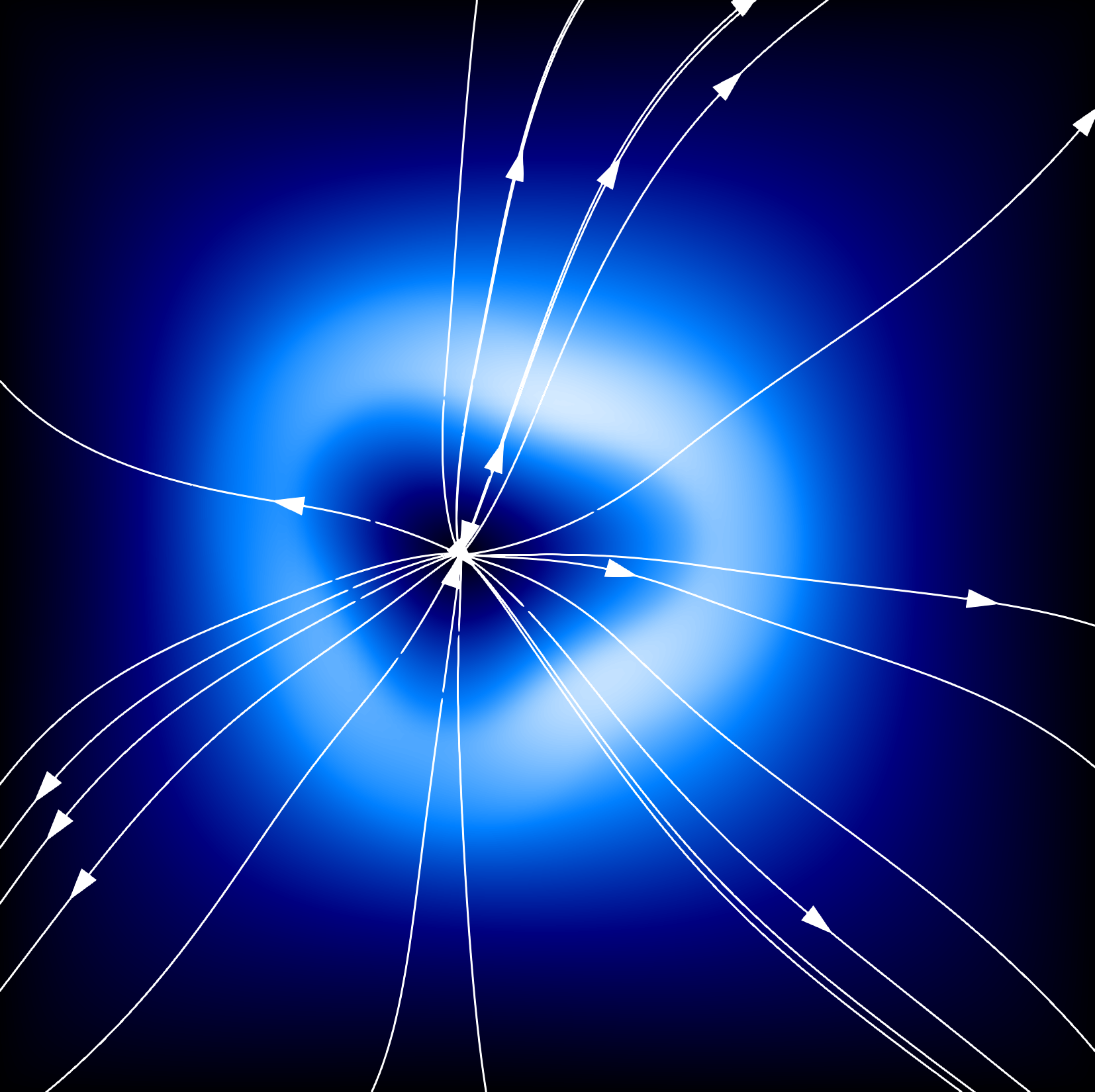}&
\includegraphics[width=.29\textwidth]{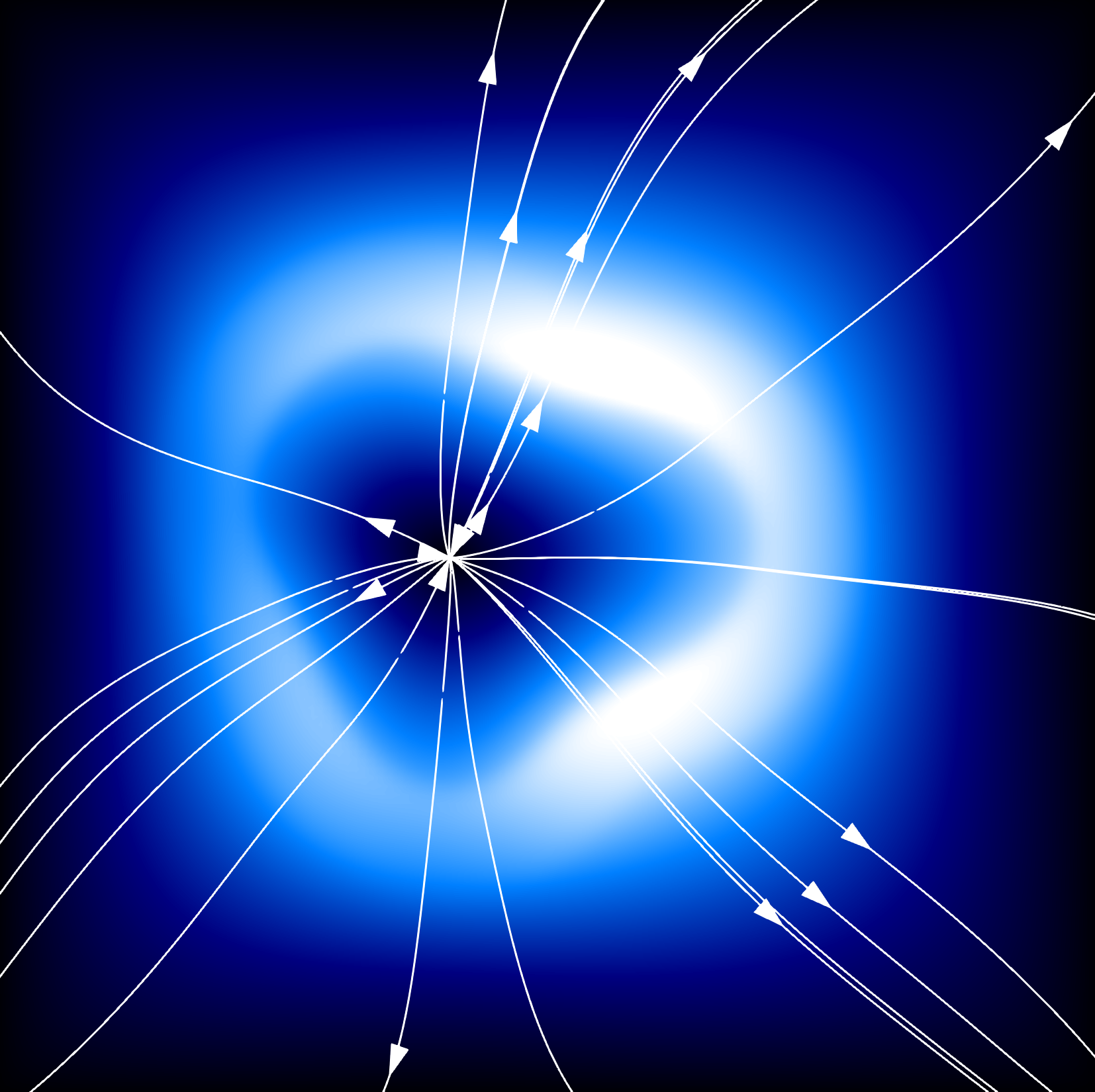}\\[-.1cm]
\!\!\!\!\!CHB\!\!\!\!\!&
\includegraphics[width=.29\textwidth]{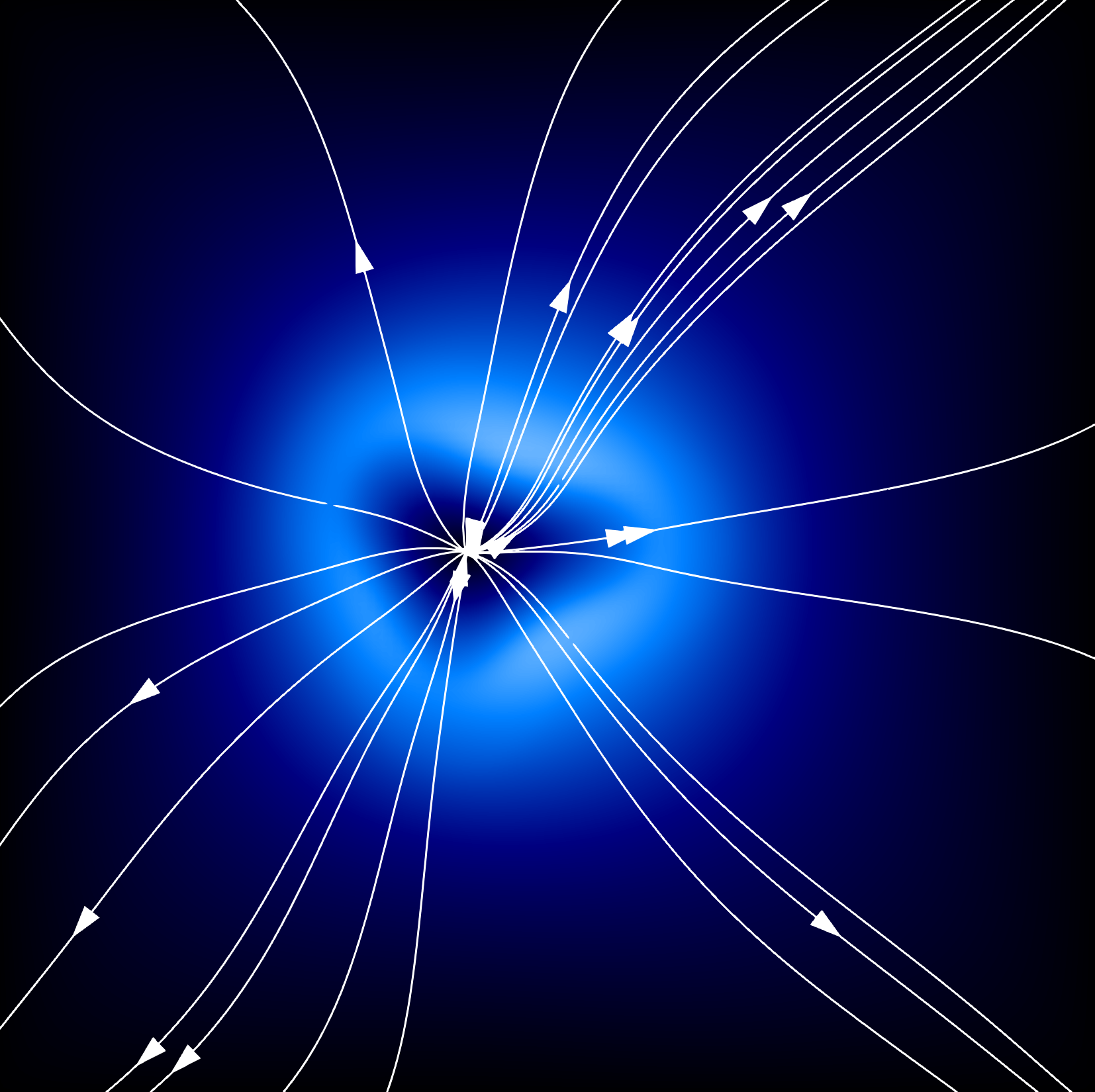}&
\includegraphics[width=.29\textwidth]{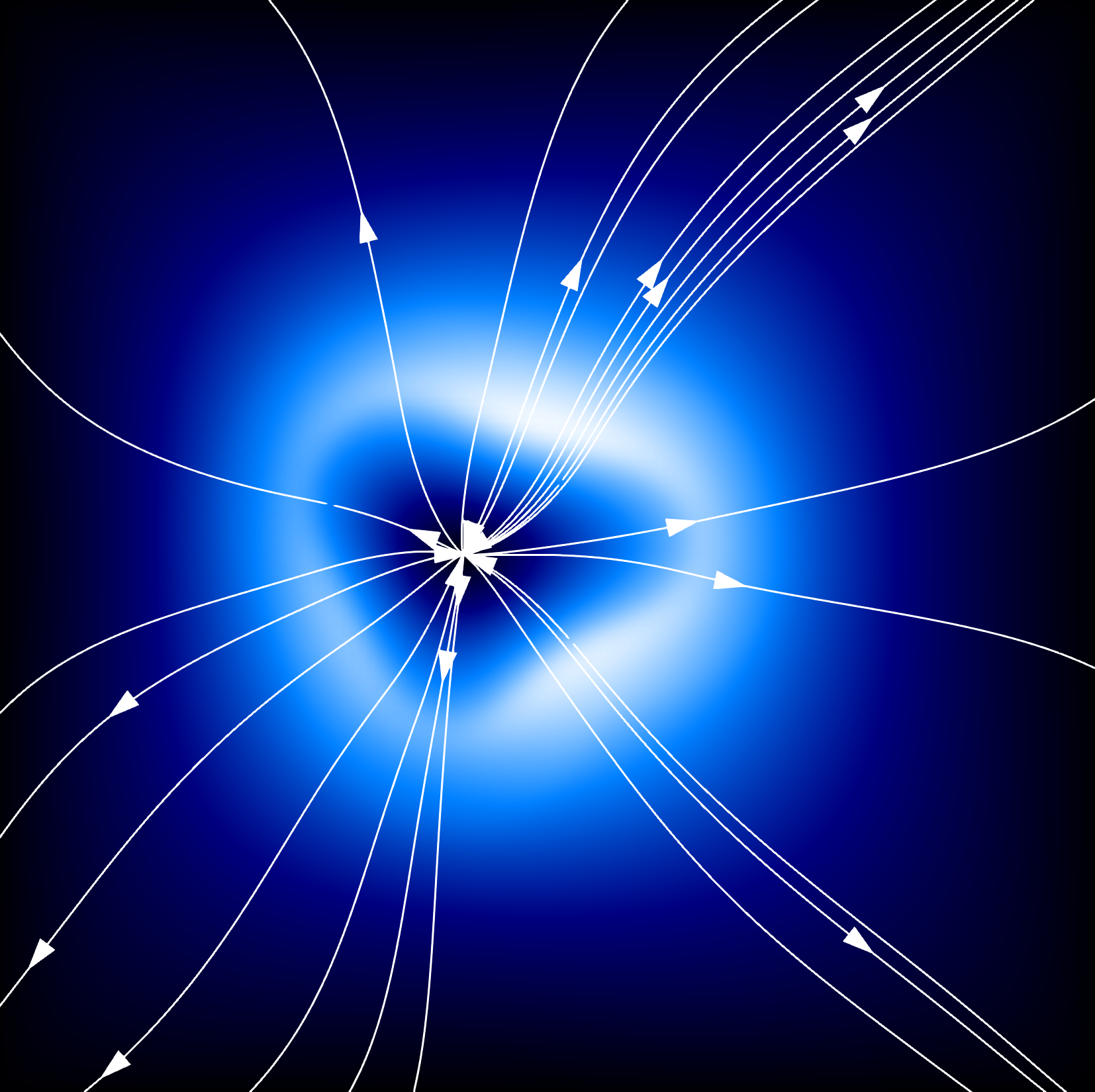}&
\includegraphics[width=.29\textwidth]{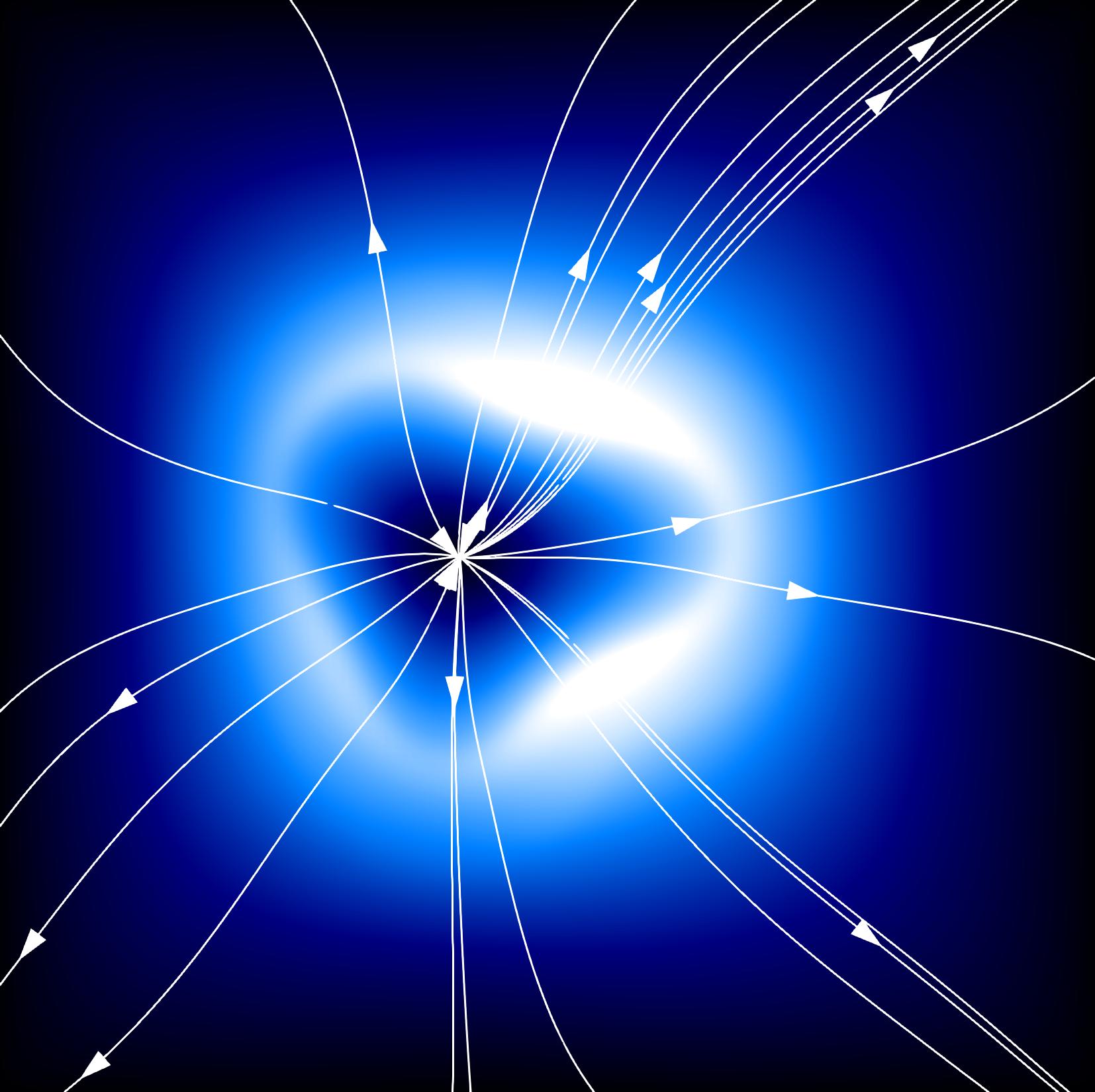}
\end{tabular} \\
\hfill\includegraphics[width=.9\textwidth]{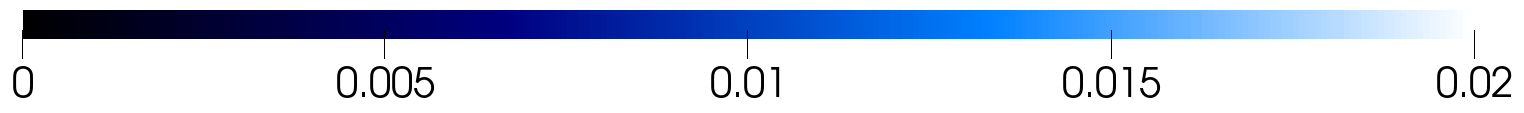}
\caption{\label{Fig:Deformation}Evolution of the deformation $u(t,x)$.}
\end{center}
\end{figure}

\begin{figure}[H] \begin{center}
\begin{tabular}{cM{.265\textwidth}M{.265\textwidth}M{.265\textwidth}}
&
$t=0.5$&$t=1.0$&$t=1.5$ \\
CHB&
\includegraphics[width=.29\textwidth]{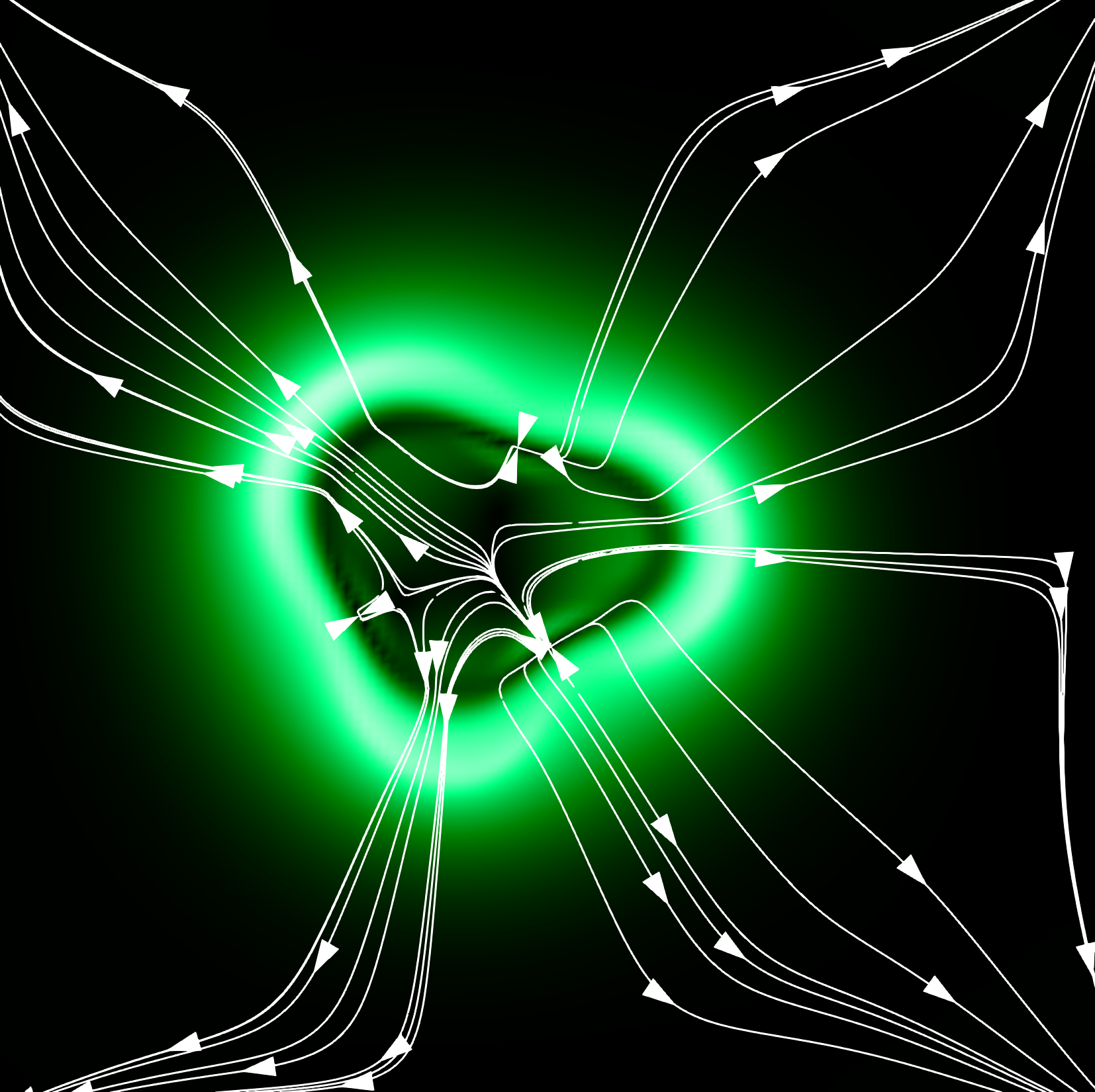}&
\includegraphics[width=.29\textwidth]{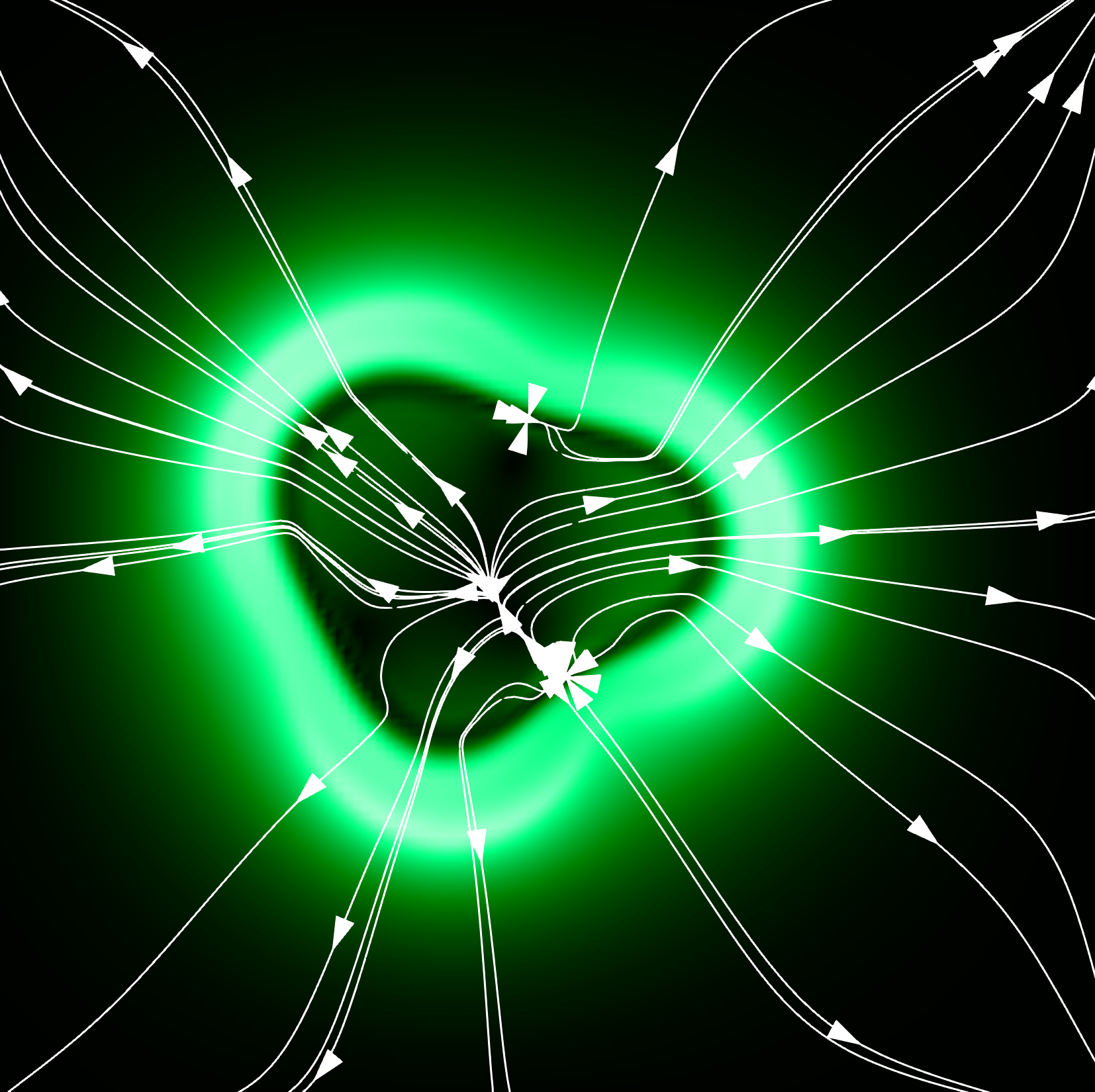}&
\includegraphics[width=.29\textwidth]{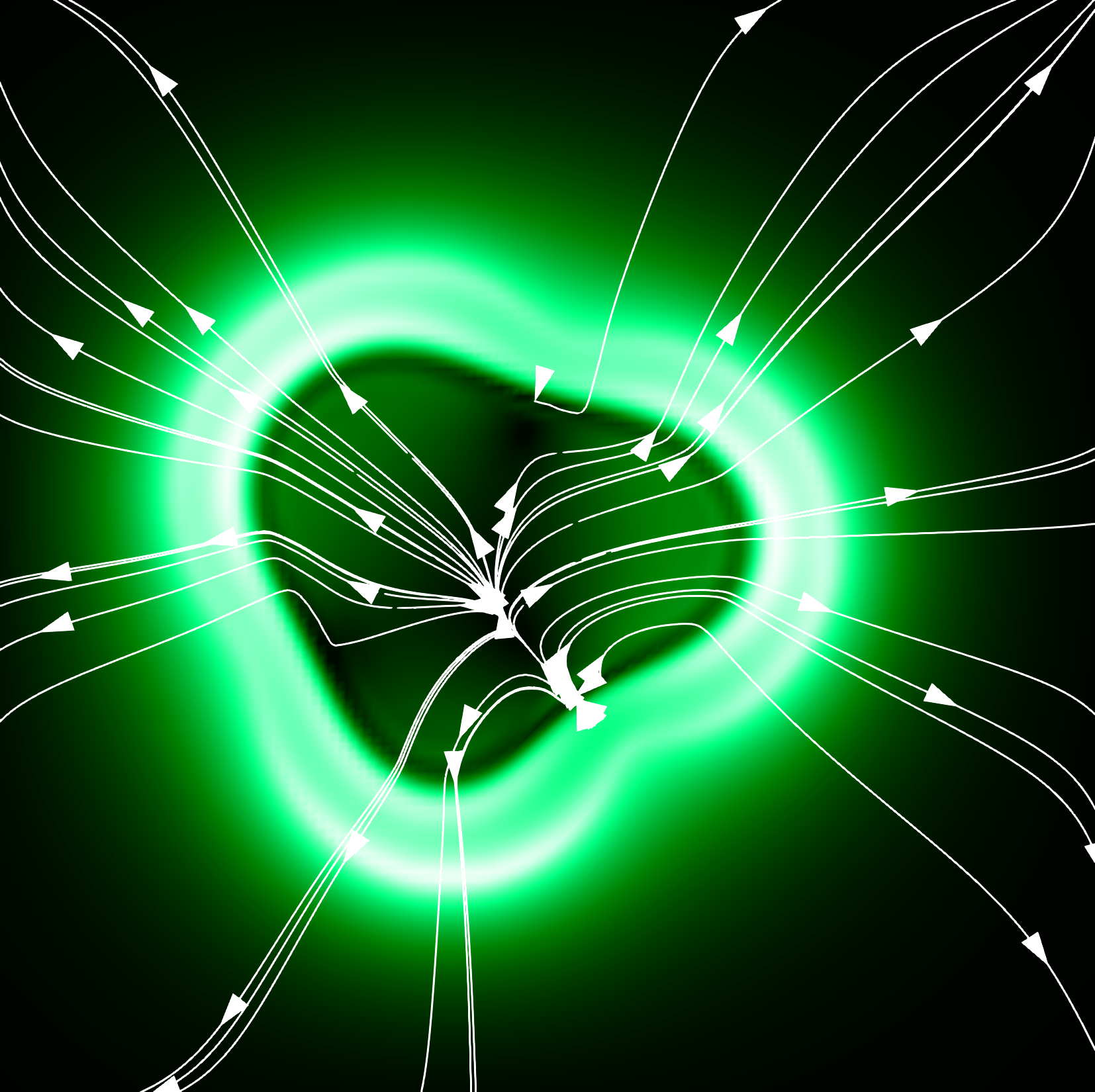}
\end{tabular} \\
\hfill\includegraphics[width=.88\textwidth]{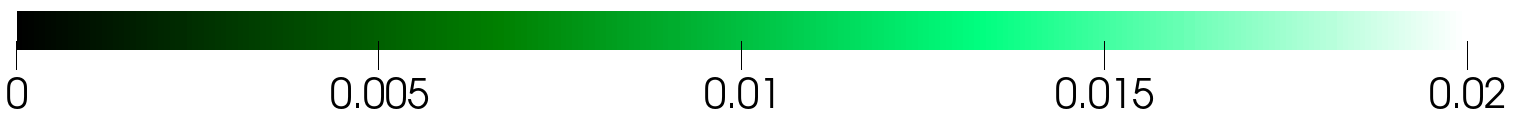}
\caption{\label{Fig:Flow}Evolution of the Darcy velocity $q(t,x)=-\kappa(\phi)\nabla p$ (which is computed in a post-processing step).}
\end{center}
\end{figure}

\cref{Fig:Deformation} provides a visual representation of the deformation profiles observed in the CL and CHB models. These profiles offer valuable information on how different factors influence the deformation patterns of the tumor and surrounding tissue. In the upper row of \cref{Fig:Deformation}, we examine the deformation profile resulting from the CL model. The deformation pattern in this model exhibits a concentric profile, that is, the deformation points outwards from the center of the domain towards the boundary. Further, we want to highlight that the most significant deformation occurs at the interface of the tumor, emphasizing the impact of elasticity on the tumor's mechanical response. Moving to the lower row of \cref{Fig:Deformation}, we consider the deformation profile of the CHB model, where the deformation equation also has coupling effects with the Biot model. The deformation exhibits a connection with the flow dynamics as depicted in \cref{Fig:Flow}, with both the Darcy velocity and the deformation profiles having similar concentrations along the tumor's interface.

\section*{Acknowledgments}
Supported by the state of Upper Austria.

	\bibliography{literature.bib}

\providecommand{\href}[2]{#2}
\providecommand{\arxiv}[1]{\href{http://arxiv.org/abs/#1}{arXiv:#1}}
\providecommand{\url}[1]{\texttt{#1}}
\providecommand{\urlprefix}{URL }
\begin{thebibliography}{10}

\bibitem{abels2024existence}
\newblock H.~Abels, H.~Garcke and J.~Haselb{\"o}ck,
\newblock Existence of weak solutions to a {Cahn--Hilliard--Biot} system,
\newblock \emph{Nonlinear Analysis: Real World Applications}, \textbf{81}
  (2025), 104194.

\bibitem{alnaes2015fenics}
\newblock M.~Aln{\ae}s, J.~Blechta, J.~Hake, A.~Johansson, B.~Kehlet, A.~Logg,
  C.~Richardson, J.~Ring, M.~E. Rognes and G.~N. Wells,
\newblock The {FEniCS} project version 1.5,
\newblock \emph{Archive of Numerical Software}, \textbf{3}.

\bibitem{areias}
\newblock P.~Areias, E.~Samaniego and T.~Rabczuk,
\newblock A staggered approach for the coupling of {C}ahn--{H}illiard type
  diffusion and finite strain elasticity,
\newblock \emph{Computational Mechanics}, \textbf{57} (2016), 339--351.

\bibitem{bociu2016analysis}
\newblock L.~Bociu, G.~Guidoboni, R.~Sacco and J.~T. Webster,
\newblock Analysis of nonlinear poro-elastic and poro-visco-elastic models,
\newblock \emph{Archive for Rational Mechanics and Analysis}, \textbf{222}
  (2016), 1445--1519.

\bibitem{bociu2023mathematical}
\newblock L.~Bociu, B.~Muha and J.~T. Webster,
\newblock Mathematical effects of linear visco-elasticity in quasi-static biot
  models,
\newblock \emph{Journal of Mathematical Analysis and Applications},
  \textbf{527} (2023), 127462.

\bibitem{boyer2012mathematical}
\newblock F.~Boyer and P.~Fabrie,
\newblock \emph{Mathematical Tools for the Study of the Incompressible
  Navier--Stokes Equations and Related Models},
\newblock Springer Science \& Business Media, 2012.

\bibitem{brunk2024}
\newblock A.~Brunk and M.~Fritz,
\newblock Structure-preserving approximation of the {C}ahn--{H}illiard--{B}iot
  system,
\newblock \emph{arXiv preprint arXiv:2407.12349}.

\bibitem{tumorstresscheng}
\newblock G.~Cheng, J.~Tse, R.~Jain and L.~Munn,
\newblock Micro-environmental mechanical stress controls tumor spheroid size
  and morphology by suppressing proliferation and inducing apoptosis in cancer
  cells,
\newblock \emph{PLoS One}, \textbf{4} (2009), e4632.

\bibitem{colli2017asymptotic}
\newblock P.~Colli, G.~Gilardi, E.~Rocca and J.~Sprekels,
\newblock {Asymptotic analyses and error estimates for a Cahn--Hilliard type
  phase field system modelling tumor growth},
\newblock \emph{Discrete \& Continuous Dynamical Systems-S}, \textbf{10}
  (2017), 37.

\bibitem{coussy}
\newblock O.~Coussy,
\newblock \emph{Poromechanics},
\newblock John Wiley \& Sons, 2004.

\bibitem{fritz2020subdiffusive}
\newblock M.~Fritz, C.~Kuttler, M.~Rajendran, L.~Scarabosio and B.~Wohlmuth,
\newblock On a subdiffusive tumour growth model with fractional time
  derivative,
\newblock \emph{IMA Journal of Applied Mathematics}, \textbf{86} (2021),
  688--729.

\bibitem{fritz2019localnonlocal}
\newblock M.~Fritz, E.~Lima, V.~Nikoli{\'c}, J.~T. Oden and B.~Wohlmuth,
\newblock Local and nonlocal phase-field models of tumor growth and invasion
  due to {ECM} degradation,
\newblock \emph{Mathematical Models and Methods in Applied Sciences},
  \textbf{29} (2019), 2433--2468.

\bibitem{fritz2023tumor}
\newblock M.~Fritz,
\newblock Tumor evolution models of phase-field type with nonlocal effects and
  angiogenesis,
\newblock \emph{Bulletin of Mathematical Biology}, \textbf{85} (2023), 44.

\bibitem{fritz2021analysis}
\newblock M.~Fritz, P.~K. Jha, T.~K{\"o}ppl, J.~T. Oden and B.~Wohlmuth,
\newblock Analysis of a new multispecies tumor growth model coupling 3{D}
  phase-fields with a 1{D} vascular network,
\newblock \emph{Nonlinear Analysis: Real World Applications}, \textbf{61}
  (2021), 103331.

\bibitem{fritz2019unsteady}
\newblock M.~Fritz, E.~A. Lima, J.~T. Oden and B.~Wohlmuth,
\newblock {On the unsteady Darcy--Forchheimer--Brinkman equation in local and
  nonlocal tumor growth models},
\newblock \emph{Mathematical Models and Methods in Applied Sciences},
  \textbf{29} (2019), 1691--1731.

\bibitem{fritz2022time}
\newblock M.~Fritz, M.~L. Rajendran and B.~Wohlmuth,
\newblock Time-fractional {C}ahn--{H}illiard equation: {W}ell-posedness,
  degeneracy, and numerical solutions,
\newblock \emph{Computers \& Mathematics with Applications}, \textbf{108}
  (2022), 66--87.

\bibitem{garcketumormechanics}
\newblock H.~Garcke, K.~Lam and A.~Signori,
\newblock On a phase field model of {C}ahn--{H}illiard type for tumour growth
  with mechanical effects,
\newblock \emph{Nonlinear Analysis: Real World Applications}, \textbf{57}
  (2021), 103192.

\bibitem{garcke2005mechanical}
\newblock H.~Garcke,
\newblock Mechanical effects in the {C}ahn--{H}illiard model: {A} review on
  mathematical results,
\newblock \emph{Mathematical Methods and Models in Phase Transitions},
  \textbf{1} (2005), 43--77.

\bibitem{garcke2016global}
\newblock H.~Garcke and K.~F. Lam,
\newblock {Global weak solutions and asymptotic limits of a
  Cahn--Hilliard--Darcy system modelling tumour growth},
\newblock \emph{AIMS Mathematics}, \textbf{1} (2016), 318--360.

\bibitem{garcke2017analysis}
\newblock H.~Garcke and K.~F. Lam,
\newblock Analysis of a {C}ahn--{H}illiard system with non-zero {D}irichlet
  conditions modeling tumor growth with chemotaxis,
\newblock \emph{Discrete and Continuous Dynamical Systems}, \textbf{37} (2017),
  4277--4308.

\bibitem{garcke2023approximation}
\newblock H.~Garcke and D.~Trautwein,
\newblock Approximation and existence of a viscoelastic phase-field model for
  tumour growth in two and three dimensions,
\newblock \emph{Discrete and Continuous Dynamical Systems-S}, \textbf{17}
  (2024), 221--284.

\bibitem{grasselli2023phase}
\newblock M.~Grasselli and A.~Poiatti,
\newblock A phase-field system arising from multiscale modeling of thrombus
  biomechanics in blood vessels: {L}ocal well-posedness in dimension two,
\newblock \emph{Discrete and Continuous Dynamical Systems-S}, \textbf{16}
  (2023), 2364--2425.

\bibitem{tumorstresshelminger}
\newblock G.~Helmlinger, P.~A. Netti, H.~C. Lichtenbeld, R.~J. Melder and R.~K.
  Jain,
\newblock Solid stress inhibits the growth of multicellular tumor spheroids,
\newblock \emph{Nature Biotechnology}, \textbf{15} (1997), 778--783.

\bibitem{lima2016}
\newblock E.~Lima, J.~T. Oden, D.~Hormuth, T.~Yankeelov and R.~Almeida,
\newblock Selection, calibration, and validation of models of tumor growth,
\newblock \emph{Mathematical Models and Methods in Applied Sciences},
  \textbf{26} (2016), 2341--2368.

\bibitem{lima2017}
\newblock E.~Lima, J.~T. Oden, B.~Wohlmuth, A.~Shahmoradi, D.~Hormuth~II,
  T.~Yankeelov, L.~Scarabosio and T.~Horger,
\newblock Selection and validation of predictive models of radiation effects on
  tumor growth based on noninvasive imaging data,
\newblock \emph{Computer Methods in Applied Mechanics and Engineering},
  \textbf{327} (2017), 277--305.

\bibitem{milosevictumor}
\newblock M.~Milosevic, S.~Lunt, E.~Leung, J.~Skliarenko, P.~Shaw, A.~Fyles and
  R.~Hill,
\newblock Interstitial permeability and elasticity in human cervix cancer,
\newblock \emph{Microvascular Research}, \textbf{75} (2008), 381--390.

\bibitem{miranville2019cahn}
\newblock A.~Miranville,
\newblock \emph{The Cahn--Hilliard Equation: Recent Advances and Applications},
\newblock SIAM, 2019.

\bibitem{riethmuller2023well}
\newblock C.~Riethm{\"u}ller, E.~Storvik, J.~W. Both and F.~A. Radu,
\newblock {Well-posedness analysis of the Cahn--Hilliard--Biot model},
\newblock \emph{arXiv preprint arXiv:2310.18231}.

\bibitem{roubivcek2013nonlinear}
\newblock T.~Roub{\'\i}{\v{c}}ek,
\newblock \emph{Nonlinear Partial Differential Equations with Applications},
\newblock Springer Science \& Business Media, 2013.

\bibitem{showalter2000diffusion}
\newblock R.~E. Showalter,
\newblock Diffusion in poro-elastic media,
\newblock \emph{Journal of mathematical analysis and applications},
  \textbf{251} (2000), 310--340.

\bibitem{storvik2022cahn}
\newblock E.~Storvik, J.~W. Both, J.~M. Nordbotten and F.~A. Radu,
\newblock A {C}ahn--{H}illiard--{B}iot system and its generalized gradient flow
  structure,
\newblock \emph{Applied Mathematics Letters}, \textbf{126} (2022), 107799.

\bibitem{storvik2024sequential}
\newblock E.~Storvik, C.~Riethm{\"u}ller, J.~W. Both and F.~A. Radu,
\newblock Sequential solution strategies for the {Cahn--Hilliard--Biot model},
\newblock \emph{arXiv preprint arXiv:2401.13358}.

\bibitem{tumorstressstylianopoulos}
\newblock T.~Stylianopoulos, J.~D. Martin, V.~P. Chauhan, S.~R. Jain,
  B.~Diop-Frimpong, N.~Bardeesy, B.~L. Smith, C.~R. Ferrone, F.~J. Hornicek,
  Y.~Boucher et~al.,
\newblock Causes, consequences, and remedies for growth-induced solid stress in
  murine and human tumors,
\newblock \emph{Proceedings of the National Academy of Sciences}, \textbf{109}
  (2012), 15101--15108.

\end{thebibliography}
	\bibliographystyle{AIMS} 

\medskip
Received xxxx 20xx; revised xxxx 20xx; early access xxxx 20xx.
\medskip

\end{document}